\documentclass[]{amsart}
\usepackage{latexsym,fancyhdr,amssymb,color,amsmath,amsthm,graphicx,listings,comment}
\usepackage[section]{placeins}
\newtheorem{thm}{Theorem} 
\newtheorem{lemma}{Lemma} 
\newtheorem{propo}{Proposition} 
\newtheorem{coro}{Corollary}
\let\paragraph\subsection

\title{On Graphs, Groups and Geometry}
\author{Oliver Knill}
\date{5/27, 2022, updated 8/23/2026}
\address{Department of Mathematics \\ Harvard University \\ Cambridge, MA, 02138 }

\begin{document}
\maketitle

\begin{abstract}
An abstract group $(G,*)$ is natural if there exists a metric structure $(G,d)$ on $G$ such 
that $(G,*)$ is up to abstract group isomorphisms the only group structure on $G$ for 
which all right translations are isometries of $(G,d)$. 
Every connected Lie group $G$ is natural. Disconnected Lie groups can be non-natural, like
$\mathbb{R} \times C_2$ or $Pin^-(2)=Dic(S^1,-1)$. 
There are also natural disconnected groups like $O(n)$. 
Every reflection group of cardinality at most the continuum is natural.
Both the infinite dihedral group $D_{\infty} =  C_2 * C_2$ and
modular group $PSL(2,\mathbb{Z}) \cong C_2 * C_3$ are natural.
The rational numbers $\mathbb{Q}$ are natural. The $p$-adic integers for prime $p$ are 
natural if and only if $p$ is odd. The $p$-adic additive groups $(\mathbb{Q}_p,+)$ of 
the fields $\mathbb{Q}_p$ are all natural. 
An arbitrary non-abelian, finitely generated group $G$ is natural 
if and only if $G$ is not generalized dicyclic. 
A finitely generated abelian group is natural if and only if it is finite 
and its Sylow 2-subgroup is elementary abelian. %$G \cong A \times C_2^r$, where $|A|$ is odd.
Examples of non-natural groups are the cyclic groups $C_{4m}$ for $m \geq 1$, 
the quaternion group $Q_8$, the integers $\mathbb{Z}$, or the doubled 
line $\mathbb{R} \times C_2$.  
Examples of natural groups are the dihedral groups $D_n$, the additive group 
$(\mathbb{R}^n,+)$, the Rubik cube, connected Lie groups and all free groups $F_n$ for $n \geq 2$ 
or all $Pin^(n)$ for $n \geq 3$, the Lorentz or Poincar\'e group.
A principal ingredient is the rigidity theory of Cayley graphs developed by 
Leemann and de la Salle.
\end{abstract}

\section{Natural metric spaces, groups and graphs}

\paragraph{}
Define a {\bf metric space} $(G,d)$ to be {\bf natural} if it carries up to 
abstract group isomorphisms a unique group structure $(G,*,1)$ such that for every $g \in G$, 
the {\bf right group translations} $R_g: x \to x*g$ 
is an isometry of the metric space $(G,d)$. A metric space is {\bf strongly natural}, if
it carries up to isomorphism a unique group structure for which both right translations $R_g$ and
left group translations $L_g: x \to g*x$ are isometries. In the strongly natural case,
also the {\bf inversion} $I: x \to i(x)=x^{-1}$ is an {\bf isometry} \footnote{See appendix}
We do not discuss the strong natural case here and focus on naturality.

\paragraph{}
A {\bf group} written additively as $(G,+,0)$ or multiplicatively as $(G,*,1)$ is called 
{\bf natural} if it emerges as the group structure from a natural metric space $(G,d)$. 
A group is {\bf strongly natural}, if it emerges as the group structure from a strongly 
natural metric space, again up to abstract group isomorphisms.
For a natural group, we only need to give the metric space to assure that there is a multiplication table 
which is unique up to an abstract group isomorphism. The algebra is forced by the geometry, 
space determines time. A natural group $X$ has a {\bf regular} (transitive and free) 
action on $X$ by isometries, but $(S,*,d)$ is not necessarily a topological group.
\footnote{The strongly natural assumption gives different results \cite{Gismatullin}. See the appendix.} 

\paragraph{}
A connected {\bf simple graph} $X=(V,E)$ is declared to be {\bf natural} if it is a natural 
metric space $(X,d)$ with the geodesic metric $d$. We will also look at infinite graphs
or {\bf weighted graphs} in which the edges $E$ are assigned a length such that the
geodesic distance produces a metric space $(V,d)$. Every
finite metric space can be seen as a weighted simple complete graph $(V,E,d)$.
Graphs are important in the context of {\bf graphical regular representation}. The later is a Cayley
graph $\Gamma$ for which ${\rm Aut}(\Gamma)$ is $R(G)$ the group of right translations.
\footnote{Edges of this graph are $(g,sg)$ so that {\bf right translations} are graph automorphisms.}
The natural graph $\Gamma$ defines a geometry on which the group is at home, hence the title
of this paper. As for a general reference for this paper, \cite{Godsil} cover
groups, transitive graphs, Cayley graphs and actions by automorphisms. The papers
\cite{LeemannDeLaSalle2021,LeemannDeLaSalle2022} are crucial for this update. 

\paragraph{}
In order to see that a metric space $(X,d)$ is natural, one has to be able to construct 
the group structure from it and verify that it is unique up to abstract group isomorphism.
A metric space $(X,d)$ is natural if and only if all regular subgroups of ${\rm Isom}(X,d)$
are mutually isomorphic and at least one a regular subgroup exists. 
Related is Sabidussi's characterization of Cayley graphs \cite{Sabidussi1964}:
a graph is Cayley if and only if its automorphism group contains a regular subgroup. 
In order to see that a group is natural, one has to find 
a natural metric space that generates the group. To show that a metric space $(X,d)$ is not natural, 
one has either to establish that no metric-compatible group structure
is possible on $X$, or then that at least two non-isomorphic competing regular group 
structures preserve the metric by right translation.

\paragraph{} 
Here is the general result for finite groups: if $S=S^{-1}$ is a generator set of a group $G$, 
denote by $\Gamma(G,S)$ the corresponding {\bf Cayley graph}. One can either stick with the geodesic 
metric (which is what graphical regular representations do) or change to a more general metric space $(G,d)$. 
By definition, any natural finite group $(G,+)$ comes from a natural metric 
space $(G,d)$ and so from a complete weighted Cayley graph $\Gamma$ for $G$
with $d(g,sg)=d(1,s)$, where $S=G \setminus \{1\}$.
Again, the assumption is that $(g,sg)$ are the edges of $\Gamma$, forcing 
{\bf right translations} to be graph automorphisms.
If such a metric finite space $(G,d)$ forces the group structure, it is natural.
New in this update is the realization (prompted from \cite{LeemannDeLaSalle2021,LeemannDeLaSalle2022})
that the Cayley group idea is much more prevalent and works also in groups which are countable
or even have the cardinality of the continuum. Natural metric spaces can be constructed from
weighted Cayley graphs with a continuum of points and a continuum of edges. 

\paragraph{}
A prototype example in the class of natural groups is the dihedral group
$D_n = C_n \rtimes_{-1} C_2$ for $n \geq 2$ or the infinite 
dihedral group $D_{\infty}$. In the finite case, represent the points in $D_n$
as $\{1,2, \dots, 2n=0\}$ and define the geodesic metric $d$ using the edge lengths
$d(x,x+1)=2$ for even $x$ and $d(x,x+1)=3$ for odd $x$. The group $D_n$ is generated
by $r(x)=x+2 \; {\rm mod} \; 2n$ and $s(x)=1-x \; {\rm mod} 2n$. Both generators act as 
isometries via right translation. With alternating edge lengths $2$,$3$, an isometry must 
preserve the two edge types. (Pairs at distance $2$ or $3$ must be adjacent because every path using two edges 
has length at least 5). Hence the weighted graph can be recovered from its vertex metric. 
The color-preserving automorphism group of this $2n$-cycle is exactly the $2n$-element 
dihedral group generated by $r,s$. In other words, the color-preserving stabilizer group of $1$
is the trivial group. The entire isometry group already has order $2n$ 
and acts regularly. A more general argument is given below for groups that are generated by 
involutions, an argument which also works also for a countable set of involutions
or even involutive groups with the cardinality of the continuum. 
The disconnected Lie group $O(n)$ for example is quickly seen to be natural because it 
is a reflection group by Cartan-Dieudonn\'e's theorem. 
\footnote{
In the 2022 draft of this paper we used products like the direct product, the semi-direct product
for groups, and Shannon products and the zig-zag product \cite{AlonLubotzkyWigderson} for graphs. 
But naturality does not survive in general the usual product as $C_6=C_2 \times C_3$ shows. 
It worked in the strongly natural case as Jakub Gismatullin pointed out. I'm grateful to
Jakub to alert  me about my original confusions with respect to
left and right invariance. Naturality has emerged now to be more accessible.
The reason is orientation rigidity of Leemann and de la Salle 
\cite{LeemannDeLaSalle2021,LeemannDeLaSalle2022}.}

\section{Finite groups}

\paragraph{}
A {\bf generalized dicyclic group} ${\rm Dic}(A,y=x^2)$ is a non-abelian group with an
abelian subgroup $A$ of index $2$ and an element $x \notin A$ such that $y^2=x^4 = 1$ and
$x^2 \neq 1$ and $x g x^{-1} = g^{-1}$ for all $g \in A$. See also \cite{MorrisSpigaVerret}. 
\footnote{Thus, $G$ is a non-abelian group that has an
abelian index-two subgroup $A$, and an element $x \notin A$ of order 
$4$ acting on $A$ as an involution. The literature can also write 
$x^4=1$ with the implicit assumption that $x$ has order $4$ but with $x^2 \neq 1$. }

\paragraph{}
A good example to have in mind is the {\bf quaternion group}
$Q_8={\rm Dic}(C_4,-1=j^2)$, where the quaternion unit $j$ is the additional 
generator $j$ outside $C_4$ generated by $i$. 
Writing $C_4=\langle 1,i \rangle = \{ 1,i,-1,-i\}$, the construction gives the 8-element sub-group
$$ Q_8=\{1,-1,i,-i,j,-j,k,-k \} $$ of the unit quaternions $SU(2)$. 
The group $Q_8$ is not natural. An other example which is the 
{\bf generalized} dicyclic group ${\rm Dic}(S_1,-1=x^2)$ 
which is also known as ${\rm Pin}^-(2)$ and which is not natural. 
Its twin ${\rm Pin}^+(2)=O(2)=S^2 \rtimes C_2$ is {\bf generalized dihedral} group 
and natural. It is the group generated by reflections at lines in $\mathbb{R}^2$. 
The naturality condition somehow revives a bit the idea of {\bf "Spiegelungs geometry"} \cite{Bachmann1959}.
\footnote{My geometry teacher Max Jeger at ETH heavily 
focused on reflection geometry \cite{JegerTransformationGeometry}.}

\paragraph{}
Here are two elementary properties which hold for a generalized dicyclic 
group $G={\rm Dic}(A,y)$ with $y=x^2$: \\ 
a) Every $g \notin A$ satisfies $g^2=y$. \\
Proof: An element $g \notin A$ is of the form $g=xa$ with $a \in A$. Now 
$g^2=(xa)^2 = (xa) (xa) = (a^{-1} x) (x a) = a^{-1} y a = y a^{-1} a = y$. \\
b) $y=x^2$ is central: $x y x^{-1} = y^{-1} = y$.  \\
Proof: by definition, $x y x^{-1} = y^{-1}$. But since $y^2=1$, $y^{-1} = y$.

\paragraph{}
The following result is stated for finite groups. We see below that it holds 
also for finitely generated groups. It already shows the main idea 
to use a weighted Cayley graph.

\begin{thm}[A corollary from Leemann-DeLaSalle]
A finite non-abelian group $G$ is natural if and only if it 
is not generalized dicyclic.
\end{thm}

\begin{proof}
(i) Assume first that $G$ is not generalized dicyclic. We show that it is natural. To do so,
construct a metric space $(G,d)$. For every unordered pair $\{ g,g^{-1} \}$ with $g \neq 1$, 
assign a distinct number $1<\lambda_{\{ g,g^{-1} \}}<2$ and define 
$d(x,y) = \lambda_{\{ y x^{-1}, x y^{-1} \}}, x \neq y$ and $d(x,x)=0$ along the diagonal. It is a metric because
all nonzero distances are in the interval $(1,2)$, so that for any three distinct points $x,y,z \in G$,
the triangle inequality $d(x,z)<2<d(x,y)+d(y,z)$ holds. An isometry $T$ fixing $1$ must also fix every distance class. 
Proof: Since every inverse-pair $\{g,g^{-1}\}$ has its own distinct metric value, $T$ 
is a color-preserving automorphism of the complete Cayley graph.
By the orientation-rigidity theorem of Leemann and de la Salle applied to the complete Cayley graph
any nontrivial such automorphism fixing 1 can occur only in the abelian or dicyclic cases.
So, for a non-abelian, non-generalized-dicyclic G, $T(x)=x$ so that $Isom(G,d)=R(G)$. 
{\bf Orientation rigidity} so forces it to be the identity, leading to ${\rm Isom}(G,d)=R(G)$.
Here comes \cite{LeemannDeLaSalle2022} in. They assure that if $G$ is non-abelian,
unless $G$ is generalized dicyclic, there is no non-trivial 
color preserving automorphism that fixes the identity. As ${\rm Isom}(G,d) = R(G)$, the metric is natural. \\ 
(ii) Now assume that $G={\rm Dic}(A,y)$ is generalized dicyclic. We show that it is not-natural.
$G$ admits the involutive isometry $i(a)=a, a \in A, i(g) =g^{-1}$  for $g \notin A$. It will lead us
to a competitive group. Write $l(g)=d(1,g)$. For a right invariant metric, $d(u,v)=l(v u^{-1})$ and 
$l(g)=l(g^{-1})$. We need to verify that the involution $i$ as defined is an isometry:
\begin{itemize}
\item If $u,v \in A$, this is obvious as $i$ fixes $u,v$.
\item If $u,v \notin A$, then $d(i(u),i(v)) = d(u^{-1},v^{-1}) = l(v^{-1} u$. But since $u^{-1} = yu$ and $v^{-1} = y v$,
this gives $v^{-1} u = y v u$ which is by centrality $v y u = v u^{-1}$. Hence $d(i(u),i(v))=l(v u^{-1}) = d(u,v)$. 
\item In the mixed case $a \in A, g \in G \setminus A$, we use that also $g^{-1}$ acts on $A$ by inversion 
$g^{-1} a^{-1} = a g^{-1}$. Now $d(i(a),i(g))= d(a,g^{-1}) = l(g^{-1} a^{-1})=l(a g^{-1})=d(g,a)=d(a,g)$.
\end{itemize}
The element $x$ is outside $A$. The group generated by both $\tau=R_x \circ i$ and $\{ R_a, a \in A \}$
produces a {\bf generalized dihedral competitor} that 
is not isomorphic to $G$. Note that $\tau^2=1$ so that $\tau R_a \tau = (R_x \circ i) R_a (R_x \circ i) = R_{a^{-1}}$. 
This is a regular copy of $A \rtimes_{-1} C_2$. 
To see that $H$ acts transitively, note that $|H|=2|A|=|G|$ and since $R(A)$ reaches $A$, while $\tau(1)=x$ reaches
the other co-set, it produces a competitor group $H$ that is not isomorphic to the generalized 
dicyclic group $G$. All elements in $H$ that are in a nontrivial co-set of $A$ are involutions, while in 
the dicyclic group $G$, they have order $4$. 
More precisely, $H$ has all $|A|$ elements outside $A$ as involutions, 
whereas $G$ has none outside $A$. Thus $H$ has exactly $|A|$ more involutions 
than $G$, showing that $H,G$ cannot be isomorphic. This finishes the proof that 
generalized dicyclic groups are non-natural. 
\end{proof} 

\paragraph{}
The property of being natural goes beyond graphical regular representations in general.
Groups like $A_4, B(2,3), P(1), Q_8 \times C_3, Q_8 \times C_4$ are natural,
even so they fail the unweighted graphical regular representation filter.
The ordinary dicyclic groups can be implemented in the Lie group $SU(2)$
as $a={\rm Diag}(e^{\pi \frac{i}{n}},e^{-\pi \frac{i}{n}})$ 
and $x=\left[ \begin{array}{cc} 0 & -1 \\ 1 & 0 \end{array} \right]$.

\begin{comment}
F        := FreeGroup("a","x");; a := F.1;;  x := F.2;;
C4       := F / [ a^4, x];;;                    #  Dic1
Q8       := F / [ a^4, x^2*a^-2, x*a*x^-1*a];;  #  Dic2
Dic3     := F / [ a^6, x^2*a^-3, x*a*x^-1*a];;  #  Dic3
Q16      := F / [ a^8, x^2*a^-4, x*a*x^-1*a];;  #  Dic4
C4semiC4 := F / [ a^4, x^4, x*a*x^-1*a ];; 
F        := FreeGroup("a","b","x");; a := F.1;; b := F.2;; x := F.3;;
Q8xC2    := F / [ a^4,b^2,Comm(a,b), x^2*a^-2, x*a*x^-1*a, Comm(x,b) ];; 
L        := [Q8, Dic3, Q16, Q8xC2, C4semiC4]; List(L, Size); List(L, IdGroup);
\end{comment} 

\paragraph{}
Together with the abelian cases, we know from all finite groups whether they are 
natural or not. Also the next result will be extended later to finitely generated 
abelian groups. We leave it because finite groups are interesting per se.

\begin{thm}
A finite abelian group $G$ is natural if and only if it 
either finite of odd order $|G|$ or $G = C_2^r$.
\end{thm}
\begin{proof}
Because $(G,+)$ is abelian, we the group operation additively and use the identity zero $0$.  \\
a) {\bf Odd order abelian groups are natural}: For every (non-oriented) pair $\Lambda_a=\{ a,-a \}$ with $a \neq 0$
choose different weights $1<\lambda_{\{a,-a\}}<2$ and again define the metric $d(x,y) = \lambda_{\{ x-y, y-x \} }
=\lambda_{\{ y-x, x-y \}}$ and $d(x,x)=0$. If $T$ is an isometry, fixing $0$, then $T(a) \in \{a,-a\}$. 
Since $|G|$ is odd, an isometry fixing $0$ and preserving every pair $\{a,-a\}$ must be globally either the 
identity or inversion. All signs must be the same so that $T(x)=x$ or $T(x)=-x$. 
Proof: if $T(a)=a$ and $T(b)=-b$ for nonzero elements in $G$, preservation of the color of $a-b$ gives
$a+b= \pm (a-b)$. The $+$ case, we get $2b=0$; the $-$ case gives $2a=0$. 
Both are impossible in an group that has odd order. Hence mixed signs cannot occur and
${\rm Isom}(G,d) = G \rtimes C_2$. Given any alternative compatible group $H$ that
is a subgroup of ${\rm Isom}(G,d)$. We show that it agrees with $G$. Since every element outside 
the translation group is an involution we have $H \subset G$, the alternative regular group $H$ has the 
same odd order as $G$, and hence contains no involutions. This implies $H=G$. \\
b) {\bf Elementary abelian $2$-groups are natural:} (This will be generalized vastly even to infinite groups
generated by involutions). In the case $G=C_2^r$ it is more pictorial. choose generators so that the
Cayley graph is the $r$-dimensional hypercube. Fix the lengths along different coordinate axes,
so that $1<\lambda_j <2$ for $1 \leq j \leq r$ are all different. Any isometry fixing $0$ must be the
identity. Proof: because every $\lambda_j$ is in the open interval $(1,2)$,
the points at distance $<2$ from $0$ are the 
coordinate neighbors; their distinct distances force each of them to be fixed. 
A metric isometry is then a hypercube automorphism fixing every coordinate direction, hence fixes every vertex.
This shows ${\rm Isom}(G,d)$ contains the group $R(G)$ of right translations of $G$, so that $G$ is natural.\\
c) {\bf Every abelian group of even order that is not of a Boolean group $C_2^r$ is not natural.}
Given any right invariant metric $d(x,y)$. It is of the form $d(x,y) = l(x-y)$ with $l(a)=l(-a)$ so that 
inversion is an isometry. As $|G|$ is even, there is subgroup $H$ of index $[G:H]=2$ by the fundamental theorem
of abelian groups. For $c \notin H$, look at the isometry $f_c(x)=c-x$. It satisfies 
$f_c^2=1, f_c R_h f_c = R_{-h}$.  It generates, together with translations in $H$,
the group $K=H \rtimes_{-1} C_2$. This group $K$ acts regularly on $G$ and is not isomorphic to $G$.\\
(i) proof that $K$ acts regularly: every element in the non-translation coset acts as
$x \to c-x +h$ for some $h \in H$. A fixed point $x$ would satisfy 
$2x=c+h$. But as $2x \in H$ since $G/H=C_2$, this would give $c+h \in H$ and so $c \in H$. 
This is a contradiction because $c$ was assumed to be outside $H$. \\
(ii) proof that $K$ is not isomorphic to $G$: distinguish two cases:
if $H$ contains an element of order larger than $2$, then $K$ is nonabelian while $G$ is abelian. 
If $H$ has exponent $2$, then $K$ is elementary abelian, whereas $G$ was assumed not to 
be elementary abelian. 
\end{proof} 

\paragraph{}
The just covered finite case is part of a more general theme.
Any non-abelian, non-finite group that is not generalized dicyclic has a graphical regular representation and
so is natural.  This covers Baumslag-Solitar groups, hyperbolic groups that are not generalized dicyclic.
{\bf Every infinite finitely generated abelian group is non-natural.}
Examples are $\mathbb{Z}^r \times F$, where $F$ is a finite abelian group.
{\bf Every infinite, finitely generated generalized dicyclic group is non-natural.}
This means that every infinite finitely presented group is already classified.
Together with the classification of all finite groups,
Naturality for all finitely generated groups is completely classified.

\section{Reflection groups}

\paragraph{}
Reflection groups $G$ are natural in great generality. All is needed is a bound on the {\bf cardinality} of $G$. 
We can not go beyond the continuum yet. We just need  part of graphical regular 
representation theory of \cite{LeemannDeLaSalle2022}. Having reflections as generators allows to 
define a metric space that is a Cayley graph. This works even if the generators are countable
or if the group $G$ has the cardinality of the continuum. The theory only stops working, if one goes beyond the 
continuum. The reason for the cardinality bound is that a metric space $(X,d)$ uses a function 
$d: X^2 \to [0,\infty)$ and that the image $[0,\infty)$ has the cardinality of the continuum. 
A group with cardinality larger than the continuum can not be forced by a metric Cayley graph
with $G$ as vertices. 

\paragraph{}
Let $G$ be a finite group generated by a set $S$ of reflections. We write $G=\langle S \rangle$. 
An example is the symmetric group $S_n$, the set of permutations of $\{1,2,\dots, n\}$.
The group $S_n$ is well known to be generated by transpositions.
This fact already shows that it is natural. While every finite group is a subgroup of a permutation group,
not all groups are reflection groups. An example is $Pin^-(1) \cong C_4$ which is the prototype of 
a non-natural group. It can not be generated by involutions as every non-identity element has order $4$. 
On the other hand, the {\bf Klein 4-group} $Pin^+(1) \cong  C_2 \times C_2$ is natural as it is generated by involutions. 

\paragraph{}
An example of a group that is generated by an infinite but countable set of involutions is 
$G=S_{fin}(\mathbb{N})= \bigcup_n S_n$, the {\bf group of permutations of $\mathbb{N}$ with finite support}.
It is countable, locally finite, non-abelian. The group is not finitely generated,
but it is generated by a countable set $S$ of transpositions. Every element in $G$ can 
be written as a product of {\bf finitely many} transpositions.

\paragraph{}
The orthogonal group $O(n)$ is a disconnected Lie group. We will see below that all connected
Lie groups are natural so that $O(n)$ is an interesting test case of a disconnected group which needs
special discussion. The group $O(n)$ is generated by reflections. It turns out that also the continuum 
is not an obstacle: we can still build a Cayley graph, even so the vertex set of this graph 
has now the cardinality of the continuum. Examples are the dihedral 
$Pin^+(2)=O(2) = S^1 \rtimes C_2$ is generated by involutions and so will be natural, while its twisted
dicyclic sibling $Pin^-(2)={\rm Dic}(S^1,-1)$ is not generated by involutions and is not natural.
For $n \geq 3$, all $Pin^{\pm}(n)$ are generated by involutions and are natural.

\paragraph{}
Every group $G$ generated by at most continuum $S$ of involutions is natural.  Lets make this 
more precise:

\paragraph{}
We say that a subset $S \subset G$ of a group $(G,*)$ {\bf generates $G$} if every 
$g \in G$ can be written as a {\bf finite} product $g=s_1*s_2* \cdots s_n$, with $s_k \in S$. 
A group $G$ is called a {\bf reflection group}, if it is generated by a subset $S$ of $G$
for which every element is a reflection. While the Cayley graph generated for such a group
has the cardinality of the continuum, it is still possible to connect any two points
$x,y$ by a finite path within the graph. 

\begin{thm}
A reflection group $G$ with the cardinality not exceeding the continuum is natural. 
\end{thm}

\begin{proof}
Let $S$ be the set of involutions generating the group. We have $S=S^{-1}$ because $S$ consists
of reflections. We just need to show that the group $G$ admits a right-invariant metric 
structure $(G,d)$ for which $R(G) = Iso(G,d)$, assures thus that $G$ is natural. To do so, we
build a graph $\Gamma=(V,E)$ in which $V=G$ are the vertices and where two 
points $g,h$ are connected, if there exists $s$ such that $s g = h$. Also the edge set is
at most of continuum size. The fact that edges of $\Gamma$ are $(g,sg)$ assures that
{\bf right translations} are graph automorphisms. \\
The $S$ has maximally the cardinality of the of the continuum (as a subset of $G$), 
there is a bijection of $S$ with a subset of the open interval $(1,2)$. If $s \to \lambda_s$ is 
this bijection, we have a weight for each generator $s \in S$. We now use this to define 
the usual {\bf path metric} using lengths $d(g,sg) = \lambda_s$. The geodesic distance is 
then a metric: symmetry and the triangle inequality work because of the assumption $\lambda_s \in (1,2)$.
Since every $g \in G$ can be written as a product of finitely many $s \in S$, we know
that the distance between any two points $d(g,h)$ is finite.  \\
This path metric is right invariant $d(ga,ha)=d(g,h)$. And here is the key:
every generator $s \in S$ is metrically recognizable because $\{ x, d(g,x) = \lambda_s = \{ s g\} \}$.
The fact that $s$ was an involution assures that there was only one neighbor and not two $\{ s g,s^{-1} g\}$.
If an isometry $T$ of $(\Gamma,d)$ fixes a vertex $g$, it must be the identity: the reason is that if $Tg=g$ then 
$T(s g)=s g$ and by induction on word length $T(g)=g$ for all $g \in G$. This shows that the 
stabilizer group of $Iso(G,d)$ is trivial.  As right translation are always transitive, the action
is regular: For any $T \in {\rm Iso}(G,d)$, choose the unique right translation $R_a$ with
$R_a(1)=T(1)$. Then $R_a^{-1}T$ fixes $1$, hence is the identity. Therefore $T=R_a$.
We have shown ${\rm Iso}(G,d)=R(G)$. 
\end{proof}

\paragraph{}
It follows that the group $S_{fin}(\mathbb{N})$ is natural as each of its blocks $S_n$ is generated
by involutions. Indeed $S_{fin}(\mathbb{N})$ has the generator set $S=\{ (12),(23),(34), \dots \}$. 
It does not stop here. And now we have seen that also countable groups like $\oplus_{x \in \mathbb{R}} C_2$ is natural. 
This is the group of finite subsets of $\mathbb{R}$ for which the addition is the symmetric difference.
An other example is $O(n)$ which by the Cartan-Dieudonn\'e theorem assures that every group element is  
a product of reflections. We have now a proof that that the disconnected Lie group $O(n)$ is natural.
An other example is the {\bf Euclidean group} $R^n \rtimes O(n)$ which for $n \geq 3$ is involution generated 
because it is reflection generated. 

\begin{thm}
Every reflection generated Lie group is natural. 
So, every Lie group containing a non-central involution is natural.
\end{thm}

\paragraph{}
It is known which Lie groups are generated by involutions. They are the Lie groups that 
contain a non-central involution $s$. Proof: 
If $s^2=1$ and $s \notin Z(G)$, there is a nontrivial conjugacy class $C=\{ g s g^{-1}, g \in G \}$
that generates a normal subgroup $H=\langle C \rangle$. 
Because $s$ is non-central, $C$ has positive dimension. The tangent directions of its 
conjugates generate a nonzero ideal of the simple Lie algebra $\mathcal{G}$, 
hence all of g. Finitely many such tangent spaces therefore span g, so some finite product
$C C \cdots C$ contains an open set. Hence H is an open subgroup of the connected group G, and consequently
$H=G$. \footnote{My college geometry teacher Max Jeger would love this as he taught geometry with focus on 
"Spiegelungsgeometrie" \cite{JegerTransformationGeometry,Bachmann1959}. }
An encyclopedic check for Groups that are generated by involutions lists
$SO(n), n\geq 3$ or $SU(n), n \geq 2$ or $Spin(n)$ for $n\geq 5$
or $PSL(2,\mathbb{R})$ ,$SL(n,\mathbb{C})$, compact $G_2,F_4,E_6,E_7,E_8$. The group $SU(2)$ is interesting
as its only involutions are $1,-1$ which are central. But in $SU(2)/C_2=SO(3)$, the 
images of elements of order $4$ become involutions. The rotation group $SO(3)$ is generated by reflections.
Similarly instructive is $SL(2,\mathbb{R})$ which has only involutions $\pm 1$ but $PSL(2,\mathbb{R})$ has
plenty of non-central involutions. If $G$ is generated by involutions, then its abelianization $G/[G,G]$ must
have exponent $2$. For its Lie algebras, this means $\mathcal{g}=[\mathcal{g},\mathcal{g}]$ implying 
that the Lie algebra must be perfect. This excludes $\mathbb{R}^n, U(n),\mathbb{T}^n, GL^+(n,\mathbb{R})$.

\paragraph{}
This provokes the question whether a positive dimensional disconnected Lie group that has no nontrivial involutions
can be natural still. The following example sheds also light on finite groups like why $C_6 = C_2 \times C_3$ was
non-natural, as some reflections $x \to a-x$ become fixed point free allowing
a competing regular group action. This is not possible for $\mathbb{R} \times C_3$, where we have naturality:
{\bf If $F$ is a finite abelian group of odd order $|F| \geq 3$, then $\mathbb{R}^n \times F$ is natural.}
Call a group $G$ {\bf uniquely 2-divisible} if $x \to 2x$ is a bijection of $G$. An obvious example is $\mathbb{R}^n$.

\begin{thm}
Every uniquely 2-divisible abelian group $G$ of cardinality at 
most the continuum is natural.
\end{thm}

\begin{proof} 
The proof is the same.
Choose distinct numbers $1<\lambda_{ \{a,-a\}} <2$ for unordered pairs $a,-a$ which is possible if the group
has cardinality at most the continuum. 
Define the metric $d(x,y)=\lambda_{ \{y-x,x-y\}}$ for $x \neq y$ and $d(x,y)=0$ if $x=y$
as for three distinct points, $d(x,z)<2<d(x,y)+d(y,z)$. All right translations are isometries. 
If an isometry $T$ fixes $0$, then the distinct distance labels imply $T(a) \in \{a,-a\}$ for all $a$.
The signs cannot be chosen independently: suppose that $T(a)=a, T(b)=b$ for two nonzero elements $a \neq b$.
Then preservation of $d(a,b)$ implies $\{b-a,a-b\}=\{b+a,-b-a\}$ which either
implies $2a=0$, or $2b=0$. But that is impossible as $G$ was assumed to have nonzero 2-torsion. Every
isometry fixing 0 is either $T(x)=x$ or or $T(x)=-x$ for all $x$. We see that
${\rm Iso}(G,d)=A \rtimes_{-1} C_2$.
Now comes the crucial naturalness argument. Every isometry outside the translation subgroup has the form
$x \to a-x$. But because multiplication by 2 is surjective, this isometry has the fixed point
$x= 2 a$. Therefore no such reflection can occur in a regular subgroup of the isometry group.
Any alternative group structure on the metric space whose right translations are isometries 
also produces a regular subgroup $H$ of $Iso(G,d)$.
Since a regular action is free, $H$ cannot contain any of the maps $x \to a-x$. Hence $H$ is a subgroup of $G$.
But $H$ must also act transitively. A subgroup of the translation group $G$ that acts transitively 
on $G$ only if it is all of $G$. Thus $H=G$ and $G$ is natural.
\end{proof}

\paragraph{}
The group $(B,+)=(2^{\mathbb{R}},\Delta) \cong \oplus_{2^c} C_2$ of all subsets of 
$\mathbb{R}$ with symmetric difference is an elementary abelian 2-group
for which every non-zero element is an involution. It is still natural
because of a result of Imrich \cite{Imrich1970}: for every infinite cardinal $\kappa$, 
there exists a graph on $\kappa$ vertices whose full automorphism group is transitive and abelian.
Every transitive abelian full automorphism group must act regularly 
and must be elementary abelian 2 so that every infinite elementary abelian 2-group 
has a graphical regular representation. The edge length labeling method is not needed. 

\paragraph{}
It appears currently to be difficult to come up with an example of a group where
we can {\bf not} decide about naturality. Also large cardinalities seem not to 
matter: Let $\kappa$ be any cardinality larger than the continuum. Define the Abelian
group $A_{\kappa} = \oplus_{\alpha < \kappa} C_4$. If $c_0$ is one of the basis elements
then $y=2e_0$ is an involution and $G_{\kappa} = {\rm Dic}(A_{\kappa},y=x^2)$ has the 
cardinality of $\kappa$. While $\oplus_{\kappa} C_2$ is natural for any cardinality, 

\begin{thm}
If $\kappa$ has the cardinality larger than the continuum, 
the generalized dicyclic group $G_{\kappa}$ is still non-natural. 
\end{thm}

\begin{proof}
We show that there is a second competing dihedral structure $H={\rm Dih}(A_{\kappa})$ on $G$
which is not isomorphic. This still works, even so that the cardinality is arbitrary large. \\
(i) As in any dicyclic case, there is a canonical automorphism $i(a)=a$ for $a \in A$ 
and $i(ax)=(ax)^{-1}$ for $a \in A$. If $d$ is a right-invariant metric on $G_{\kappa}$, define $l(g)=d(g,1)$.
Then $d(g,h)=l(gh^{-1})$ and $l(g)=l(g^{-1})$. Since $i(g) \in \{g,g^{-1} \}$ 
$l(i(g))=l(g))$ for all $g$ we have $d(i(g),i(h)) = l(i(gh^{-1}) =l(g h^{-1})=d(g,h)$ 
so that $i$ is an isometry for every right invariant metric on $G_{\kappa}$. \\
(ii) There is a second isometry $T(g) = i(g) x$. It still satisfies $T^2=(g) = i(i(g)x) x = g i(x) x = g x^{-1} x=g$
for all $g \in G$ and $T(a=ax=T(ax)=a$ so that $T$ swaps $A$ with $Ax$ and can not have a fixed point.
If $R_a$ is right translation by $a \in A$, then $T R_a T = R_{a^{-1}}$ so that $H=\langle \{ R_a, a \in A\},T \rangle$ is
the generalized dihedral group $A \rtimes_{-1} C_2=Dih(A)$. 
$H$ acts regularly on the same underlying set $G$ because it is transitive and has trivial stabilizers
We see that every right-invariant metric on $G_{\kappa}$
admits two regular groups of isometries, $Dic(A,y)$ and $Dih(A)$. \\
(iii) Also as seen before, the two groups are not isometric:
In $Dic(A,y)$ every element outside $A$ has square $y \neq 1$ and every involution of $Dic(A,y)$ is central.
In the competing generalized dihedral group $A \rtimes_{-1} C_2$,
the element $T$ is an involution and is non-central if $A$ has an element with $a \neq a^{-1}$. 
$Dic(A_{\kappa},y)$ is not isomorphic to $Dih(A_{\kappa}$. No no metric can make the dicyclic group natural.
\end{proof} 

\paragraph{}
This reminds of the very start of this research where we have seen that $D_4$ is natural but $Q_8$ the
quaternion group is not natural. That example was the smallest instance of a now
{\bf general dihedral-dicyclic dichotomy}. In $D_4=Dih(C_4)$ there are many non-central involutions
so that it is natural. The quaternion group $Q_8= Dic(C_4,a^2)$ with $a^2 = -1$
has all involutions central and so is not natural. The same two groups reappear in the proof of 
non-naturalness of every generalized dicyclic group. If $G=Dic(A,y)$, every right invariant metric on $G$
also admits a competing regular dihedral group $Dih(A) = A \rtimes_{-1} C_2$. This argument survived for general $A$, 
whether $A$ is finite, countable, continuum sized or arbitrarily large. The small example of $D_4$ versus $Q_8$ 
already displayed the entire mechanism. 

\paragraph{}
At the moment it appears even to be difficult to come up with groups for which naturality can not
be decided. As for cardinality, the Leemann-de la Salle orientation rigidity theorem works for arbitrary groups so 
that {\bf every non-abelian, not- generalied dicyclic group $G$ with cardinality of the continuum or less is natural.}
This includes for example the rational Heisenberg group which is not finitely generated. For non-abelian
groups, unexplored territory therefore goes beyond the continuum. A concrete example needs to be huge. 
Let $F$ be a field of characteristic 3 and cardinality larger than the continuum. We can look at the 
Heisenberg group 
$$ H(F) = \{ \left[ \begin{array}{ccc} 1 & a & c \\ 0 & 1 & b \\ 0 & 0 & 1 \end{array} \right], a,b,c \in F \} \; . $$
It is non-abelian and nilpotent of class $2$ meaning $(1+g)^3=1$ so that it has no involutions. It can not be 
generalized dicyclic as dicyclic groups contain elements of order $4$. Every finitely generated subgroup is finite. 
It requires $\kappa$ generators as its abelienization contains $F \oplus F$. In this example, the complete weighted
metric tricks breaks down as we could use only a continuum of metric pairs. There is way to turn every Cayley color
into its own real distance. This difficulty is related with an existing open problem in the 
Cayley-graph theory. Watkins conjectured that, apart from the abelian and generalized-dicyclic obstructions 
and the finite exceptional cases, every group has a GRR \cite{Watkins1975}.
The finite case is settled, and Leemann-de la Salle settled the finitely generated infinite case. 
It looks as if Leemann's current research program explicitly identifies the extension 
to groups that are not finitely generated as the remaining problem \cite{Leemann2023}.

\paragraph{}
One line of question which could be useful is within the class of non-metabelian groups because for 
finitely generated groups the non-metabelian ones are natural. What about Lie groups? \\
One can ask, is every non-metabelian Lie group natural? 
For example, $O(n)$ is non-metabelian and natural. We will answer this in the Lie group section.

\paragraph{}
We end with the remark that the Lorentz group $O(1,3)$ is reflection generated and so
natural. It is also non-metabelian. We discuss it and the Poincar\'e group more
in the Lie group section. 

\section{Finitely generated groups}

\paragraph{}
A literature search carried out with the assistance of GPT-5.6 Sol
August 2026 led to the work of Leemann and de la Salle \cite{LeemannDeLaSalle2021,LeemannDeLaSalle2022},
from which the following classification theorem emerged.  It is
especially remarkable because the finite group case illustrates
that {\bf graphical regular representation theory} is not equivalent to
the notion of {\bf natural groups}: there are finite groups such as the alternating
group $A_4$ and some small dihedral groups, which have no graphical regular representation
but nevertheless are natural. Among infinite finitely generated groups,
the groups without a graphical regular representation are precisely the abelian and 
generalized dicyclic groups. Corollary 1.2 of \cite{LeemannDeLaSalle2022} 
combines this with the classical finite classification. But the absence of 
a graphical regular representation does not yet produce non-naturality. 
For abelian and generalized dicyclic groups, some additional input are needed.

\paragraph{}
The following result illustrates that the notion of natural itself is
"natural". It aligns with graphical regular
representation theory, but the resulting classification is smoother 
because the sporadic finite {\bf GRR exceptions} disappear. 
The work of Leemann and de la Salle supplies essentially all of the difficult 
positive existence part.

\begin{thm}[Naturality Classification 
(using the Leemann-de la Salle GRR classification)]
A finitely generated nonabelian group $G$ is natural $\Leftrightarrow$ $G$ is not generalized dicyclic.
A finitely generated abelian group $G$ is natural $\Leftrightarrow$ $G$ is finite of odd order or $G \cong C_2^r$.
\end{thm}

\paragraph{}
This for example implies that if a finitely generated group $G$ has a nonabelian commutator $G'$ then 
it is natural. To check that a group is natural, one only needs to find four elements $a,b,c,d$
such that $[ [a,b],[c,d] ]$ is not the identity and have established that it is natural. 
This is not necessary however as $G=S_3=D_3$ shows. That group is natural but $G''=1$ so that 
$[[a,b],[c,d]]=1$ for all $a,b,c,d$. For $D_{\infty}=C_2 * C_2$, the 
commutator is $\mathbb{Z}$ so that again $D_{\infty}''=1$ while $D_{\infty}$ is natural. 
One can say therefore:

\begin{coro}
Every finitely generated non-metabelian group is natural. 
\end{coro}

\paragraph{}
The proof breaks into four pieces (1) Finite groups, (2) infinite abelian groups, 
(3) infinite generalized dicyclic groups, and (4) infinite non-abelian groups which are not generalized dicyclic.
The following two lemmas cover the non-natural cases if $G$ is infinite
but the first lemma does not require $G$ to be infinite:

\begin{lemma} Every finitely generated generalized dicyclic group is non-natural. \end{lemma}

\begin{proof}
Let $G={\rm Dic}(A,x)$ with $y=x^2 \neq 1$. We have seen already that $y$ is an involution in $A$.
Let $d$ be any right-invariant metric on $G$. Define $l(g)=d(g,1)$.
Then $d(u,v)=l(uv^{-1})$ and, by symmetry, $l(g)=l(g^{-1})$.
Define $i(a)=a$ for $a \in A$ and $i(g)=g^{-1}$ for $g \notin
A$. (That is $i(xa)=(xa)^{-1} = x a y$.)
Leemann and de la Salle note that this canonical extra permutation preserves every Cayley graph 
of a generalized dicyclic group. But we need that it preserves every right-invariant metric. For $u,v \in G$, 
a direct four-case calculation gives
$$ i(u) i(v)^{-1}=\left\{ \begin{array}{ll}  u v^{-1}       &  u,v \; {\rm in \; the \; same \; coset \; of \; A} \\ 
                                            (u v^{-1})^{-1} &  u,v \; {\rm in \; different cosets} \\ \end{array} \right. \; . $$
Hence $d(i(u),i(v))=l(i(u) i(v)^{-1})=l(uv^{-1})=d(u,v)$. So $i$ is an
isometry for every right-invariant metric on $G$. Now let $R_g$ denote the ordinary
right translation and define $s=R_x \circ i$. Then $s(a)=ax=x a^{-1}$, $s(x a)=a^{-1}$ implying $s^2=1$. Moreover,
$s R_a s = R_{a^{-1}}$ for $a \in A$. Therefore, $H=\langle R_a, a \in
A, s \rangle$ is isomorphic to the generalized dihedral group $H \cong A \rtimes_{-1} C_2$. 
It acts regularly on the same underlying set $G$. Indeed, the
translations $R_a$ send $1$ to all elements of A, while the elements $sR_a$ 
send $1$ to all elements of $xA$, uniquely. So, every right-invariant metric for
$G$ automatically admits a competing regular group $H$. We still must make
sure that this competitor is not abstractly isomorphic to $G$.
For every $a \in A$, $(xa)^2 = xaxa=x^2=y \neq 1$. Thus every element in
the co-set $xA$ has order $4$. In the competing generalized-dihedral
group, every element outside $A$ has order $2$. 
Let $B=\{a \in A| a^2=1 \}$. In $G$, all involutions are in $A$ while in 
$H=A \rtimes_{-1} C_2$, every element in the nontrivial coset is an involution.\\
(case i) If $G$ is finite then $G$ has $|B|-1$ involutions, while
$H$ has $|A|+|B|-1$ involutions. The two finite groups are not isomorphic. \\
(case ii) If $G$ is infinite and finitely generated, then $A$ is a finitely 
generated abelian group, so that $B$ is finite.
Hence $G$ has only finitely many involutions, whereas the competing $H$
has infinitely many. Again the groups $G,H$ are not isomorphic. \\
The group ${\rm Dic}(A,x)$ is non-natural.
\end{proof}

\begin{lemma}
Every infinite finitely generated abelian group $G$ is non-natural.
\end{lemma}

\begin{proof}
Let $G$ be an infinite finitely generated abelian group. Then $G \cong
\mathbb{Z}^r \times F$, where $F$ is a finite abelian group. 
Reduce one of the $\mathbb{Z}$-coordinates modulo 2
using $\chi: G \to C_2, \chi(n_1,\dots, n_r,f) = n_1 \; ({\rm mod} \; 2)$.
Its kernel $A={\rm ker}(\chi)$ satisfies $[G:A]=2$ so that $A$ contains
elements of infinite order. For any translation-invariant metric, inversion
$u \to -u$ is an isometry. Choose $b \notin A$. Then $u \to u+a,a \in A$,
together with $u \to b+a-u,a \in A$ form a regular group $H = A \rtimes_{-1} C_2$.
The translations send $0$ bijectively onto $A$,
while the maps $u \to b+a-u$ sends $0$ bijectively onto the other coset
$b+A$. The action is therefore regular. Since $A$ contains an element of
infinite order, $H$ is non-abelian, whereas $G$ is abelian. Hence every
translation-invariant metric on $G$ admits a non-isomorphic competing 
regular group $H$ of isometries. The group $G$ is non-natural.
\end{proof}

\paragraph{}
And here is the infinite non-abelian case proving naturalness:

\begin{lemma}
Any infinite finitely generated nonabelian 
which is not generalized dicyclic is natural.
\end{lemma}

\begin{proof}
Suppose $G$ is infinite, finitely generated, non-abelian, and not a
generalized dicyclic group. The work of Leemann and de la Salle proves
that such a group admits a finite-degree graphical regular
representation. More precisely, their classification says that among infinite finitely
generated groups, the groups without GRRs are precisely the abelian and
generalized dicyclic groups. A graphical regular representation immediately
naturalizes $G$. If its Cayley graph $\Gamma=\Gamma(G,S)$ has 
${\rm Aut}(\Gamma)=L(G)$ acting regularly, use its graph metric
$\delta$ to define the metric $d(g,h)=\delta(g^{-1},h^{-1})$.
Then $d(ga,ha) = \delta (a^{-1} g^{-1},a^{-1} h^{-1}) = \delta(g^{-1},h^{-1})=d(g,h)$
so that right translations are all isometries. 
Since inversion conjugates left regular representation
to the right regular representation. 
Indeed, if  $J(g)=g^{-1}$ then the graph metric determines adjacency so that
${\rm Iso}(G,d)=J {\rm Aut}(\Gamma) J^{-1}= J L(G) J^{-1} = R(G)$.
If another group structure on the same set has all right translations as
isometries, its right regular representation is a regular subgroup of
$R(G)$. But a regular subgroup of a regular permutation group must be equal to
the whole group: for every $g \in G$, there is exactly one element of
$R(G)$ sending $1$ to $g$. Transitivity forces that element $R_g$ to 
belong to the subgroup. Hence, the competing translation group must 
be $R(G)$.
\end{proof}

\paragraph{}
Together with the already established finite case, this finish the proof.

\section{Lie groups}

\paragraph{}
Lets look at the Euclidean case in which case Gismatullin also shows that $\mathbb{R}^n$ is strongly 
natural. For Lie groups, Gismatullin also introduced {\bf completely strongly natural}, meaning that
there exists a bi-invariant metric on which the generating metric is compatible with the topology 
of the group. Gismutallin also shows that in general if $G,H$ are strongly natural groups, then the
product is a strongly natural group. We have seen with small examples like $C_2 \times C_3$ that this
does not work in the natural case.
Since we do not the product operation at hand, we can not conclude that $(\mathbb{R}^n,d)$
with the usual Euclidean $l^2$ metric is natural. The space $\mathbb{R}^3 = \mathbb{R}^2 \times \mathbb{R}$ 
admits a group structure that is right invariant. Let $R_s$ denote the rotation
by angle $s$ in $\mathbb{R}^2$. Now, define the non-abelian group operation 
$(v,t) * (w,s) = (R_s v + w,t+s)$. The left translation $L_{(u,r)}(v,t) = (u,r) * (v,t) = (T_t u + v,r+t)$
is not an isometry: take $x=(0,0)$ and $y=(0,\pi)$ which has distance $d(x,y)=\pi$. But
$L_{(u,r)}(v,t) = (R_{\pi} u,r+\pi=(-u,r+\pi)$ satisfies 
$d(L_{(u,r)} x,L_{(u,r)} y) = \sqrt{\pi^2+4|u|^2} > \pi$. This example still shows that constant 
Riemannian metric is not a natural metric space. 

\paragraph{}
Still, $\mathbb{R}^n$  is natural. It actually falls within the class of connected Lie groups so 
that the following older argument is now not necessary any more. It is still illustrative to have
a more intuitive metric.

\paragraph{}
\begin{thm}
{\bf $\mathbb{R}^n$ is natural}
\end{thm}

\begin{proof}
We define a metric $X=(\mathbb{R}^n, l^4 \; {\rm metric})$ that natural, forces the usual translation group. 
In a nutshell, the reason is that the isometry group of $X$ is $\mathbb{R} \rtimes H$ for a finite group 
but that $\mathbb{R}^n$ has no proper subgroup of finite index.
Here are more details: It is not enough to break the rotation symmetry using a constant Riemannian metric. But we can
take $l^p$ metric different from $p=2$, which then forces the Euclidean group structure. Pick the
metric $d(x,y) = (\sum_{k=1}^n |x_j-y_j|^4)^{1/4}$. Pick any group structure $*$ on
$\mathbb{R}^n$ that has right translations as isometries. The identity element in
$(\mathbb{R}^n,*)$ can be placed at the identity $0$ of $(V,+)$. The
right translation group $K=\{ R_a, a \in \mathbb{R}^n\}$ of right translations is a regular subgroup of the isometry group of
$(\mathbb{R}^n,*)$. It is $\mathbb{R}^n \rtimes H$ with $H = Z_2^n \rtimes S_n$.
Let $V$ denote the ordinary translation subgroup of the affine isometry group of $(\mathbb{R}^n,+)$
and let $N=K \cap V$. The map $K \to H$ taking an affine isometry to its linear part has a finite image
$A$ so that $N=K \cap V$ has finite index in $K$.
Evaluation $k \to k(0)$ gives a bijection from $K$ to $V$.
Because $K$ acts regularly and $N$ is a normal subgroup of $K$, the finitely many cosets of $N$
give finitely many co-sets in the corresponding additive subgroup $A$ of $V$. Hence
$[V:A]=m$ is finite and $(R^n,+)$ is divisible, meaning $x=m(x/m)$, we must have $K=V$.
We have shown that the right translation group of $(\mathbb{R}^n,*)$ is isomorphic to $\mathbb{R}^n$.
\end{proof}

\paragraph{}
The following result is especially surprising: 

\begin{thm}
Every connected Lie group is natural.
\end{thm}

\begin{proof}
Finite-dimensional Lie group have the cardinality of the continuum.
One can therefore fix for each inverse-pair $\{g,g^{-1}\},\qquad g\neq 1$
a distinct weight $1<\lambda_{\{g,g^{-1}\}}<2$. Define the metric
$d(x,y) = \lambda_{\{ xy^{-1}, y x\} }$ if $x \neq y$ and $d(x,y)=0$ else. 
\footnote{This generates the discrete topology but naturalness allows this.}
It is right invariant: $d(x g,y g)=d(x,y)$ for every $g \in G$.
This edge-colored complete Cayley graph of $G$ has produced a metric space.
A theorem of Byrne-Donner-Sibley has been extended to arbitrary infinite groups  
by Leemann-de la Salle: if the complete Cayley graph of a group has a
nontrivial color-preserving automorphism fixing the identity, then
the group is either abelian or generalized dicyclic.
Because all $\lambda_{\{g,g^{-1} \}}$ are distinct, an isometry 
is a color-preserving automorphism of this complete Cayley graph. \\
Case (i): {\bf G is non-abelian}. A nonabelian connected Lie group can 
not be generalized dicyclic because a generalized dicyclic group has an abelian subgroup       
$A$ of index 2. In a Lie group, a finite-index subgroup must be open and
an index-two subgroup would disconnect $G$.
By the theorem, $G$ the action is regular and {\rm Isom}(G,d)=R(G). 
The metric forces $G$. Every nonabelian connected Lie group is natural.\\
Case (ii): {\bf G is abelian} A connected abelian Lie group has the form
$G \cong \mathbb{R}^n \times \mathbb{T}^m$. It is divisible; in particular $2G=G$.
For the metric above, the color-preserving theorem gives only the familiar extra inversion symmetry
$i(x)=x^{-1}$. Thus ${\rm Isom}(G,d)=G \rtimes \langle i \rangle$.
In additive notation, every isometry outside the translation subgroup has the form
$x \to a-x$. But since $2G=G$, there exists x with 2x=a. Hence a-x=x: every such reflection has a fixed point.
A non-identity element of a regular action cannot have a fixed point. 
Therefore any regular subgroup H of ${\rm Isom}(G,d)$ contains 
no element from the inversion coset. Hence $H$ is a subgroup of $G$. Transitivity then forces H=G.
Therefore the abelian case is natural as well.
\end{proof}

\begin{center}
\begin{tabular}{|c|c|c|} \hline
   Group               &   Natural?   & Strongly natural?  \\  \hline
  $\mathbb{R}^n$        &     YES         &   YES    (Gismatullin)    \\
  Compact connected Lie  &    YES         &   YES    (Gismatullin)    \\ \hline
\end{tabular} \end{center} 

\paragraph{}
The case $G=\mathbb{R}^n$ is covered but there is an easier proof here using
the {\bf Mazur-Ulam theorem}, which tells that  
every surjective isometry of a real normed vector space must be affine: $f(x)=Ax+b$
where $A$ is a linear isometry of the norm. Mazur-Ulam can be rephrased as telling
that the {\bf metric defines the affine structure}.
Choose a norm on $\mathbb{R}^n$ so that only the $\pm Id$ belong to the linear maps
that preserve it. Mazur-Ulam then immediately implies 
${\rm Isom}(\mathbb{R}^n,d) = \mathbb{R}^n \rtimes \{\pm I\}$.
Every isometry is either $x \mapsto x+a$ or $x\mapsto a-x$
Since the second type always has a fixed point: $x=a-x$ if and only if $x=\frac{a}{2}$
no non-identity map $x\mapsto a-x$ can belong to a regular subgroup of the isometry group. 
Hence every regular subgroup lies entirely in the translation group and $H\subseteq(\mathbb{R}^n,+)$. 
Since H must also be transitive, necessarily $H=\mathbb{R}^n$
This was the same mechanism which worked for a uniquely 2-divisible abelian group: the unavoidable inversion
$x\mapsto-x$ does no harm because every translated inversion $x\mapsto a-x$ has the fixed point a/2.

\paragraph{}
Connectedness was assumed as the general case is more subtle. In the abelian case,
disconnected Lie groups like $\mathbb{T}^n \times Z_4$ are non-natural. But for a case like 
$\mathbb{T}^n \times Q_8$, the case is not yet settled. Here are some examples with independent proofs:

\begin{thm}
$O(n)$ is natural for all $n \geq 1$. 
\end{thm}

\begin{proof}
The key is the Leemann-de la Salle's theorem, which tells that a group is 
orientation-rigid if and only if it is neither generalized dicyclic nor 
abelian with an element of order larger than $2$. For $n \geq 3$, the 
orthogonal group $O(n)$ cannot be generalized dicyclic. If it had an
abelian subgroup $A$ of index $2$, Then $O(n)/A \cong C_2$. 
But $SO(n)$ is perfect for $n \geq 3$, meaning $[SO(n),SO(n)]=SO(n)$.
Therefore, every homomorphism $SO(n) \to C_2$ is trivial implying $SO(n) \subset A$.
As $SO(n)$ itself already has index $2$ in $O(n)$, this forces $A=SO(n)$
which is impossible because $SO(n)$ is non-abelian. 
For $n=2$, the group $O(2)= S^1 \rtimes C_2$ is generalized dihedral, 
not dicyclic and so natural. Finally, $O(1)=S^1 \times C_2$ is an 
elementary $2$-group and is natural by the abelian classification.
\end{proof}

\paragraph{}
For the Pin groups, one has

\begin{thm}
$Pin^+(n)$ is natural for every $n \geq 1$.
$Pin^-(n)$ is natural if and only if $n \geq 3$.
\end{thm}

\begin{center}
\begin{tabular}{|c|c|c|} \hline
   Group               &   Structure                  &  natural?  \\  \hline
  $Pin^+(1)$           &   $C_2 \times C_2$ (Klein 4)  & YES        \\  \hline
  $Pin^-(1)$           &   $C_4$                       & NO         \\  \hline
  $Pin^+(2)$           &   $O(2)$                      & YES        \\  \hline
  $Pin^-(2)$           &   $Dic(S^1,-1)$               & NO         \\  \hline
  $Pin^{\pm}(n), n>2$  &   Cover of $Spin(n)$          & YES        \\  \hline
\end{tabular} \end{center}

This is interesting as $Pin^-(2) = S^1 \cup j S^1$ is a subgroup of $SU(2)$ 
and contains all the dicyclic groups $Dic_m=Dic(C_{2m},j) = \langle e^{\pi i/m},m \rangle$ ,
where $j$ is the basic quaternion unit. The argument for $n \geq 3$ is similar to the 
proof in the $O(n)$ case, it is just $Spin(n)$ will play the role of $SO(n)$:

\begin{proof}
Suppose $A \subset Pin^{\pm}(n)$ was an abelian subgroup of index $2$. 
Its identity component $Spin(n)$ is a connected compact semi-simple group 
and is perfect. Hence the quotient homomorphism $Pin^{\pm}(n) \to C_2$ 
must be trivial on $Spin(n)$. Consequently $Spin(n)\subset A$. 
But $Spin(n)$ is non-abelian for $n \geq 3$ leading to a 
contradiction. Therefore $Pin^{\pm}(n)$ for $n \geq 3$
is neither abelian nor generalized dicyclic. It is therefore 
orientation-rigid and consequently natural.
\end{proof}

\paragraph{}
For $n=4$, $Spin(4) \cong SU(2)\times SU(2)$, so it is semi-simple 
and perfect. For $n=3$, $Spin(3) \cong SU(2)$.
Thus $n=3$ is exactly the right cutoff. For $n=2$,
$Spin(2) \cong S^1$ is abelian and certainly not perfect, which 
illustrates why the proof stops precisely at $n=3$

\paragraph{}
The group $\mathbb{R} \times C_2$ is an other example of a non-connected Lie group 
that is not natural. As we can use the jumping between branches to implement
also reflections. 

\paragraph{}
It is nice to see the naturality question becomes interesting in groups that pop up in physics
like $Spin(1) \cong C_2$ (spin up or spin down) or $Spin(2) \cong S^1$ and $Spin(3)=SU(2)$
appearing as gauge groups. All these three cases are also Euclidean spheres with group 
operation. $Spin(4) \cong SU(2) \times SU(2)$ which covers $SO(4)$ which is an Euclidean
version of the Lorentz group.

\section{Three routes to naturality}

\paragraph{}
We have seen that connected Lie groups and Lie groups generated by reflections are
both natural. A third route to naturality comes with the notion of the commutator 
$G'=[G,G]$. A group is metabelian if $G''=1$.  Non-metabelian Lie groups satisfy
$G'' \neq 1$.  As usual, the assumption is that Lie groups are second countable.

\begin{thm}
Every non-metabelian Lie group is natural.
\end{thm}

\begin{proof}
If $G'' \neq 1$, the group $G$ is non-abelian and cannot be 
generalized dicyclic, (generalized dicyclic group has abelian $G'$).
A Lie group has cardinality at most the continuum. 
Now use the Leemann-de la Salle orientation-rigidity result with 
the complete symmetric generating set $S=G \setminus \{1\}$. 
Their general orientation-rigidity result applies to groups 
that are neither abelian of exponent larger than $2$ nor 
generalized dicyclic. The finitely generated theorem 
in the 2022 paper is the finite-generating-set specialization. 
The distinct-distance construction on inverse pairs
$\{g,g^{-1} \}$ works as $G$ has the cardinality not exceeding the continuum.
Thus an isometry fixing 1 is forced to be the identity, and
$Iso(G,d)=R(G)$.
\end{proof}

\begin{coro}
A non-natural Lie group must be disconnected and metabelian.
\end{coro}

\paragraph{}
Examples:
\begin{itemize}
\item $\mathbb{R} \times C_2$, the double line which is non-natural
\item $Pin^-(2) = S^1 \rtimes C_2$ which is non-natural 
\item $G=O^2$ is metabelian as $G'=SO(2)$ and it is natural.
\item $O(n)$ is not metabelian for $n \geq 2$ as $O(n)'=SO(n)$. It is natural, also because reflection generated.
\item $E(n) = R^n \rtimes O(n)$ is disconnected. There is a second reason now why it 
is natural: it is not-metabelian. Indeed it contains for $n \geq 2$ the perfect $SU(n)$ (meaning $SU(n)'=SU(n)$).
\end{itemize}

\paragraph{}
Since  $E(1)$, the infinite dihedral Euclidean group is natural (as it is generated by involutions)
and $E(2) = \mathbb{R}^2 \otimes O(2)$ is natural, all Euclidean groups $E(n) = R^n \rtimes O(n)$
are natural. We have in the example of Euclidean groups
{\bf three independent routes to naturality}: 

\begin{itemize}
\item The Euclidean groups $E(n)$ are natural by reflection generation.
\item The special Euclidean $SE(n)$ are natural by connectedness.
\item The Euclidean groups $E(n)$ for $n \geq 3$ are natural because they are non-metabelian.
\end{itemize}

\paragraph{}
The {\bf three routes to naturality} are nicely seen in the example of the Lorentz group $O(1,3)$
and the Poincar\'e group $P(1,3) \mathbb{R}^{1,3} \rtimes O(1,3)$. 
The group $O(1,3)$ which is 4-connected but by Cartan-Dieudonn\'e reflection generated as well. 
It is also non-metabelian. Its identity component $O^+(1,3)$ is called the {\bf proper orthochonous Lorentz group}. 
It is natural by connectedness. The Poincar\'e group is also generated by involutions.
Proof: its Lorentz subgroup is generated by involutions by Cartan-Dieudonn\'e. Look also the
involutions $P(x)=-x$ and $Q(x)=a-x$ which are both Poincar\'e transformations that are involutions.
But $QP(x)=a+x$ is translation so that every translation is a product of two involutions.
The Poincar\'e group is also non-metabelian.
Both the full Lorentz and full Poincar\' groups are disconnected but natural for two 
independent reasons: by reflection generation and non-metabelianity. 
Their identity components are natural by connectedness.

\begin{tabular}{|c|c|c|c|} \hline
  Group                              &  Involution generated? &  Connected? & Non-metabelian? \\ \hline
  $O(1,3)$                           &   YES                  &    NO       & YES             \\ 
$SO^+(1,3)$                          &   NO                   &    YES      & YES             \\ 
$\mathbb{R}^{1,3} \rtimes O(1,3)$    &  YES                   &    NO       & YES             \\
$\mathbb{R}^{1,3} \rtimes SO^+(1,3)$ &   NO                   &    YES      & YES             \\ \hline
\end{tabular}

\section{$p$-adic integers}

\paragraph{}
The original attempts to decide about the naturality using density or using the 
natural p-adic metric do not hold. Also other approaches do not hold. But the
arguments already given, work also again. We repeat them because we have flipped 
so many times in deciding whether $\mathbb{Z}_p$ are natural or not.

\begin{thm}
$\mathbb{Z}_p$ is natural if and only if $p$ is odd. 
\end{thm}

\begin{proof}
For every prime $p$, the usual $p$-adic metric does not make $\mathbb{Z}_p$ natural.
$x=a_0+a_1p+a_2p^2+\cdots, a_j\in\{0,\cdots,p-1\}$. As a metric space,
$\mathbb{Z}_p \simeq \prod_{n=0}^{\infty} \mathbb{Z}/p\mathbb{Z}$
with $d(x,y)=p^{-k}$, where $k$ is the first coordinate at which the two digit sequences differ.
Now put coordinate-wise addition modulo $p$ on the digits. This gives the group
$\prod_{n=0}^{\infty} C_p$. Its translations preserve the first coordinate at which two sequences 
differ, hence are isometries of the usual $p$-adic metric.
But $(\mathbb{Z}_p,+)$ is torsion-free, whereas $\prod_{\mathbb{N}} C_p$ has exponent p. 
Thus they are not isomorphic. The usual $p$-adic metric never forces the additive group
$\mathbb{Z}_p$. \\
(i) Assume first $p$ is odd. We construct a metric which does force $\mathbb{Z}_p$.
The crucial algebraic fact is $2G=G$. Multiplication by 2 is 1:1 because $2$ is a unit in 
$\mathbb{Z}_p$. We now construct a very rigid metric similarly as in other proofs.
Consider the unordered inverse pairs $\{g,-g\}$ for $g \neq 0$
There are continuum many of them, so choose distinct real numbers
$\lambda_{\{g,-g\}} \in (1,2)$ for all such pairs. Define
$d(x,y)=0$ for $x=y$ and $d(x,y)=\lambda_{\{x-y,y-x\}}$ else. 
This is a metric as seen before. It is also translation invariant because $d(x+g,y+g)=d(x,y)$
for all $g$. Let $T$ be an isometry fixing $0$. Since all the numbers $\lambda_{\{g,-g\}}$
are distinct, $d(0,T(x))=d(0,x)$ implies $T(x)\in \{x,-x\}$.
Thus every point can apparently only be fixed or inverted.
But the choices cannot be made independently. Suppose for two nonzero elements
$T(x)=x$, $T(y)=-y$. Since distances uniquely determine differences up to sign,
$T(x)-T(y)=x+y$ must equal either $x-y$ or $-(x-y)$.
The first possibility gives $x+y=x-y$ showing $2y=0$
The second gives $x+y=-x+y$ showing $2x=0$. Both are impossible.
Consequently an isometry fixing 0 is either $T(x)=x$ for all $x$
or $T(x)=-x$ for every $x$. We have now ${\rm Iso}(G,d)=G \rtimes C_2$,
where $G$ acts by translations and $C_2$ acts by inversion.
This is exactly the smallest isometry group one could hope 
for with an abelian group, since inversion is unavoidable for every invariant metric.
Any competitor group $H \leq {\rm Iso}(G,d)$ would act 
regularly on $G$. Such a group $H$ would correspond to another group structure 
compatible with the metric. Every element of the isometry group is either
$x\mapsto x+a$ or $x \mapsto a-x$. But for odd $p$, every map of the second 
type has a fixed point, namely $x=\frac{a}{2}$ because division by $2$ 
is possible in $\mathbb{Z}_p$. A non-identity element of a regular action 
cannot have a fixed point. Therefore $H$ contains no reflection $x \mapsto a-x$.
Thus $H \subset G$, where the right-hand G denotes the translation subgroup.
But $H$ must be transitive. A subgroup of the translation group is transitive 
on G only if it is the whole translation group. Hence $H=G$. \\

(ii) Assume now $p=$. We have $2 \mathbb{Z}_2 \neq \mathbb{Z}_2$, because
$2 \mathbb{Z}_2$ has index $2$. For every translation-invariant metric on the abelian group 
$\mathbb{Z}_2$, both translations and inversion are isometries.
Consider the competitor $H= \{x \mapsto x+a | a \in 2\mathbb{Z}_2\} 
\cup \{x \mapsto b-x | b \notin 2\mathbb{Z}_2\}$. 
This acts regularly. Indeed, a reflection $x \mapsto b-x$ with $b \notin 2\mathbb{Z}_2$ 
has no fixed point, because a fixed point would require $2x=b$ which has no solution.
The resulting regular group is non-abelian. It is the
generalized dihedral group of $2\mathbb{Z}_2$so it cannot be isomorphic 
to the additive group $\mathbb{Z}_2$. As this construction works 
for every invariant metric, the group $\mathbb{Z}_2$ is not natural.
\end{proof}

\paragraph{}
Again, there is useful statement hiding here:
for an abelian group G of cardinality at most the continuum, 
with no 2-torsion, the same construction 
gives that unique $2$-divisibility is the dividing line: 
if $2G=G$, the rigid metric makes $G$ natural; if 
$2 G\neq G$, an index-two quotient gives the 
generalized-dihedral competitor $H$.

\section{Examples}

\paragraph{}
This section appears redundant now (2026),as the class of groups for which naturalness is settled
has grown. The section can still be useful pedagogically. For me, it is a bit nostalgic also 
to look at the first steps (2008) in a mathematical research area which was mostly new and where 
the learning curve was steep for me (and so treacherous) but where it has been an exciting adventure. 
I got into the topic originally by looking for a mathematical topic that looked unexplored and still 
was elementary so that an aspiring mathematician can learn in it and make discoveries. 
I certainly had a lot of fun, especially as there were struggles and the 2022 to 2026 transition 
illustrate the struggle. Finding unexplored topics might (for a while at least) the only thing in which 
we are still better than machines. Already in 2008, search engines had become very powerful and asking a standard
problem, say in real analysis or linear algebra, would be answered well. Today it is
almost impossible to look something up, without given expert answers. Search has become deeper and
finding sources has become easier.
\footnote{Around 2004-2007 some members of our department have investigated, how much flavors of 
inquiry based teaching could enrich larger service courses.}
While inquiry based learning is tough in larger courses as it requires to rely on students not to look things up,
I personally would recommend a beginning mathematician to explore at least one pet project without 
checking literature or even discussing with peers, just to practice discovery.
This paper about natural graphs had been such a pet project (at least until 2022, today in 2026, under pressure
to fix earlier mistakes (and rescue a stranded former me,
I have also relied more on looking up literature. The machines have
become more effective in finding the right sources and finding mistakes or even answer research
questions (as I illustrate below). It has become much more difficult to think independently.

\paragraph{}
An example of a natural graph is the {\bf complete graph} $K_p$ 
for which the number $p$ of vertices is a prime number. 
The corresponding metric generates the {\bf discrete topology} on a finite set with $p$ elements. 
There is only one group structure because there is up to isomorphism only one finite group $G$
of prime order $p$. It is the {\bf cyclic group} $C_p=\mathbb{Z}/(p \mathbb{Z})$ with $p$ elements. 
This follows from Lagrange's theorem as if we take an arbitrary element $g \in G$ different
from the identity, it generates a subgroup $\langle g \rangle$ whose order must divide $p$.
As it contains at least $2$ elements $\{1,g\}$ and its order must divide $p$, it must be equal
to $p$, so that the cyclic group agrees with the full group $G$. There would be other metric 
spaces that force the same group. We could take the Cayley graph of the cyclic group $C_p$
with generator set $S=\{1,-1\}$, which is the cycle graph. We have \\

{\bf Any finite group with prime order is natural. }

This is much weaker than the general classification mentioned in the previous section, but it
is a statement that can be checked immediately and does use any distinction whether the group is
abelian or not. 
\footnote{The notation $C_p$ is more common than $Z_p$ in group theory. We use it to avoid confusing
it the group of $p$-adic integers $\mathbb{Z}_p$.}

\paragraph{}
The smallest example of a non-natural group is the cyclic group $C_4$. We can think of it 
as ${\rm Dic}_1$, because dicyclic groups have been seen cases in the non-abelian case as
examples of non-natural groups. 
The reason why $C_4$ is not natural is because any metric space with four points which admits this
group as a symmetry group also allows to host the {\bf Klein 4-group} $D_2 = C_2 \times C_2$.
The Klein group $D_2$ is the smallest {\bf dihedral group}. It is natural. In order
to see that, we just need a metric that breaks the symmetry to exclude the second group $C_4$ or
order $4$. We replace therefore the geodesic metric on the cyclic graph $C_4$ so that only $C_2 \times C_2$ remains as 
a possible group structure. This can be done for example by taking alternating the distance $2$ and $3$ so
that we take the Euclidean metric of a $2 \times 3$ rectangle restricted to the vertices.

{\bf The groups $C_{2^k}$ are not natural for $k \geq 2$. }

\paragraph{}
The case of the prime $p=2$ is very special because all groups $C_{p^k}$ for odd 
primes $p$ have odd order and so are natural. The {\bf ``make it Coxeter!" construction} like 
$C_4 = \langle a, a^4=1 \rangle \to D_4 = \langle a,b, a^2=b^2=(ab)^2=1 \rangle$ 
is not available as a tool to build an other non-equivalent but compatible group. 
P.S. $C_4$ is also not strongly natural because all groups of order $4$ are abelian.
There is no distinction between natural and strongly natural for groups with 4 elements.

\paragraph{}
Up to isomorphism, there are two groups of order $6$. Only $D_3$ is natural.
It is also known as the symmetric group $S_3$ of all permutations on a $3$-point set. 
This non-abelian group is also the smallest non-abelian dihedral group and a semi-direct
product $D_3 = C_3 \rtimes C_2$. In the case $D_3$, the {\bf zig-zag product} of $C_3$ and $C_2$ 
with a suitable weighted metric produces the natural metric space. 
The zig-zag product of $C_3$ and $C_2$ is isomorphic to the {\bf utility graph} $K_{3,3}$. 
While the utility graph both admits $D_3$ and $C_2 \times C_3$ as symmetry group, it is
not a natural graph. It carries a modified weighted metric that breaks the symmetry.
Excluding $C_6 = C_2 \times C_3$ as such, it forces the group structure $D_3$. 

\begin{figure}[!htpb]
\scalebox{0.4}{\includegraphics{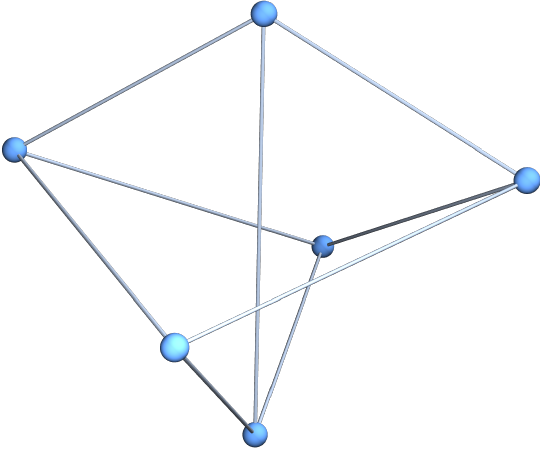}}
\label{Figure 1}
\caption{
The zig-zag product of the Cayley graphs of $C_3$ and $C_2$ relates
to the semi-direct product $D_3=C_3 \rtimes C_2$ is the bipartite  
graph $K_{3,3}$ which is also known as the {\bf utility graph}. The graph
itself is not natural, as it features both the direct product $C_6=C_2 \times C_3$ as well as 
the semi-direct product $D_3$. However, it can be weighted to admit only the $D_3$ 
group structure. The group $C_2 \times C_3$ is not natural. 
}
\end{figure}

\paragraph{}
The second group of order $6$ is the cyclic group $C_6$ which is the {\bf direct product}
$C_2 \times C_3$ by the {\bf fundamental theorem of abelian groups}. In the direct product case,
$C_3 \times C_2$, we could also try the {\bf Shannon product} of the Cayley graphs $C_3$ and 
$C_2$ of $C_3$ and $C_2$ with a weighted metric. Still, also the Shannon product still 
admits $D_3=C_3 \rtimes C_2$ so that $C_3 \times C_2$ is not natural. 
The Shannon product is also called {\bf strong product}, but Shannon \cite{Shannon1956}
was the first who defined the product in 1956.
Together with a group completed disjoint addition, it produces the {\bf Shannon ring}, 
a commutative ring structure on {\bf signed graphs} 
(see \cite{RemarksArithmeticGraphs} for more on this). \\

{\bf The group $C_{2} \times C_n$ is not natural for any odd $n$. }

\paragraph{}
A finite simple connected graph $(V,E)$ for which the cardinality condition 
$|V|=|{\rm Aut}(V,E)|$ is satisfied and for which the automorphism group acts regularly, is natural. 
It is called a {\bf graphical regular representation} of its automorphism group. 
(The graphical regular representation of a group is a Cayley graph whose only graph symmetries
are the groups right multiplication actions.) 
A characterization of Cayley graphs which are graphical regular representations is not known. 
Because a potential group structure on the vertices is part of the automorphism group of
the graph, graphs for which $|{\bf Aut}(V,E)|<|V|$, are never natural. \\

{\bf Graphs for which $|{\bf Aut}(V,E)|<|V|$ are never natural.}

The {\bf complete graph} $K_p=(V,E)$ for an odd prime $p$ is an example of a natural graph with 
$p=|V|<|{\rm Aut}(V,E)|=p!$. It is a natural graph because $p$ being prime forces a cyclic group 
structure. The graph $K_p$ is not a graphical regular representation because
the automorphism group ${\rm Aut}(K_p)=S_p$ is the symmetric group $S_p$
which has larger order $p!$ in comparison to $p$ for $p>2$. 

\paragraph{}
{\bf The {\bf integers} $\mathbb{Z}$ are not natural.}

Proof: given any right translation invariant metric 
$(\mathbb{Z},d)$ (there are many like for example the metric induced when placing the integers
on a {\bf helix} 
$$n \to [\cos(a n),\sin(a n),b n] $$ 
in $\mathbb{R}^3$ and taking the induced 
Euclidean metric from $\mathbb{R}^3$), there is a second, non-abelian group structure on this metric space.
It is the {\bf infinite dihedral group} $D_{\infty} = \langle a,b , a^2=b^2=0 \rangle$.
One can take $a(n)=1-n,b(n)=-1-n$ which satisfy $ab(n)=n+2,ba(n)=n-2$. 
The action is transitive, since translations give all even integers and the operation $a$
gives the odd integers; it is free because every reflection has its fixed point
at a half-integer. Thus it is a regular action that gives a non-isomorphic group structure for 
every right-invariant metric on $\mathbb{Z}$.
The group $D_{\infty}$ preserves the metric space also but it is non-isomorphic to $\mathbb{Z}$ 
because it is non-abelian. The zero element is the {\bf empty word}, the addition is done by 
{\bf concatenating} words, but there are two branches. \\

{\bf The infinite dihedral group $D_{\infty}$ is natural.} 

Proof: define a metric on the set $\mathbb{Z}$ in which 
$d(0,a) \neq d(0,b)$. Now, the $\mathbb{Z}$ group structure is no more possible. 
There is only group structure for which the right translation is an isometry. 
We see this more clearly in the general classification result on non-abelian finitely 
generated groups. 

\paragraph{}
The transition of going from $\mathbb{Z}$ to $D_{\infty}$ is an example
of a {\bf Higman-Neumann-Neumann} (HNN) extension. Unlike for finite groups, for infinite groups, 
this is possible without changing the cardinality, illustrated by 
 the {\bf Hilbert hotel metaphor}. 
The pair of examples $\mathbb{Z}, \mathbb{D}_{\infty}$ was the starting point for our 
research that started in 2008 while exploring topics which would be suitable for an inquirey 
based course involving mathematical that have not appeared yet so that students can not google the
proofs. It was fascinating already then to see that we can have an infinite metric space 
that determines a unique group structure up to abstract group isomorphism.

\paragraph{}
In July 2026, Jakub Gismatullin pointed out that my original 
definitions were inconsistent. He saw that {\bf A) $\mathbb{Z}$ is {\bf strongly} natural.} 
(in contrast to the fact that $\mathbb{Z}$ is not natural).\\
{\bf Proof:} Take the metric $d(m,n)=|m-n|$. Its isometry group is $D_{\infty}$
generated by the generators $a(x) = 1-x$ and $b(x) = -1-x$. They are
two non-commuting isometries. Both $a$ and $b$ preserve the distance metric.
A regular subgroup of this isometry group is either the translation group $\mathbb{Z}$, 
or the half step action of $D_{\infty}$. But the latter group structure is not bi-invariant 
as $d(1,2)=1,d(1,-2)=3$. The only strongly compatible group structure is $\mathbb{Z}$. \\
{\bf B) $D_{\infty}$ is not {\bf strongly} natural:}
{\bf Proof:} write $D_{\infty} = \langle r,s | s^2=1, srs=r^{-1} \rangle$. If $d$ is any
bi-invariant metric on $D_{\infty}$,
then the left translation $L_r$ and the right translation $R_s$ are both isometries
that commute and so generate $\mathbb{Z} \times C_2$. But the later group
is abelian. It acts regularly because $RL_r^n R_s^l = r^n s^l$.
So, every bi-invariant metric on $D_{\infty}$ also supports a compatible 
abelian group structure.

\paragraph{}
{\bf Fundamental groups} of surfaces can be both natural or non-natural. 
The {\bf projective plane} $\mathbb{P}^2$ has the natural fundamental 
group $C_2=\mathbb{Z}/(2\mathbb{Z})$. The $n$-torus has the fundamental group 
$\mathbb{Z}^n$, which is non-natural. \\

{\bf $\mathbb{Z}^n$ is not natural.}

\paragraph{}
More interesting is the fundamental group $G$ of the {\bf Klein bottle}
is $G=\langle a,b | a b a b^{-1} \rangle$ which is the
semi-direct product $H \rtimes K = \mathbb{Z} \rtimes \mathbb{Z}$, where
$H$ is the normal subgroup generated by $b$ and $K$ is the 
subgroup generated by $a$. While every element in $G$ can be written as
$b^n a^m$, the multiplication is not the direct product. The later 
would lead to the fundamental group $\mathbb{Z}^2$ and not $\mathbb{Z} \rtimes \mathbb{Z}$. 
Indeed, the Klein bottle is a {\bf non-trivial fiber bundle} of the circle, 
while the 2-torus is a {\bf trivial fiber bundle}.
The classification theorem shows that the group $G$ is natural. We just have to 
exclude that it is generalized dicyclic. But the Klein bottle group 
is torsion-free. Indeed, if we identify $a^m b^n$ with $(m,n)$, with multiplication
$(m,n)(r,s)=(m+(-1)^n r,n+s)$, then if $(m,n)^k=(0,0)$ for some k>0, the 
second coordinate gives $kn=0$, hence $n=0$; then the first coordinate gives 
$km=0$, hence $m=0$. Having no nontrivial finite-order element shows that
we do not have a generalized dicyclic case. \\

{\bf The Klein bottle group is natural.}

\paragraph{}
The fact that $\mathbb{R}^n$ is natural has been covered before.
The original 2022 argument used
the product property. The product property can be used when using {\bf strongly natural}.
Gismatullin could show, using the product 
property that $\mathbb{R}^n$ is is strongly natural.) We still have.
that {\bf $\mathbb{R}^n$ is natural}
As covered in a previous section, a connected Lie group is natural. \\
The original argument that compact Lie groups are natural does not apply 
as it used a bi-invariant metric and so refers more to {\bf strong natural}. 
Gismatullin covers the {\bf strongly natural} case.

\paragraph{}
Our old argument used that if a group admits a {\bf bi-invariant Riemannian metric} 
$d$ it must be a {\bf topological group} in which the metric and the 
arithmetic are compatible. By Gleason and Montgomery-Zippin, a topological group $(X,d,+)$ 
for which $(X,d)$ is a topological manifold, has a unique Lie group structure,
because any other compatible group operation again will be a Lie group so that there is an 
isomorphism. (Gleason-Montomery-Zippin more generally establish that
if a locally compact topological group $G$ is a projective limit of 
Lie groups and $G$ has no small subgroups then $G$ is a Lie group).
The fact that a {\bf compact connected Lie group is strongly natural} is also constructive. 
Given such a metric space, the group structure can actually be constructed: 
while there can be different metrics, 
if $G$ is not simple, any choice of a {\bf Riemannian manifold} 
$(G,d)$ for which $d$ is bi-invariant fixes the group structure: pick a point $0$ and call it 
``zero". Now we use the metric (using the {\bf fundamental theorem of Riemannian geometry}) to
construct a unique connection $\nabla$ and so get a {\bf Lie algebra structure} 
$[X,Y]=\nabla_X Y - \nabla_Y X$. The Lie algebra determines the group 
operation on a small ball $B_r(0)$ (with $r$ smaller 
than the radius of injectivity in the Riemannian manifold) by the {\bf Baker-Campbell-Hausdoff
formula}. Knowing how to add group elements in a small ball $B_r(0)$ 
defines then the group structure globally, because the set of group elements in $B_r(0)$ generate $G$. 

\paragraph{}
Among the {\bf Euclidean spheres} in $\mathbb{R}^n$ we now that the $0$-sphere, $1$-sphere and $3$-sphere are
natural and that the even dimensional spheres are not natural. For the others, we do not know more yet.
For the $0$-sphere $C_2$, things are settled because this is a finite group 
of prime order forcing the unique group $C_2$ of order $2$. 
The group structure on positive-dimensional cases $S^1=U(1),S^3=SU(2)$
comes from the multiplicative group structure given by the list of 
{\bf real associative normed division algebras} $\mathbb{R}, \mathbb{C},\mathbb{H}$
classified by the {\bf Frobenius theorem} \cite{Frobenius1878}. 
Since they are compact Lie groups, the multiplication on them is unique. 
For all other Euclidean spheres we so far we only know about the {\bf strong natural question}
they admit a topological group structure at all \cite{Samelson1940}. 

\paragraph{}
The fact that the circle $\mathbb{T}^1$ is natural but the integers $\mathbb{Z}$ is not natural 
already shows that the {\bf Pontryagin duality operation} $G \to \overline{G}$ 
does not preserve the class of natural groups. An other example is $\mathbb{Z}_2$ which is
not natural. The Pontriagyn dual, the Pruefer $2$-group is natural. 
The $\mathbb{Z}_p$ groups for odd $p$ are natural however. 
Unique 2-divisibility is the dividing line: 
if $2G=G$, the rigid metric makes $G$ natural; if $2G \neq G$, an index-two quotient 
gives the generalized-dihedral competitor. This is the reason why $\mathbb{Z}_2$ is not
natural, while $\mathbb{Z}_p$ is natural for odd $p$. \\

\paragraph{}
We have seen that the {\bf alternating groups} $A_n$ and {\bf symmetric groups} $S_n$
are all natural. One can get the metric by realizing the metric space
as a weighted Cayley graphs using two generators in the group. 
The small $A_n$ and $S_n$ ones are especially interesting. 
It appears that non-natural non-abelian finite groups are harder to find. 
But there are some: the {\bf quaternion group} $Q_8$ is not natural. We discuss this below
and see that any metric an $8$-point set which is compatible with $Q_8$ 
also will allow for the group $C_2^3 =C_2 \times C_2 \times C_2$. 

\begin{figure}[!htpb]
\scalebox{0.2}{\includegraphics{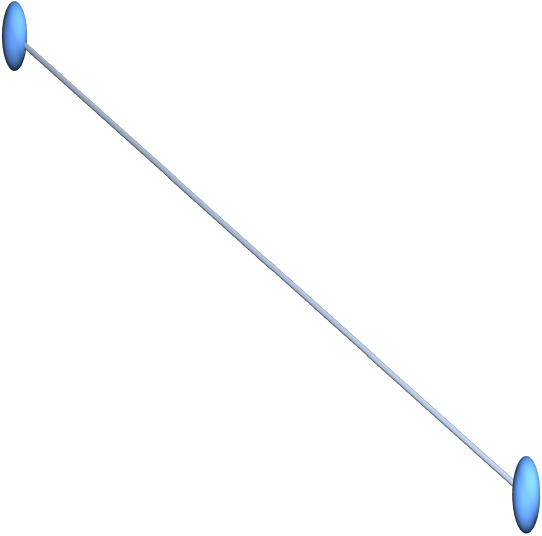}}
\scalebox{0.2}{\includegraphics{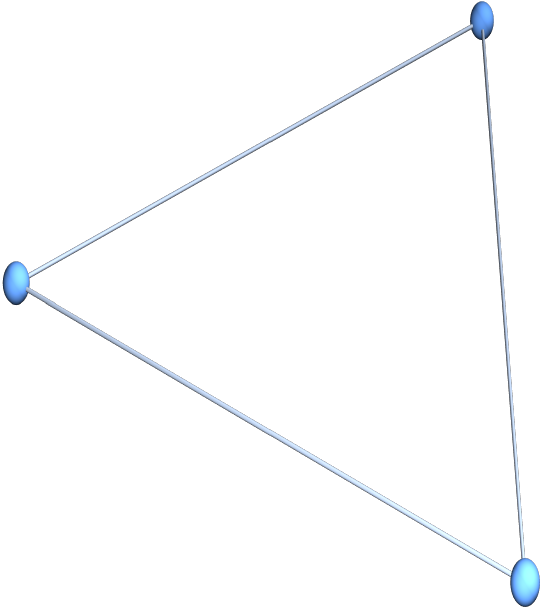}}
\scalebox{0.2}{\includegraphics{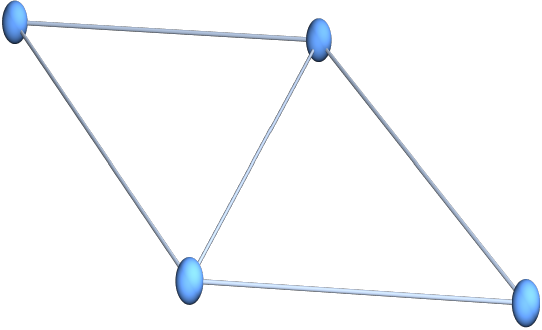}}
\scalebox{0.2}{\includegraphics{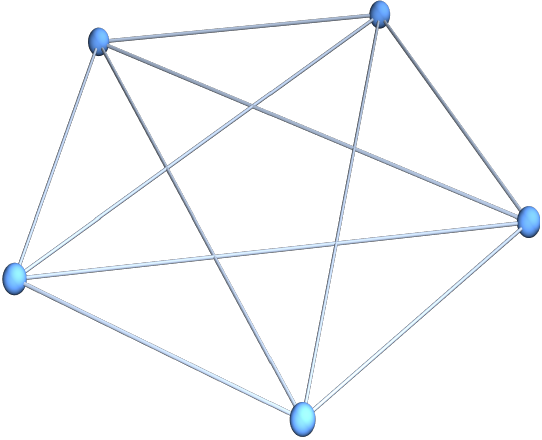}}
\scalebox{0.2}{\includegraphics{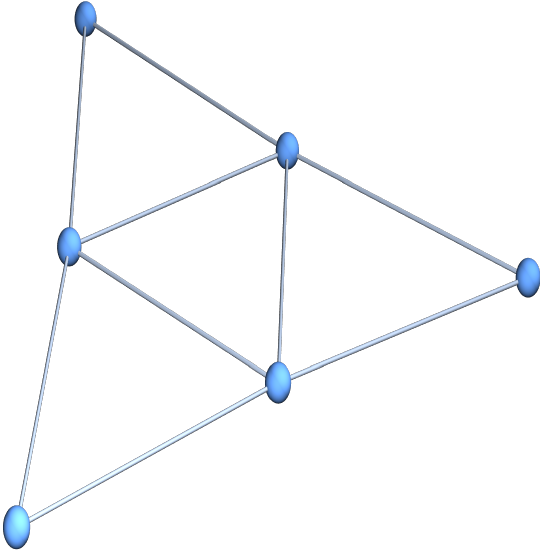}}
\scalebox{0.2}{\includegraphics{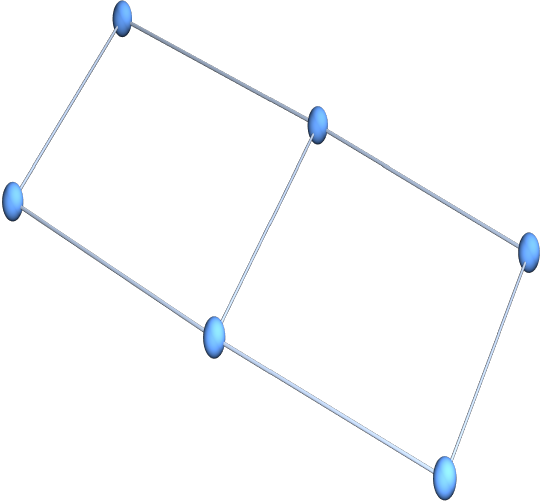}}
\label{Figure 1}
\caption{
Some natural metric spaces defining small groups 
$C_2,C_3,C_2^2,C_5,C_2 \times C_3,S_3$. Up to order $6$, 
the finite group $C_4=\mathbb{Z}/(4 \mathbb{Z})$ and 
$C_2 \times C_3$ are non-natural.
}
\end{figure}

\paragraph{}
For abelian groups, naturalness extends to the metric completion:

\begin{thm}
If $G$ is an abelian natural group $G$, then the completion $\overline{G}$ (using the natural 
metric generating it) is natural. 
\end{thm}

\begin{proof} 
Assume there was an alternative group operation $\oplus$ preserving 
the metric on $\overline{G}$ and the group $(\overline{G},\oplus)$ was not isomorphic to 
$(\overline{G},+)$, then, when the restriction $(G,\oplus)$ would not be isomorphic
to $(G,+)$ and $G$ would not be natural. 
The reason is that if the restriction would suddenly be isomorphic to the abelian $(G,+)$,
then it would also be abelian. An invariant metric on an abelian group is however
bi-invariant and its operations are uniformly continuous and extend uniquely to the 
metric completion, contradicting that $(\overline{G},+)$ and $(\overline{G},\oplus)$ are 
not isomorphic.
\end{proof}

\paragraph{}
The reverse is not necessarily true as $\mathbb{Z}$ is not natural but $\mathbb{Z}$ is 
dense in $\mathbb{Z}_p$ and $\mathbb{Z}_2$ is natural, while $\mathbb{Z}_p$ is not natural.
We need therefore look at the rational numbers $\mathbb{Q}$:

\begin{thm}
a) The rational additive group $(\mathbb{Q},+)$ is natural.
b) The rational circle $\mathbb{Q}/\mathbb{Z} \subset (\mathbb{T},+)$ is natural. 
c) The $p$-rational numbers $\mathbb{Z}[1/p] \subset (\mathbb{R},+)$ are natural. 
d) The {\bf Pr\"ufer p-group} $\mathbb{Z}[1/p]/\mathbb{Z} = \{ e^{2\pi i k/p^n}, k,n \in \mathbb{Z} \} \subset (\mathbb{T},+)$ is natural.
\end{thm}

\begin{proof}
The naturality in each case follows from the fact that the group is abelian, torsion-free, 
and uniquely 2-divisible group. The argument uses a weighed Cayley graph. We only do it for one case:
Choose, for every pair $\{g,-g\}$, $g \in \mathbb{Q} \setminus \{0\}$, a distinct real number
$1<\lambda_{ \{q,-g \} }<2$. Define $d(x,y)=0$ if $x=0$ and $d(x,y)=\lambda_{ \{x-y,y-x \} }$ if 
$x \neq y$. This is a symmetric translation-invariant metric. The triangle inequality is 
automatic because every nonzero distance is $<2$, while the sum of two nonzero distances is $>2$.
Now let $T$ be an isometry fixing $0$. Since all inverse pairs have different distances from $0$, 
$T(g) \in \{g,-g\}$. Suppose there are nonzero $x,y$ with $T(x)=x$, $T(y)=-y$. 
Preservation of $d(x,y)$ implies $\{T(x)-T(y),T(y)-T(x) \} =\{ x-y, y-x \}$, so 
$\{x+y,-x-y\}=\{x-y,y-x\}$. Hence either $x+y=x-y$, giving $2y=0$, or $x+y=-x+y$,
giving $2x=0$. Since $\mathbb{Q}$ has no 2-torsion, both are impossible.
Therefore every isometry fixing $0$ is either $g \mapsto g$ for all $q$, or $g \mapsto -g$
for all $g$. Thus ${\rm Isom}(\mathbb{Q},d) = \mathbb{Q} \rtimes C_2$, with $C_2$ acting by inversion.
Now let $H$ be any regular subgroup of this isometry group. Every element outside the translation 
subgroup has the form $g \mapsto a-g$. But $\mathbb{Q}$ is divisible by $2$, so $a/2 \in \mathbb{Q}$, and
$a-\frac{a}{2}=\frac{a}{2}$. Thus every such reflection has a fixed point. 
No non-identity element of a regular action can have a fixed point, so that $H$ contains no reflection. 
Hence $H \subset (\mathbb{Q},+)$, transitivity forces $H=(\mathbb{Q},+)$.
Therefore $(\mathbb{Q},+)$ is natural.
\end{proof}

\paragraph{}
The non-natural dyadic integers $\mathbb{Z}_2$, are Pontryagin dual to the Pr\"ufer p-group
$\mathbb{Z}[1/p]/\mathbb{Z}$ that is natural. While $\mathbb{Z}_2$ was not natural, its dual
is natural. This parallels that the natural Lie group 
$\mathbb{T}$ is dual to the non-natural $\mathbb{Z}$. 
It is not a surprise not to have Pontryagin stability.
Indeed $\mathbb{Z}_p$ for odd $p$ does not have the group operation forced from the 
usual $p$-adic metric. It has become natural for a different metric. 

\paragraph{}
Let us look at the {\bf modular group} $PSL(2,\mathbb{Z})$
which is isomorphic to the free product $C_2 * C_3$ \cite{Alperin1993}. 
(To see the free $PSL(2,\mathbb{Z}) \cong C_2 *C_3$, notice
that $A=\left[\begin{array}{cc} 1 & 1 \\ 0 & 1 \end{array} \right]$,
$B=\left[\begin{array}{cc} 0 & -1 \\ 1 & 0 \end{array} \right]$ 
generate $PSL(2,\mathbb{Z})$ so do $x=BA$ and $y=B$ 
but $x^3=1$ in $PSL(2,\mathbb{Z})$ and $y^2=1$ in $PSL(2,\mathbb{Z})$ and no words in 
$x,y$ produces the identity in $PSL(2,\mathbb{Z})$ so that
no further relations exist and $PSL(2,\mathbb{Z}) = \langle x,y | x^2=y^3=1 \rangle
=C_2 * C_3$. )

\paragraph{}
The Cayley graph of $PSL(2,\mathbb{Z})$ is the line graph of the $(2,3)$ biregular 
Bass-Serre tree. If we look at $\langle a,b | a^2=b^3=1 \rangle$. It has triangles
and is not a tree but its line graph is a $(2,3)$ is a bi-regular Bass-Serre tree. 
Lets look at the original Cayley graph of $PSL(2,\mathbb{Z})$ which has triangles. \\

{\bf The Cayley graph $\Gamma$ of $C_2 * C_3$ is a natural graph.}

Let $(\Gamma,*)$ be a group structure for which right translations are isometries. 
We can pick any of the points as the identity $e$. Its unit sphere consists of 3 
points $S(e) = \{\alpha,\beta,\gamma\}$ where $\alpha$ is isolated in $S(e)$ and
$(\beta,\gamma)$ are en edge in $S(e)$ (producing the triangle). Every isometry 
must preserve this structure. We therefore have $\alpha^{-1}=\alpha$. In order
that the edge $(\beta,\gamma)$ is preserved we must have $\gamma=\beta^{-1}$ 
producing $\beta^3=1$ (going around the triangle). We have shown that $C_2 * C_3$ is
the only group structure on that graph. 
so that the modular group $PSL(2,\mathbb{Z})$ is natural. It is a finitely generated
non-abelian group so that the general classification result could be used too.

% Bass-Serre theorem: (Gamma,Gcal) group graph. Exists G, tree T, s.t. quotient group graph is (Gamma,Gcal). 
% https://arxiv.org/abs/2002.10639 general definition of free product of graphs, not necessarily vertex transitiv

\begin{figure}[!htpb]
\scalebox{0.4}{\includegraphics{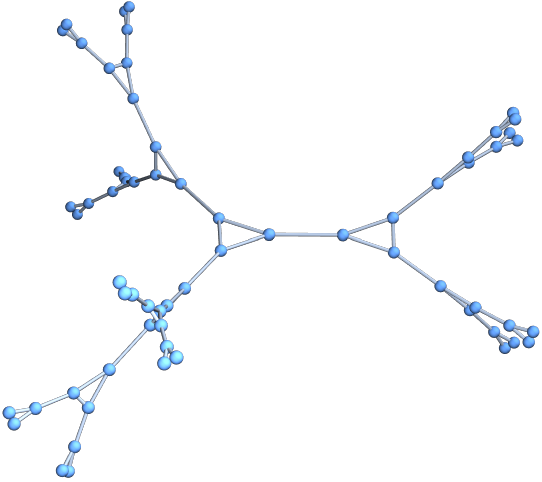}}
\label{Figure 1}
\caption{
The figure shows part of the Cayley graph of $PSL(2,\mathbb{Z}) = \mathbb{Z}_2 * \mathbb{Z}_3$.
}
\end{figure}

\paragraph{}
None of the $p$-adic group $G=\mathbb{Z}_p$ {\bf with the p-adic metric} are natural metric spaces.
For any metric which is invariant under $\mathbb{Z}_2$ however, there is a different
non-abelian group structure which still preserve the metric via right translations. 
In the case $p=2$, the construction is similar to the $\mathbb{Z} \to D_{\infty}$ transition:
one can see $\mathbb{Z}_2$ as the boundary of an infinite rooted Bethe lattice tree.
\footnote{This picture is often invoked as a more intuitive picture of the abstract 
{\bf projective limit} construction. A non-standard analysis trained person would 
take $n$ non-standard large and deal the sphere of distance $n$ to the origin 
of the root. It is a finite group but every element in the Cantor set is 
infinitesimally close to an element in that finite group.} 
If the generator $a$ flips the two main sub-trees and the generator
$b$ flips the sub-trees and additional induces the group addition by $1$ 
on the right tree and group addition by $-1$ on the left tree, then $a^2=b^2=1$
and $ab=T$ is the group addition by $1$ on half of the space and subtraction
by $1$ on the other half. The construction for $\mathbb{Z}_p$ is similar. 
There is no other metric on $\mathbb{Z}_p$ which forces the group structure,
because any such $G$-invariant metric is also invariant under a non-abelian
modification. 

\paragraph{}
The case of the Pauli group $P(1)$ appeared in \cite{MoebiusKantor} (where we still were
under the spell of the old definition). The result mentioned there that $P(1)$ is natural 
remains true, but its proof needed to be revisited and the general structure theorem covers it.

The {\bf Pauli group} $P(1)$ of order $16$ has the Moebius-Kantor graph MK as
Cayley graph. If $X,Y,Z$ are the Pauli matrices, then $P(1)$ is the group with elements
a b, where a is in $\{1,-1,i,-i\}$ and b is in $\{1,X,Y,Z\}$.
Its automorphism group ${\rm Aut}(MK)$ is the Tucker group.
We can write $P(1) = E \rtimes C_2$, where $E = C_4 \times C_2$.
It is remarkable that also the semi-dihedral $SD_{16}$ and the dihedral group $D_{16}$ 
have MK as a Cayley graph. There are 4 nonabelian groups of order $16$. They are
$P(1),D_{16},SD_{16}$ and $Q_{16}=Dic_4$ and only $Dic_4$ is non-natural.
The Pauli group $P(1)$ is interesting as it is one of the small examples in the GRR 
classification theorem.

{\bf The Pauli group $P(1)$ is a natural group.}

\begin{proof}
As we use right translations, pick the Cayley convention $(g,sg)$ with
generators $S \in \{ X,Y,Z\}$. As usual, think of the generators also as edge colors. 
Right translation $R_h$ maps the edge $(g,gh)$ to the edge $(gh,sgh)$. Assume the three edge lengths are 
$3,4,5$. We break so the symmetry and the automorphism group of the Moebius Kantor graph (which is the 
Tucker group) becomes a smaller group: it will be $P(1)$. 
The weighted geodesic metric $d$ on $MK$ recognizes which pairs are edges
and to which generator they belong. Hence, every isometry preserves each of the three edge colors.
Every isometry that fixes the identity must be the identity so that the stabilizer of 1 in the isometry group
is trivial. The 16 right translations are isometries and act transitively so that the right translation
group is P(1). Any other group structure $H$ on the 16 points for which
right translations are isometries must be a subgroup of P(1) and since its
cardinality is 16, it must be P(1).
\end{proof}

\paragraph{}
This example illustrates also that if a group $G$ is presented with a finite generating set $S$
consisting of {\bf involutions}, then it is natural and the weighted Cayley graph is a natural metric space.
The argument implies for example that the symmetric group $S_n$ is a natural group, as it is
generated by $n-1$ involutions. And $A_n$ is involution generated if $n \geq 5$.
The group $A_2$ is the trivial group which is natural. 
The group $A_3=C_3$ is natural because it its order is prime. 
The case $A_4$ needs to be treated separately. Write 
$A_4 = \langle a,b | a^2=b^3=(ab)^3=1 \rangle$. Its Cayley graph with generator set $S=\{a,b,b^{-1} \}$
is the truncated tetrahedron graph. It is a natural graph. Its isometry group is $S_4$ and
$S_4$ has a unique subgroup of order $12$. So, also $A_4$ is natural. \\

{\bf The groups $S_n$ are all natural. The groups $A_n$ are all natural. }

\paragraph{}
The special case $A_4$ illustrates that if a metric space $X$ with $n$ elements admits a
regular subgroup $G$ of order $n$ and every subgroup of ${\rm Isom}(X)$ is isomorphic to $G$,
then $X$ is a natural metric space and $G$ is a natural group.
The article \cite{MoebiusKantor} from 2026 still used the stronger definition. Its 
proof actually works for the right translation but not for strong naturality condition. 
It needed a revisit. \\

\paragraph{}
The case of all groups of order 16 illustrates well the general
naturality classification theorem for finite groups. 
{\bf Among the 14 groups of order 16, there are 7 that are natural
and 7 that are non-natural. } There are $9$ non-abelian groups of order $16$. 
Three of them are generalized dicyclic and so non-natural.
The dihedral $D_{16}$, the semi-dihedral $SD_{16}$ and the Pauli group $P(1)$
are natural. Then there is $M_{16}$, $D_8 \times C_2$ and $C_2^2 \rtimes C_4$ 
which are natural. 
The dicyclic $Q_{16}={\rm Dic_4}$, the generalized dicylic 
$Q_8 \times C_2$ and the generalized dicyclic $C_4 \rtimes_{-1} C_4$ are not natural.
The group $Q_8 \times C_2$ is $Dic(C_4 \times C_2,i^2)$ if we understand
$i$ as part of $Q_8$. The group $C_4 \rtimes_{-1} C_4$ is 
$Dic(C_4 \times C_2,b^2)$ if we write $G=\langle a,b | a^4=b^4=1, bab^{-1} = a^{-1} \rangle$.
Among the 5 abelian groups of order $16$, only $C_2^4$
is natural. The others four, $C_4 \times C_4, C_2 \times C_8$ or $C_4
\times C_2^2$ are all non-natural. This illustrates a difference
between graphical regular representation an naturalness. The group $P(1)$
does not have a GRR and still is natural. One key distinction between
having a GRR and naturalness is that an ordinary uncolored Cayley graph 
may retain unavoidable-looking extra symmetries. A metric can distinguish
different inverse pairs by imposing different distances for them.

\paragraph{}
An interesting question is to compare the number $Nat(n)$ of natural groups of order
$n$ with the number $NonNat(n)$ of non-natural groups of order $n$. We have infinitely
often equality and both directions also appear infinitely often
\begin{itemize}
 \item $Nat(2p)=Nat(2p)=1$
 \item $Nat(4p) \in \{1,2\} < NonNat(4p)=3$ 
 \item $Nat(4p^2)=NonNat(4p^2)=6$, $p \geq 7, p=4k+3$.
 \item $Nat(2n+1) > NonNat(2n+1)=0$ for all $n$.
 \item $Nat(2p^2)=3 > NonNat(2p^2=2$ for all odd primes
\end{itemize}
For example, $n=4p$, thee are the 2 Abelian $C_{4p}$ and $C_2 \times C_{2p}$ and
the $Dic_p$ which are non-natural so that $NonNat(4p)=3$. 
When $p=4k+3$, there is only one remaining group $C_2 \times D_{2p}$ and that is regular.
When $p=4k+1$, there is an additional natural group. For $n=2^k$, the number 
of non-natural groups are extremely sparse in comparison to natural ones.

\begin{comment}
A={{1,1,0}, {2,1,0}, {3,1,0}, {4,1,1}, {5,1,0}, {6,1,1}, {7,1,0}, {8,2,3}, {9,2,0}, {10,1,1},
{11,1,0}, {12,2,3}, {13,1,0}, {14,1,1}, {15,1,0}, {16,7,7}, {17,1,0}, {18,3,2}, {19,1,0}, {20,2,3},
{21,2,0}, {22,1,1}, {23,1,0}, {24,10,5}, {25,2,0}, {26,1,1}, {27,5,0}, {28,1,3}, {29,1,0}, {30,3,1},
{31,1,0}, {32,39,12}, {33,1,0}, {34,1,1}, {35,1,0}, {36,8,6}, {37,1,0}, {38,1,1}, {39,2,0}, {40,9,5},
{41,1,0}, {42,5,1}, {43,1,0}, {44,1,3}, {45,2,0}, {46,1,1}, {47,1,0}, {48,43,9}, {49,2,0}, {50,3,2},
{51,1,0}, {52,2,3}, {53,1,0}, {54,12,3}, {55,2,0}, {56,8,5}, {57,2,0}, {58,1,1}, {59,1,0}, {60,10,3},
{61,1,0}, {62,1,1}, {63,4,0}, {64,246,21}, {65,1,0}, {66,3,1}, {67,1,0}, {68,2,3}, {69,1,0}, {70,3,1},
{71,1,0}, {72,40,10}, {73,1,0}, {74,1,1}, {75,3,0}, {76,1,3}, {77,1,0}, {78,5,1}, {79,1,0}, {80,43,9},
{81,15,0}, {82,1,1}, {83,1,0}, {84,12,3}, {85,1,0}, {86,1,1}, {87,1,0}, {88,7,5}, {89,1,0}, {90,8,2},
{91,1,0}, {92,1,3}, {93,2,0}, {94,1,1}, {95,1,0}, {96,217,14}, {97,1,0}, {98,3,2}, {99,2,0}, {100,10,6}}
f[{x_,y_,z_}]:={y,z}; B=Map[f,A]; ListPlot[Transpose[A], Joined->True,PlotRange->All] 
\end{comment}

\paragraph{}
This looked like a nice research problem. Some of it can be delegated. Here is the answer
to my question "In the case $n=2^k$ there are lots of groups and the number of natural groups vastly outnumbers the number of 
non-natural ones. Like for $n=64$, where there are already $246$ natural and only $21$ non-natural groups. 
Can one say something about this case? In general, the number of groups of a fixed order is not known 
in closed form, but for $n=2^k$ it could be possible to say something."
The language model came up by itself with the theorem. It first explored the question experimentally for small
cases, looked up the total number \cite{Burrell2022} of groups of the order $2^{10}$ and concluded that there
are $49487367182$ natural groups of order $2^{10}$, because
there are only $107$ non-natural ones, where $66$ are generalized dicyclic and $31$ are non-natural abelian. 
Then it formulated by itself and proved that naturality is overwhelmingly generic among 2-groups.

\begin{thm}[GPT 5.6 Sol]
For $k \to \infty$, almost every group of order $2^k$ is natural.
\end{thm}
\begin{proof} 
There are $p(k)$ abelian groups of order $2^k$ among which exactly one of them is abelian.
The number of generalized dicyclic groups of order $k$ is (after some simplification) 
$\sum_{j=0}^{k-2} p(j)-1$ so that $NonNat(2^k)=p(k) + \sum_{j=0}^{k-2} p(j)-2$.
Using the Hardy-Ramanujan estimate $p(k)=e^{\pi \sqrt{2k/3}}/(4 k \sqrt{3})$ one can 
get $\sum_{j \leq k} p(j) = \sqrt{5k} p(k)/\pi = exp(O(\sqrt{k}))$.
The Higman-Sims theorem \cite{Sautoy2000} for the total number $g_k$ of 2-groups gives 
$\log_2(g_k) \sim 2 k^3/27$. Therefore, $NonNat(2^k)/g_k \to 0$ and $Nat(2^k)/g_k \to 1$. 
\end{proof} 

\section{Dynamical aspects} 

\paragraph{}
The following discussion is more relevant to {\bf strongly natural}. The reason
is that a natural group $(G,*)$, when upgraded to a metric space $(G,*,d)$ is not
always a topological group. For a topological group, both left and right translations need
to be continuous.

\paragraph{}
Metric groups $G$ are examples of {\bf topological groups}
\cite{MontomeryZippin,PontryaginGroups,Husain}, objects located in the 
heart of harmonic analysis, physics or operator algebras. Examples are 
{\bf compact Lie groups} $G$ which are special as 
they feature an even bi-invariant Riemannian metric $d$.
A finite group $G$ with a selected set of generators 
$S = \{ g_1,g_1^{-1}, \dots ,g_k,g_k^{-1} \}$
satisfying $S^{-1}=S$ define an {\bf undirected Cayley graph} with {\bf geodesic 
distance} as metric, where the distance between two group elements $d(g,h)$ is the 
minimal {\bf word length} of $g^{-1}h$. 

\paragraph{}
Assume $G$ is a natural group. It has both a metric as well as an algebraic structure. 
The action of $G$ by right translations can be through of as {\bf ``time"} acting on the 
{\bf metric space} $(G,d)$ by right translations. The metric structure is considered 
{\bf ``space"}. In this context, time is associated with the algebraic structure coming from right translation
and space with the topological structure. The right translations $R_g(x)=x g$ is {\bf transitive} by definition.
For a natural metric space $(X,d)$, the metric alone produces a {\bf natural group of 
symmetries on its own metric space} without any need to define the group operation. 
The {\bf group multiplication table} is determined
up to isomorphisms from the topology. One can now define notions like the 
{\bf periodic part of the group} and the complement, the {\bf aperiodic part of $G$}:
the periodic part of $G$ is the set of $x \in G$ for which there exists an $n>0$ such that 
$x^n=1$. The aperiodic part is the complement in $G$. 

\paragraph{}
In the case of a {\bf compact natural metric space} the topological dynamical system 
can be upgraded to a canonical {\bf measure theoretical dynamical system}. 
There is {\bf uniquely ergodicity}: there is a unique $G$-invariant probability 
measure $m$ on the {\bf Borel $\sigma$-algebra} of $G$. 
It is the normalized {\bf Haar measure} and invariant under
automorphisms of the probability space given by the right 
multiplication with an element in $G$. 
Natural groups define {\bf uniquely ergodic transformation groups} for which the metric alone 
determines the arithmetic as well as all the expectations of random variables. 

\paragraph{}
Note that the existence of a Haar measure does not even need naturality. It follows from a 
transitive action of a group by isometries. Since the right translation group applied to 
the identity reaches the entire group and so is transitive, there is a unique Haar measure. 

\paragraph{}
For a natural group $G$, the {\bf group of isometries} ${\rm Aut}(G,d)$ of the metric space 
$(G,d)$ defining $(G,+)$ must contain the group $G$. But it can be larger: 
in the case of the complete graph $(K_p,d)$ for example with the discrete metric $d(x,y) = 1$
for all $x \neq y$, the isometry group is the 
{\bf symmetric group} $S_p$. It is still natural if $p$ is prime. 
So, $K_7$ is an example of a natural graph, even so it is not a graphical regular 
representation of a group.
It would be interesting to know the complexity of deciding whether a 
finite group of order $n$ is strongly natural or not.

\paragraph{}
A natural compact metric space $(X,d)$ defines a 
dynamical system $(G,X,d,m)$, where $G=X$ acts on $X$
as isometries by right translation. 
It still preserves a unique {\bf Haar probability measure} $m$.
Because the action is {\bf transitive}, it is {\bf ergodic} (every $G$ invariant
Borel set has measure $0$ or $1$). It is even uniquely ergodic (there is only one 
invariant measure). Given a specific point $x$ in the
metric space $(G,d)$, we have right transformations $R_x: X \to X$
which are measure-preserving. 

\paragraph{}
Let us call a point $x \in X$
{\bf ergodic}, if $R_x$ is ergodic with respect to the invariant measure $m$.
This happens if and only if the subgroup generated by $x$ is dense in the entire group.
If the orbit of $R_x$ is dense, then by the {\bf Krylov-Bogolyubov theorem}, 
the Birkhoff functionals $f \to \limsup_n n^{-1} \sum_{k=1}^n f(T_x^n y)$ on the 
Banach space $C(G)$ of continuous functions on $G$ define an invariant measure, the
Haar measure on $G$. The metric space $X$ 
obviously splits into an ergodic and a non-ergodic part. For $\mathbb{T}=\mathbb{R}/\mathbb{Z}$ 
for example, all rational $x$ are non-ergodic, while the irrational $x$ are ergodic. 

\paragraph{}
As an other example is obtained by looking at the 
{\bf $k$-dimensional torus} $X=\mathbb{T}^k$ for which $(X,+)$ is a natural group 
as a compact connected Lie group. Almost all points $x \in X$ are ergodic. 
A point $x=(x_1,x_2, \dots x_k) \in X$ is ergodic if the coordinates are 
{\bf rationally independent} meaning that there is no non-zero integer vector 
$n=(n_1,n_2, \dots, n_k)$ such that 
$n \cdot x = n_1 x_1 + n_2 x_2 + \cdots + n_k x_k = m$ is an integer. 
This assures transitivity. Ergodicity can be establihed by Fourier theory. 
Lebesgue almost all points $x$ are ergodic because the set of non-ergodic $x$ 
is a countable union of sub-manifolds in $\mathbb{T}^k$ which all have 
measure zero so that also the union has zero measure. 
In the case $\mathbb{T}^2$ for example, the set of non-ergodic
points is a union of projections from $\mathbb{R}^2 \to \mathbb{T}^2$ 
of countably many lines $a x + b y = c$, where $a,b,c$ 
are integers and $(x,y) \in \mathbb{T}^2 = \mathbb{R}^2/\mathbb{Z}^2$. 

\paragraph{}
If $(V,E)$ is a natural finite simple graph with vertex set $V$, 
then the Haar probability measure is the {\bf uniform counting measure} which assigns 
to every point the constant weight $1/|V|$. Once we have fixed a specific implementation 
of the by assumption unique group structure, every vertex $v \in V$ has an associated
automorphism $T_x=R_x$ given by right translation. The {\bf Lefschetz number} $\chi(G,T_x)$
of such a transformation is defined as the {\bf super trace} of $T_x$ on the
cohomology groups $H^k(G) = {\rm ker}(L_k)$, where $L_k$ are the blocks in the 
{\bf Hodge Laplacian} $L = D^2=(d+d^*)^2 = d d^* + d^* d$. We write this as
$$ \chi(G,T_x) = \sum_{k=0}^{d} (-1)^k {\rm tr}(T_x | H^k(G)) \; . $$ 
Note that by Hodge theory, the vector spaces $H^k(G)$ are just the kernels of 
concrete matrices $L_k$ acting on $k$-forms (functions on all $k$-dimensional 
simplices of $(V,E)$. 
The {\bf Lefschetz fixed point theorem} for graphs \cite{brouwergraph} 
\footnote{It holds of course for any finite abstract simplicial complex. 
The Whitney complex of a graph is just one example and the most intuitive one.}
assures then that $\chi(G,T_x)= \sum_{y \in Fix(T_x)} i_T(y)$, where ${\rm Fix}(T_x)$
is the {\bf fixed point set} consisting of complete sub-graphs of $G$ which are 
fixed by $T_x$.

\paragraph{}
We also know that the average Lefschetz number over 
all automorphisms of a graph is $1$: 
$$ \frac{1}{|V|} \sum_{x \in V} \chi(G,T_x) = \chi(G/G)=1 \; . $$
This is essentially a Riemann-Hurwitz formula (see \cite{brouwergraph}).
We have now a necessary algebraic condition which allows to probe whether 
a graph is natural: 

\begin{propo}
Let $(V,E)$ be a natural graph. The sum over all
Lefschetz numbers over all transformations $T_x$ is equal to $|V|$.
\end{propo}

\paragraph{}
The Lefschetz fixed point theorem for graphs also allows to get 
a grip on counting simplices which are fixed
by the automorphisms. These fixed simplices become points in the 
Barycentric refinement and correspond to fixed points of homeomorphisms
in continuum frame works. 
For a complete subgraph $y$ that is fixed by $T_x$ and a vertex $x$,
let $i_x(y)$ be the {\bf index} of $T_x$, defined as
$$ i_x(y) = (-1)^{{\rm dim}(x)} {\rm sign}(T_x|y) \; . $$

\begin{coro}
For a natural graph, $\sum_{x \in G} \sum_{y \in {\rm Fix}(T_x)} i_x(y) = |G|$.
\end{coro}

\paragraph{}
For the cyclic graph $C_n$ for example, where the isometry group is is 
the natural group $G=D_n$, every non-zero group element $x$ defines a translation 
which has zero index and the identity has total index $0$, the Euler characteristic of the graph. 
There are always $m$ reflections which have index $2$ and $m$ reflections which have index $0$. 
We see that the total index $\sum_{x} \sum_{y \in {\rm Fix}(T_x)} i_x(y) = 2n = |D_n|$. 

\section{Graphical regular representations}

\paragraph{}
A simple criterion for a graph $G=(V,E)$ to be natural 
involves the {\bf automorphism group} ${\rm Aut}(G)$ of $G$. It is the group of all 
{\bf graph automorphisms} of $G$. An automorphism of a graph could also be defined as
a bijections on the vertex set $V$ which preserves the set $E$ of edges and the set 
$E^c$ of non-edges. 

\paragraph{}
The {\bf graphical regular representation}
of a group $G$ is a {\bf Cayley graph} $(V,E)$ with generator set $S \subset G$, for which 
the automorphism group of $(V,E)$ are the groups right-multiplication actions. 
A group $G$ is said to have a {\bf graphical regular representation}, if it is the 
automorphism group of a graph with a graphical regular representation.

\begin{propo}
If $(V,E)$ is a finite simple graph for which the automorphism group $G={\rm Aut}(V,E)$ has
order $|V|$, then the graph $(V,E)$, the group $(G,+)$ as well as the metric space
$(V,d)$ are all natural. 
\end{propo}

\begin{proof}
Because right translations in a natural group are graph isomorphisms of the Cayley graphs, 
a natural group $G$ must be a subgroup of the finite group ${\rm Aut}(V,E)$. 
Of course, there is only one group with $|V|$ elements which is a subgroup ${\rm Aut}(V,E)$ 
and which has $|{\rm Aut}(V,E)|=|V|$ elements: it is $G={\rm Aut}(V,E)$. 
In order to attach group elements to vertices of the graph, choose one of the vertices
and call it $0$. For any automorphism $T \in {\rm Aut}(V,E)$, attach it to the 
vertex $x=T(0)$. This pairs up $G={\rm Aut}(V,E)$ in a 1:1 manner with the vertex set $V$. 
\end{proof}

\paragraph{}
Almost equivalent is the case of finite metric spaces. The automorphism group now becomes
{\bf the group of isometries}:

\begin{propo}
Assume $(G,d)$ is a finite metric space for which the group of right translation isometries 
has order $|G|$, then $(G,d)$ is natural.
\end{propo}

\begin{proof}
Any group compatible with the metric in the given way must be a subgroup of the group of isometries
and if it is the group of symmetries, it is uniquely determined. We can build up the
group structure by labeling any point as $0$, then take an isometry $T$ in the isometry
group and assign the isometry group element $T$ to the point $T(0)$ in the metric space.
\end{proof}

\paragraph{}
One can also reverse this. We see that most graphs are non-natural because most graphs
lack symmetry. Especially, every graph with trivial automorphism group (meaning 
there is no symmetry except the identity) is non-natural if it 
has more than two vertices. For example, the linear graph $L_2$ of diameter $2$ with 
$3$ vertices can not be natural because its automorphism group is $C_2=\mathbb{Z}/(2 \mathbb{Z})$. 
This is not large enough to host a group structure with $3$ group elements. 

\begin{propo}
If $|V|$ does not divide the order of the automorphism group of a finite
simple graph $(V,E)$, then $(V,E)$ is not natural.
\end{propo}
\begin{proof}
By the {\bf Lagrange theorem} in group theory, the order of a subgroup of ${\rm Aut}(V,E)$ 
must divide the order $|{\rm Aut}(V,E)|$ of ${\rm Aut}(V,E)$. For a natural graph $(V,E)$, the 
order of the induced group is $|V|$. 
\end{proof}

\paragraph{}
While the topic of natural groups has relations with 
graphical regular representations of a group, it is
different. Graphical regular representations of a finite group render
the group natural but there are natural graphs which
do not have a graphical regular representation and there
are natural groups which do not come from graphs alone but
need weighted graphs. A general, finite metric space can always
be described by a weighted Cayley graph, as we are allowed to
take $S=G$ as the generating set. 

\begin{propo}
If a finite Cayley graph is a graphical regular representation, then it is natural. 
If a finite group has a graphical regular representation then it is natural.
\end{propo}

\begin{proof}
This follows from the definitions as the Cayley graph defines then
a metric space $(V,d)$ on which the group $G$ has the same 
order $|V|$ which is a subgroup of the automorphism group. 
\end{proof}

\paragraph{}
eany graphs are natural even so they have no graphical regular representation.
For natural graphs, the automorphism group can be much larger than the number
of vertices of the graph. For example, every complete graph $K_p$ with prime $p$
is natural despite that it is far from being a graphical regular representation. 
In the case of $(V,E)=K_p$, we have $|{\rm Aut}(K_p)|=p!$ which is for $p>2$ always larger
than $|V|$. Still, there is a unique group which can be planted on $K_p$ because
there is only one group of prime order $p$. 

\paragraph{}
If $G$ has a graphical regular representation, then only the
identity element fixes the set $S$ of generators. But even if
$G$ induces non-trivial symmetries on $S$, it is still possible
that the graph is natural. We just put different distances on the
different generator edges which breaks the symmetries and so that 
$G$ induces no symmetry any more on $S$. Still, we can use results
from the theory of graphical regular representation, \cite{LiSim2000}
for example gets a result implying that all 
{\bf metacyclic $p$-groups} are natural. 

\paragraph{}
As the above proposition shows, the concept of graphical regular
representation comes close to the notion of ``naturalness". There is a difference
because a group can be natural without coming from a bare graph. It can be natural 
with respect to a metric space like for example a weighted Cayley graph. 

\paragraph{}
The concept of graphical regular representations started with
Sabidussi \cite{Sabidussi1964}, who calls a graph 
{\bf vertex transitive} if the automorphism group is transitive. At that time, Cayley
graphs were still called {\bf Cayley color-groups} as one can imagine the different generators
assigned different colors. Examples of follow up work is
Nowitz in 1968 and Imrich, Watkins, Bannai \cite{Bannai1977} in the 1970ies \cite{Watkins1974}.
Modern sources are \cite{Godsil2015,LiSim2000}.

\section{Graph operations}

\paragraph{}
The {\bf graph complement} of $(V,E)$ is the graph $(V,E^c)$, where
$E^c$ is the complement of the set of all possible edges on the
complete graph $K_{|V|}$ with $|V|$ elements. It is here convenient to 
extend graph isometry also to the {\bf disconnected case}, where $d(x,y) = \infty$
rewrites that $x,y$ are in different connectivity components.
With this assumption, every permutation of the vertex set $V$ is a graph isometry of the
complete graph $K_n$ as well as of the graph complement $P_n=K_n^c$,
the {\bf point graph} with $n$ vertices and no edges.

\paragraph{}
The graph complement of a Cayley graph $\Gamma(G,S)$ 
is sometimes a Cayley graph again with $\Gamma(G,S^c)$ like if $S^c=G-S-\{0\}$ generates the
same graph. The group $G=C_7$ can be generated by $a(x)=x+1,a^{-1}=x-1$
leading to the Cayley graph $C_7$, the cyclic graph with $7$
elements. Its graph complement $C_y^c$ which is a discrete Moebius
strip generated by $a(x)=x+2,a^{-1}(x)=x-2$
and $b(x) =x+3,b^{-1}(x)=x-3$ on $\mathbb{Z}/(7 \mathbb{Z})$. 
Graph complements of cyclic graphs have a rich structure \cite{GraphComplements}. 
Let us not use the extended definition for metric spaces allowing $d(x,y)=\infty$
and restrict to graphs for which both graph and complement are connected.

\begin{propo}
If $(V,E)$ is a connected natural finite simple graph such that its graph 
complement is connected, then his graph complement $(V,E^c)$ is natural too.
\end{propo}

\begin{proof}
Every isometry of a finite simple graph equipped
with the geodesic metric necessarily maps edges into edges and non-edges to non-edges.
Graph automorphism and isometries are the same. Now we can use the fact that the group of
graph automorphisms of a graph $G$ and the group of graph automorphisms
of the graph complement $G^c$ is the same. If $G$ is natural, there is an up to isomorphism
unique group structure on $G$ so that all group operations are isometries. This group structure
is also possible on the graph complement. Any different, non-isomorphic group structure on 
the complement would produce on on $(V,E)$.
\end{proof}

\paragraph{}
Let us comment on the disjoint union of two metric spaces which is some sort
of problematic notion as it requires to define the distance between two separate points
to be infinite. In the case when both metric spaces have finite diameter, we can postulate
a finite distance between disconnected parts. 
The {\bf disjoint union} of two compact metric spaces is a metric space again, if we
allow the distance function to be extended. 
But it is important to note that the disjoint union of two natural
metric spaces $(X_1,d_1) \cup (X_2,d_2)$ is not necessarily natural in the strict sense
simply because we would need to extend the distance function to take infinite values. 
The disjoint union of two copies of $\mathbb{R}$ for example is not natural. Place $0$
on one branch, then look at the point $a$ of minimal distance to $0$ on the second branch.
Now because of the isometry assumption $a=-a$ and every point $x$ on one branch has
a mirror point $x+a$. This allows for the group structure $\mathbb{R} \times C_2$ 
as well as a semi-direct non-abelian group structure $\mathbb{R} \rtimes \mathbb{Z}_2$
where every element on the mirror branch is a reflection. 

\paragraph{}
Also the {\bf join} of two natural graphs is in general not natural.
While we could put a weighted distance on the cyclic graph $C_4$ such that it becomes
natural, just produce a rectangular shape. The group then is the
{\bf Klein four group} $C_2 \times C_2 = D_2$ which is natural. 
The graph $K_2 \oplus K_2 = K_4$ is not natural (as a graph) while $K_2$ is natural
and $K_2 \oplus K_2$ carries a weighted metric that is natural. 

\paragraph{}
The {\bf discrete 3-sphere graph} $C_4 \oplus C_4$ is the 16-cell. 
and smallest discrete $3$-sphere and has $8$ vertices. 
It is {\bf not a natural graph} as it admits both the quaternion group as well
as the elementary abelian group $C_2^3$. 
Also the 2-sphere, the {\bf octahedron graph} $C_4 \oplus C_2$, which is a $2$-sphere with $6$ vertices,
is not natural. There are two groups of order $6$. The cyclic group
$C_6=\mathbb{Z}/(6\mathbb{Z})$ and the symmetric group $S_3=D_3$. The automorphism group of the 
octahedron is $S_4 \times C_2$ which has 48 elements and this admits both $C_6 \sim C_3 \times C_2$
as well as $S_3$. If we look at the graph as a metric space and break the symmetry of the metric
by changing the length of a single triangle, then we have
only the $S_3$ symmetry left so that $S_3$ is natural. 

\begin{figure}[!htpb]
\scalebox{0.4}{\includegraphics{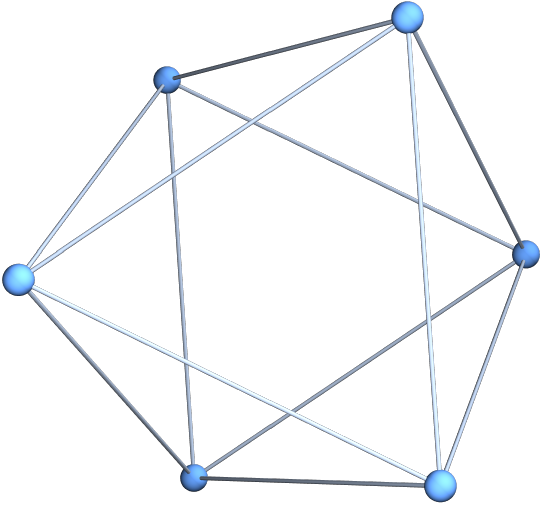}}
\label{Figure 1}
\caption{
The octahedron graph is not natural as it admits two different
group structures. One can change the metric however to produce a natural
metric space: shrink the lengths of a single triangle without
changing any other lengths. This produces a metric space which is natural
and induces the group $S_3$. We knew already that $S_3=D_3$ is natural
because it is a dihedral group, and because dihedral groups are natural. 
}
\end{figure}

\paragraph{}
There are various graph multiplications. For all of them, the vertex
set of the product is the {\bf Cartesian product} $V(G) \times V(H)$ of the
vertex sets of $G$ and $H$. Among abelian multiplications are 
the {\bf weak product}, the {\bf strong product} (=Shannon product \cite{Shannon1956}), 
the {\bf large product} (the dual of the strong product \cite{Sabidussi}) and the 
{\bf tensor product} in which $E(G) \times E(H)$ is the edge set. A non-abelian
product is the {\bf lexicographic product}.
A general question is to determine which graph products preserve natural graphs. 
An other question is which graph products allow for a new metric on the product so that
it becomes a natural metric space. 

\paragraph{}
Each graph multiplication can be extended to weighted graphs and produce
again a weighted graph. An edge $( (a,b), (c,d) )$ in the product has then
the weight $d(a,c) + d(b,d)$. Any finite natural group defines weighted
natural graphs and by adapting the weights one can assure that the 
{\bf direct product} of two groups is natural. We will come to that. 

\paragraph{}
For the subclass of Cayley graphs (in particular for graphs which must be
vertex transitive), there is an other product, the {\bf Zig-Zag product}. 
It is associated with the {\bf semi-direct product} of the corresponding groups and 
when weighted can be used to see that the semi-direct product of two finite
natural groups is natural. 

\paragraph{}
There is also a {\bf non-abelian} multiplication,
the {\bf lexicographic product of graphs} $G \cdot H$ 
(introduced by Hausdorff in 1914 and further studied by Harary and Sabidussi). 
The edges in the later are $( (a,b), (c,d))$ if $(a,c) \in E(G)$ or $a=c$ and
$(b,d) \in E(H)$. The non-abelian lexicographic product is interesting because it is
{\bf self dual} in that the product of the graph complements is the complement of the
product and that the {\bf factorization problem} for graphs is complexity-wise equivalent
to the graph isomorphism problem. This means that 
testing whether a graph of size $n$ is a lexicographic product of two smaller graphs
is polynomially equivalent to decide whether two graphs of size $n$ are isomorphic. 
It has been asked in \cite{Zerovnik}, whether 
prime factors of the strong product can be computed in polynomial time. 

\paragraph{}
The {\bf lexicographic product} of $P_n$ with a graph $G$ is the $n$-fold copy of $G$.
However, the lexicographic product of $G$ with $P_n$ is connected. 
The graphs $K_n \cdot P_2$ are $(n-1)$-spheres
by induction because $K_n \cdot P_2$
is the suspension of $K_{n-1} \cdot P_2$. The Betti vectors are $b=(1,0,0, \dots,1)$
leading to Euler characteristic $0$ or $2$ depending on parity. 
%Table[s=NormalizeGraph[LexicographicProduct[CompleteGraph[k],PointGraph[2]]];Betti[s],{k,6}]
Interesting is the case $K_n \cdot P_m$ which leads to $n$-dimensional graphs
with Betti vectors $(1,0,0, \dots,(m-1)^n)$. This produces explicit examples
of graphs, where the Euler characteristic grows exponentially with the 
dimension and the vertex size. The $f$-vector $f=(f_0,f_1,f_2, \dots, f_d)$
of $K_n \cdot P_m$ is 
$f_k = B(n,k+1) m^k$ and the $f$-function 
$f(t)=1+f_0 t + f_1 t^2 + ... f_{d} t^{d+1}$ is $(1+m x)^n$. 

\begin{propo}
The Lexicographic products $A \cdot B$ and $B \cdot A$ of two finite simple connected natural 
graphs can be weighted to become a natural metric space. 
\end{propo}

\begin{proof}
Again just take the product metric where the distances in the first
coordinate and the second coordinate are different. We need connected in order to have a natural
metric space in the product. We can not conclude that spheres $K_n \cdot P_2$ are natural for
example. 
\end{proof}

\begin{comment}
n = 5; m = 11; s1=CompleteGraph[n];
s = NormalizeGraph[ LexicographicProduct[s1, PointGraph[m]]];
a = Fvector[CompleteGraph[n]]; 
Fvector[s]==Table[a[[k]] m^k, {k, Length[a]}]
n = 7; m = 13; s1 = CompleteGraph[n];
s = NormalizeGraph[LexicographicProduct[s1, PointGraph[m]]];
Factor[Ffunction[s, x]]
\end{comment}

\begin{figure}[!htpb]
\scalebox{0.4}{\includegraphics{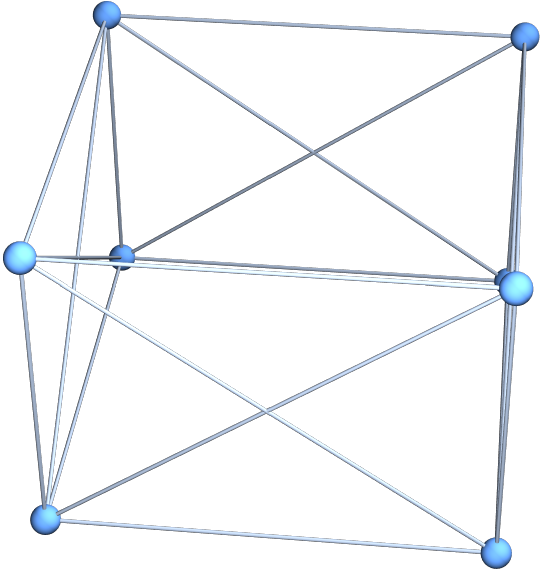}}
\scalebox{0.4}{\includegraphics{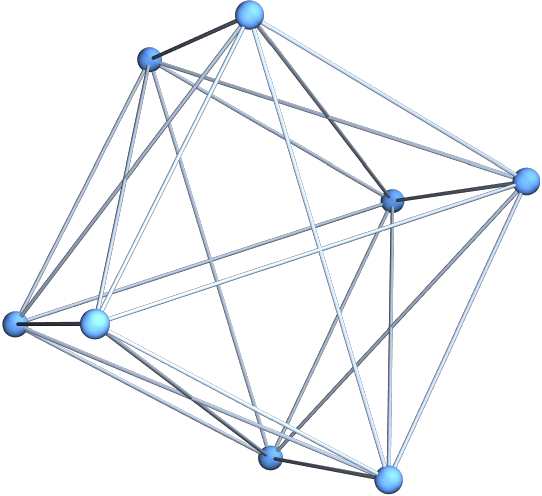}}
\label{Figure 1}
\caption{
The lexicographic product $C_4 \cdot K_2$ of $C_4$ with $K_2$
and the lexicographic product $K_2 \cdot C_4$ of $K_2$ with $C_4$
}
\end{figure}

\section{Completions}

\paragraph{}
The group $\mathbb{Z}$ is dense in the 3-adic group of integers $\mathbb{Z}_3$. While the former is not natural, 
the later is natural. The $2$-adic group of integers $\mathbb{Z}_2$ is not natural but
$\mathbb{Z}_p$ is natural for $p \neq 2$. The group $\mathbb{Q}$ is natural as it is a 
dense subgroup of the natural group $\mathbb{R}$. 
By Ostrowski's theorem, the only topological completion of the field of rational numbers $\mathbb{Q}$
are $\mathbb{Q}_2$ or $\mathbb{Q}_p$ with odd $p$ or $\mathbb{R}=\mathbb{Q}_0$.

\paragraph{}
{\bf Update 2026:} Unlike stated in the 2022 version of this paper, 
there is in general no relation between a group and its completion. The point is that 
naturality is a statement about regular actions which does not survive the completion
in general. If $G$ is dense in its completion $\overline{G}$, then every isometry 
extends uniquely to an isometry in $\overline{G}$. But an action that is free and transitive
(and so regular) on $G$ does not longer have to be free or transitive on its completion. 
The example of $\mathbb{Z} \subset \mathbb{Z}_p$ for odd prime $p$ shows this. 

There is a conceptional problem in the trying to complete a group topologically as a priori
we do not have a metric on the group. Sometimes we have a natural completion like given by 
Ostrowski's theorem. For a natural group, we have however a metric space and we can talk about
the {\bf topological completion of a natural group}. 

\begin{thm}
\bf $\mathbb{Q}_p$ is natural for every prime $p$.
\end{thm}

\begin{proof}
The reason is that $(\mathbb{Q}_p,+)$ and $(\mathbb{R},+)$ are both 
vector spaces over $\mathbb{Q}$. 
Since $\mathbb{Q}$ is countable and both $\mathbb{R}$ and $\mathbb{Q}_p$ 
have cardinality $c$ of the continuum, their dimensions as 
$\mathbb{Q}$-vector spaces are both $c$. 
Hence they are isomorphic as $\mathbb{Q}$-vector spaces, and 
therefore as abstract additive groups.
(Their Hamel dimension over the countable vector space $\mathbb{Q}$ 
is the continuum). This shows that $(\mathbb{Q}_p,+)$ and 
$(\mathbb{R},+)$ are isomorphic as abstract groups. 
If one of them is natural and $\mathbb{R}$ is, then 
all the others are. 
\end{proof}

\paragraph{}
{\bf Remarks.}\\
1) One can see an important distinction between $\mathbb{Z}_2$ and $\mathbb{Q}_2$. 
While the dyadic group of integers is not divisible by $2$, the group $\mathbb{Q}_2$ 
is divisible by $2$. The parity obstruction has disappeared. 
The group $\mathbb{Z}$ allowed to be generated by involutions that had no 
fixed points, producing the dihedral competitor. This was not possible in the 
case $\mathbb{Q}_2$. \\
2) The statement that the abstract group $\mathbb{Q}_p$ is natural is not 
equivalent to the question whether the $p$-adic metric $d_p(x,y) = |x-y|_p$ is 
a natural metric space. 

\paragraph{}
The situation is actually very interesting but would go too far. I asked Chat GPT 5.6 about
an example where the completion of a natural group for which the naturalizing metric is 
biinvariant and for which the completion is not natural. There are indeed examples.
$G=(\mathbb{Q}_2/\mathbb{Z}_2)^{(\mathbb{N})}$ is natural. 
The metric $d$ can be chosen bi-invariant and naturalizing, but 
$\overline{G}=H,H$ is not natural. 
See \cite{NaturalityCanBeLostUnderMetricCompletion}.

\section{Direct products}

\paragraph{}
The {\bf Shannon product} \footnote{Also called strong product} \cite{Shannon1956}
of two arbitrary graphs $G,H$ is the graph $G*H$ in which the Cartesian product
is the vertex set and where two pairs $(a,b), (c,d)$ are connected if the projection on both 
factors is either an edge or vertex. There are three possibilities therefore, connections
of the form $((a,b),(a,c))$ with $(b,c)$ being an edge in $H$. Or then, a connection is of the form 
$((a,c),(b,c))$ where $(a,b)$ is an edge in $G$. Or then a connection is of the form
$((a,b),(c,d))$ where $(a,c)$ is an edge in $G$ and $(b,d)$ is an edge in $H$. 

\paragraph{}
In the special case of a product of two Cayley graphs, we get again a Cayley graph. 
The Shannon product of the Cayley graphs $\Gamma(G,S_G), \Gamma(H,S_H)$
is the Cayley graph $\Gamma(G \times H, S_G S_H \cup S_G \cup S_H)$. If the
generator set is {\bf enhanced} by adding $1$ $\overline{S}_G=S_G \cup \{1\}$ then
$\overline{S}_G \overline{S}_H$ is an enhanced generating set for the Cayley graph.
The enhanced generator set is only used for the multiplication. It is not
that the Cayley graphs get equipped with a self-loop so that Cayley graphs remains
an undirected simple graph.

\paragraph{}
The Shannon product also led to the notion of Shannon capacity \cite{Shannon1956} 
(see \cite{ComplexesGraphsProductsShannonCapacity} for a more expository account). 
It produces a natural arithmetic \cite{ArithmeticGraphs,RemarksArithmeticGraphs}. 
It is compatible with curvature \cite{GraphProducts} and cohomology, also when 
looking at Lefschetz numbers.

\begin{propo}
The Cayley graph $\Gamma(G \times H, \overline{S}_G \overline{S}_H)$ is
the Shannon product of the Cayley graphs
$\Gamma(G,\overline{S}_G), \Gamma(H,\overline{S}_H)$.
\end{propo}
\begin{proof}
Given two points $(a,b), (c,d) = (xa,yb)$ where $x \in \overline{S}_G$ and $y \in \overline{S}_H$. 
They are connected if either $(a,c)$ is connected in $G$ and $(b,d)$ is connected in $H$ or then 
$a=c$ and $(b,d)$ is connected in $H$ or $b=d$ and $(a,c)$ are connected in $G$. This is exactly 
what the Shannon product does.
\end{proof} 

\paragraph{}
The {\bf direct product} $G \times H$ of two groups $G,H$ is a special case of the 
{\bf semi-direct product} $G \rtimes H$. In the direct product case, 
the action of the base on the fibers is trivial.
The direct product $G_1 \times G_2$ of two natural finite groups is not natural in general.
This generalizes to non-necessarily finite groups and to
{\bf semi-direct products} $G=H \rtimes K$ which is a group 
in which $H$ is a normal subgroup. Think of $K$ as the base and 
the group $H$ as the fiber. The base elements $k \in K$ now introduce 
group automorphisms $\phi_k: H \to H$ on the fibers and 
$(k_1,h_1) * (k_2,h_2) = (k_1 k_2, h_1 \phi_{k_1} k_2)$ is the group
operation on the semi-direct product. 
The semi-direct product of two natural groups is not necessarily natural. 

\paragraph{}
The prototype of a semi-direct product is the {\bf Euclidean
group} $G=\mathbb{R}^n \rtimes SO(n)$ generated by orthogonal transformations
(rotations) and translations. This is a vector bundle with
base manifold $SO(n)$ where at every point the fiber is a vector space 
$\mathbb{R}^n$. The operation is best seen when writing a group element
as $a: x \to Ax + b$ with $A \in O(n)$ and $b \in \mathbb{R}^n$. If
$b: x \to B x + c$ is an other group, then $a b(x) = A B x + Ac+b$, showing
how the base rotation group has acted on the fibers. Now, this Euclidean
group is already interesting as it is a non-abelian Lie group which does
not admit a bi-invariant metric.

\section{An other angle to Lie groups}

\paragraph{}
Here is an other approach to Lie groups.
There is a result of Niemiec of 2014 \cite{Niemiec2014} 
that establishes naturality of quite many groups including $\mathbb{R}^n$
and compact Lie groups.  {\bf The Euclidean group $E(n)=\mathbb{R}^n \rtimes SO(n)$ is natural for every $n \geq 1$. }
More generally using a result of Niemiec \cite{Niemiec2014} (applying in the non-abelian case) and 
Mazur Ulam \cite{MazurUlam1932,Vaisala2003,Nica2012} (in the abelian case). It again gives naturality. 
Even so the argument given before are more natural, it is still illustrative to see an other approach

\begin{thm}
Every connected non-abelian Lie group is natural.
\end{thm}

\begin{proof}
Niemiec uses left invariance but it can be turned into right invariance. 
He gives conditions for which a compatible proper right-invariant metric $d$ 
such that the full isometry group is exactly the right-translation group 
${\rm Iso}(G,\rho)=R(G)$. Every connected non-abelian Lie group has a trivial center
and thus has such a metric. 
The result which Niemiec reaches is stronger than naturality. If ${\rm Iso}(G,d)=R(G) \cong G$
and another group $(G,\circ)$ has all its right 
translations as $d$-isometries, then $R(G,\circ) \subset {\rm Iso}(G,\rho)$.
But $R(G,\circ)$ acts transitively on the underlying set. 
The original $R(G)$ acts regularly, so for every $x \in G$ there is exactly one element of 
$R(G)$ sending the identity to x. Transitivity therefore forces $R(G,\circ)=L(G)$.
Hence $(G,\circ) \cong G$. Thus $d$ is a naturalizing metric.
For $E(n)$, the assumptions work because the group is non-abelian for every $n\ge1$. Its center is trivial:
$Z(E(n))=\{1\}$ because element commuting with all translations must have orthogonal part $I$, 
and a translation commuting with every orthogonal transformation must have translation vector $0$. 
\end{proof} 

The following uses Mazur-Ulam. (I learned about the use of Mazur Ulam 
also from a preprint of Jakub Gismatullin. ).

\begin{thm}
Every connected abelian Lie group is natural.
\end{thm}

\begin{proof}
The non-abelian connected case has already been covered. In the abelian case, use that every
connected abelian Lie group has the form $G=\mathbb{R}^n \times \mathbb{T}^m \cong V/\Gamma$
where $V=\mathbb{R}^{n+m}$ and $\Gamma \cong \mathbb{Z}^m$ is a discrete lattice spanning $\mathbb{R}^m$.
As in the $\mathbb{R}^n$ case, chose a $l^4$ norm so that the linear isometry group is $\{\pm I\}$.
The Mazur-Ulam theorem assures that every surjective isometry of a finite-dimensional 
normed vector space is affine. Now put on $G=V/\Gamma$, the quotient metric $d(\bar{x},\bar{y})$
defined by $\inf_{\gamma \in \Gamma} |x-y+\gamma|$. It induces the usual Lie-group topology.
Because $\Gamma$ is discrete, the projection $V \to V/\Gamma$ is locally an isometry. 
Hence every isometry of $G$ lifts to an isometry of the universal cover $V$. By Mazur-Ulam,
the lift is affine $x \to Ax+b$. The linear part $A$ must be a linear isometry of the norm and 
preserve the deck group $\Gamma$. But the linear isometry group was chosen to be only
$\{\pm I\}$. Therefore every isometry of G has the form $\bar{x} \to \bar{x}+\bar{b}$
or $\bar{x} \to -\bar{x}+\bar{b}$. Thus ${\rm Isom}(G,d) = G\rtimes\{\pm 1\}$.
Now use the same argument as in the odd $p$-adic integer $\mathbb{Z}_p$ case.
Since $G$ is connected abelian, it is divisible; in particular $2G=G$. 
This implies that every reflection $S_b(x)=b-x$ has a fixed point.
A right-translation group of any competing group structure must act freely. 
Hence it cannot contain any reflection. 
It must lie entirely in the translation subgroup G. Since it is also transitive, 
it must be equal the entire translation subgroup.
Therefore the quotient metric naturalizes $G$:
The statement $\mathbb{R}^n \times \mathbb{T}^m$ is natural, has been proven.
\end{proof}

\paragraph{}
This second proof is still also elegant: the nonabelian connected case is supplied 
by Niemiec, while the abelian connected case boils down to Mazur-Ulam.
Things are related.  It illustrates the theme:

\begin{itemize}
\item We ask: Given a group $G$. Can a metric on $G$ 
force its group structure? 
\item Niemiec studied questions like:
given $G$. Can we construct a metric space, whose full isometry group is $G$?
\end{itemize}

He manages to do this on the underlying set $G$ itself, with the regular 
action of $G$. This immediately is relevant to naturality. If
${\rm Iso}(G,d)=R(G)$, then there is simply no room inside the isometry group 
for a competing regular group structure.

\paragraph{}
While Niemiec's work aligns with strong naturality, but the conclusion 
"every connected Lie group is strongly natural" does not follow from Niemiec + Mazur-Ulam.
The reason is structural. Niemiec gives, for a connected nonabelian Lie group $G$, 
a suitable one-sided invariant metric whose full isometry group is 
the regular action. After inversion, this proves ordinary naturality exactly 
as we used it. Strong naturality asks about the bi-invariant situation, which 
is much more restrictive.
A connected Lie group admits a compatible bi-invariant metric only when it is a 
"small invariant neighborhoods" group. For connected Lie groups this is essentially 
equivalent to the adjoint action having relatively compact image; 
at the Lie-algebra level, one needs an ${\rm Ad}(G)$-invariant positive-definite 
inner product. Thus groups such as $\mathbb{R}^n,\mathbb{T}^m$ or compact
connected Lie groups and products of these do  work. However, many connected Lie groups 
do not. The group $SL(2,\mathbb{R})$ for example does not allow for a bi-invariant
metric. 

\paragraph{}
For finite groups the Cayley graphs of the semi-direct product is the 
{\bf Zig-Zag product} $\Gamma(H,S) \rtimes \Gamma(K,T)$ 
for the corresponding Cayley graphs $\Gamma(H,S),\Gamma(K,T)$.
This is defined as the Cayley graph of 
$H \rtimes K$ generated by $\{ t_1 s t_2 , t_1,t_2 \in T, S \in S \}$.
The name is justified in that 
the three group operations $t_1, s, t_2$ generate a zig-zag curve. 

\paragraph{}
There is a condition which is satisfied if the direct product $(X,d)$ of 
two metric spaces $(X_1,d_1),(X_2,d_2)$ cite{BridsonHaefliger1999}: 
an isometry $T$ of $(X,d)$ decomposes
as a product of isometries if for every $x_1 \in X_1$, there is
a point $T_1(x_1)$ such that $T(x \times X_2) = (T_1(x) \times X_2)$. 

\section{More examples}

\paragraph{}
A {\bf Bethe lattice} can be seen as a Cayley graph. Lets look at the
Bethe lattice of degree $4$ which is the {\bf Cayley graph} of the {\bf free group} 
$F_2 = \{ a,b | \}$ with 2 generators. It allows for an additional group structure 
$C_2*C_2*C_2*C_2=\{ a,b,c,d | a^2=b^2=c^2=d^2=1 \}$. 
Unlike in the 2022 version, where we claimed that this shows that $F_2$ is not natural. 
We still have 

{\bf $F_n$ is natural for $n \geq 2$. } 

\paragraph{}
Watkins proved in \cite{Watkins1976} that a free product of at least two, and at most countably many, 
at-most-countably-generated groups admits a graphical regular representation. Since
$F_2 \cong \mathbb{Z}*\mathbb{Z}$, this applies directly to $F_2$. 

The ordinary Bethe lattice is too symmetric. Naturalizing $F_2$ requires not merely 
weighting the four branches of that same tree, but changing the Cayley graph 
enough to eliminate all automorphisms except the regular $F_2$-action. 
Watkins' graphical regular representation does this. 

{\bf For odd $n$, the group $C_n$ is natural and the graph $C_n$ is natural} 

\begin{proof}
If $n$ is odd, take as the metric space $C_n$ with geodesic metric. Its group of isometries is
the dihedral $D_n$ which has $Z_n$ as a rotation subgroup of index $2$. Any other group structure
on $C_n$ which has all right translations as isometries must be $Z_n$. If $n$ is odd, then $H$
has odd order but any element in $D_ \setminus Z_n$ is a reflection of odd order.
\end{proof}

\begin{proof}
The argument is similar to the argument that $\mathbb{Z}$ is not natural.
Any metric space admitting $Z_n=Z_{2m}$ as a space on which $Z_n$
acts by right translation must be $C_n$ with a rotation invariant metric. 
But this metric space admits the non-abelian group $D_n$ too. 
\end{proof}

\paragraph{}
We just have seen that the group $Z_{27}$ is natural. Among all groups of order 27,
there is also the {\bf Burnside group} $B(2,3)$, the quotient of the free group $F_2$
by the subgroup generated by all cubes of $F_2$. It is also known as the 
{\bf Heisenberg group modulo 3} and that it can be represented as $(Z_3 \times Z_3) \rtimes Z_3$.
Based on this semi-direct product, we claimed in 2022 that $B(2,3)$ is natural. But since this 
was based on the semi-direct product construction, we need to revisit it. It still remains true:

{\bf $B(2,3)$ is natural.} 

The reason is that stronger. There is an argument which both covers 
$Z_{p^n}$ and $B(2,3)$. Every group of order $p^n$ for odd $p$ is natural. 

{\bf Every finite $p$-group with an odd prime $p$ is natural.}

\begin{proof}
Assume $p$ is odd. A $p$-group has order $|G|=p^n$. Partition the set $G \setminus \{1\}$ (it has
even cardinality) into inverse pairs $\Lambda_j=\{a,a^{-1}\}$. Attach distance numbers $\lambda_j$
to each of these pairs $\Lambda_j$ We take them all different and in the open interval $\lambda_j \in (1,2)$. 
We can now define a metric on $G$ by defining $d(x,y)$ to be $\lambda_j$ if $x y^{-1}$ is
in the inverse pair $\Lambda_j$. This is a metric, as symmetry $d(x,y) = d(y,x)$ follows from 
$(x y^{-1})^{-1}=y x^{-1}$
and because $d(x,z) <2 < d(x,y) + d(y,z)$ implies the triangle inequality. We also have $d(x,y)=0$ if and
only if $x y=1$. The metric $(G,d)$ is now right translation invariant because $(y h) (x h)^{-1} = y x^{-1}$
assures $d(xh,yh)=d(x,y)$. The right translation group $R(G)$ is a regular subgroup of order $p^n$, the isometry 
group $I={\rm Isom}(G,d)$. Look at the stabilizer $I_1$ of the identity in $I$. Given $T \in I_1$, then 
$d(1,Tg)=d(1,g)$ forces inverse pairs to be mapped into each other: $T( \{ g,g^{-1}\}) = \{ g,g^{-1} \}$.
Since $T$ can only change members in each of the $m=(p^n-1)/2$ pairs, $I_1$ must be a subgroup of $Z_2^m$.
But $I$ is transitive on the $|I|=p^n$ points and $|I_1|=p^n 2^m$, we must have $m=0$ and the stabilizer is trivial.
Let $H$ be any other group structure on the $p^n$ points $X$ of $G$. The right-translation group $R(H)$ is a regular
subgroup of $I$ so that it is a Sylow p-subgroup of I. Both $R(H)$ and $R(G)$ are Sylow p-subgroups of $I$.
By the {\bf second Sylow theorem} they are conjugated meaning $R(H)$ is isomorphic to $R(G)$, 
so that $H$ is isomorphic to $G$.
\end{proof}

\paragraph{}
Burnside already knew in 1902 that all groups $B(d,3)$ are finite. 
In 2022, we did not know yet which $B(d,3)$ are natural. But since $|B(d,3)|$ 
is known to have order 
$3^{d+\left( \begin{array}{c} d \\ 2 \end{array} \right) +\left( \begin{array}{c} d \\ 3 \end{array} \right)}$,
also $B(d,3)$ are all natural. 

\paragraph{}
The {\bf quaternion group} $Q_8$ has as automorphism group the symmetric
group $S_4$, a group with $6!$ elements. This large symmetry suggests that $G$ can
host itself and become natural. However, the
metric of the Cayley graph of $Q_8$ generated by $S=\{i,j,k,-i,-j,-k\}$
($S=S^{-1}$ gives an undirected graph) allows only to choose $4$
lengths as parameters. And these are the multiplications by $i,j,k,-1$.
Any such a metric also admits the group structure $Z_2^3 = Z_2 \times D_2$, the product of the
$2$-point group $Z_2$ and the {\bf Klein four group} $D_2=Z_2 \times Z_2$.
The quaternion group $Q_8$ is not natural. (Additionally, we know it is 
an example of a non-simple group which does not split.)

{\bf Every right invariant metric on $Q_8$ also admits a regular $Z_2^3$ action. 
The quaternion group $Q_8$ is not natural and not strongly natural.}

Lets give a more formal proof:

\begin{proof}
Write $Q_8 = \{ \pm 1, \pm i, \pm j, \pm k\}$ and pick an arbitrary right invariant metric $d$.
Define $l(g)=d(g,1)$ the distance to the origin $1$. Right invariance shows $d(x,y) = l(x y^{-1})$ 
and the symmetry assumption $d(x,y)=d(y,x)$ shows $l(g)=l(g^{-1})$. Every right invariant metric
is so determined by $l(-1),l(i)=l(-i),l(j)=l(-j),l(k)=l(-k)$. We can visualize the metric group
as a cube $Z_2^3$ in which the width, length and height are given by the lengths. We have just seen
that this metric is invariant under translations of $Z_2^3$. The metric compatible with $Q_8$ must
also be compatible with the group $Z_2^3$. But since $Q_8$ is non-abelian and $Z_2^3$ is abelian
they can not be isomorphic. Hence, $Q_8$ admits both group structures. It can not be natural.
\end{proof} 

\paragraph{}
The property of being a non-simple but non-split group is shared also by
higher {\bf dicyclic groups} ${\rm Dic}_n$ and we have $Q_8=Dic_2$. Even-so
$1 \to Z_{2n} \to Dic_n \to Z_2 \to 1$ is a short
exact sequence, we do not have $Dic_n = Z_{2n} \rtimes Z_2$.
The reason is that the possible choices of semi-direct products are either the
dihedral group $D_{2n}$ or then the direct product $Z_{2n} \times Z_2$.
$Dic_n$ has order $4n$ and can be presented as
$\langle a,b | a^{2n}=1,b^2=a^n, b^{-1} a b = a^{-1} \rangle$.
In the quaternion case $n=2$ one has the familiar generators $a=i,b=j, ab=k$.
We know now that all ${\rm Dic}_n$ are non-natural.

\paragraph{}
The group $Q_{16}$ is not natural as every metric
always also admits the group $D_8$ and since $Q_{16}$ has a unique involution 
and $D_8$ has many involutions the groups are not isomorphic. \\

{\bf The generalized dicyclic groups $Q_8,Q_{16},Q_{32},...$ are all non-natural.}

They always admit a dihedral group also. 
Generalized dicyclic groups are precisely one of the two infinite exceptional 
families that cannot possess graphical regular representations \cite{LeemannDeLaSalle2022}.
Indeed, we know now, thanks to that paper that the list of
generalized dicyclic groups form the complete list of all non-natural groups.

\begin{figure}[!htpb]
\scalebox{0.4}{\includegraphics{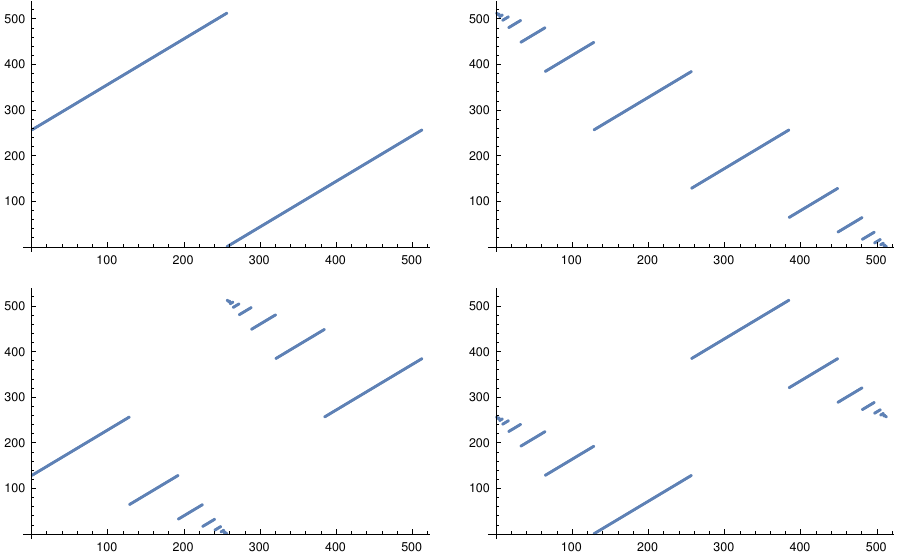}}
\scalebox{0.4}{\includegraphics{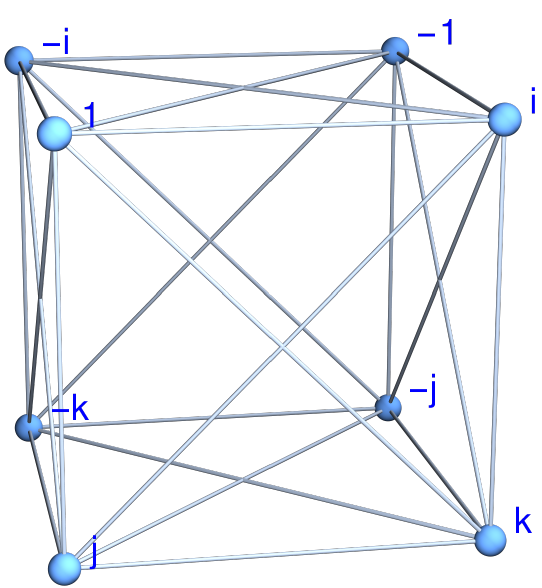}}
\label{Figure 1}
\caption{
The graph to the left generates the natural group $D_5$.
The quaternion group $Q_8=\langle i,j,k | i^2=j^2=k^2=ijk=-1 \rangle$ on the other hand is
not natural.
}
\end{figure}

\paragraph{}
The {\bf Burnside groups} $B(d,n)=F_d/\{ x^n, x \in F_d\}$ generalize the cyclic groups 
$B(1,n) = Z_n = F_1/\{ x^n, x \in F_1=\mathbb{Z} \}$ but they can be infinite
as Golod and Shafarevich established in 1964. 
The {\bf general Burnside problem} asking to identify the Burnside groups $B(d,n)$ which are finite, 
is still open and appears to be a tough problem in 
combinatorial group theory. One already does not know even whether 
$B(2,5)$ is finite. We have not looked at the question of naturalness here.

\paragraph{}
Related to the {\bf Burnside problem} is the recursively defined 
{\bf Grigorchuk group} $G$, a subgroup of the automorphism group of the rooted
Bethe lattice $T$, an infinite tree with a degree $2$ root and all other vertex degrees being 
$3$. It was historically the first example of a finitely generated 
group with {\bf intermediate growth rate} $|B_n(x)|$ of the size of balls of 
radius $n$ in the {\bf Cayley graph} of the group. It is between {\bf polynomial growth}, happening for
example for the free abelian group $\mathbb{Z}^n$ and {\bf exponential growth},
happening for example for the free group $\mathbb{F}_n$. The tree naturally defines a hierarchy of 
trees $T=T_1 \cup T_2, T_1=T_{11} \cup T_{12}$ etc. 
% See \cite{Harpe2022} for some of the history the growth of groups. 
If $T=T_1 \cup T_2$ splits $T$ at the root into two trees, then $a$ flips $T_1 \leftrightarrow T_2$, 
$b$ flips $T_{11} \leftrightarrow T_{12}$ flips $T_{211} \leftrightarrow T_{212}$ 
and acts in the same way on $T_{222}$; the generator $c$ flips $T_{11} \leftrightarrow T_{12}$, 
leaves $T_{21}$ and lets $b$ act on $T_{22}$. The fourth generator $d$ leaves $T_1$, flips 
$T_{211}$ with $T_{212}$ and lets $b$ act on $T_{222}$. 
The four generators all are involutions satisfying $a^2=b^2=c^2=d^2=bcd=1$.
The Cayley graph for $G$ equipped with different distances to the generators $a,b,c,d$ then 
produces a natural group. \cite{DudkoGrigorchuk} produced a continuum of such examples with 
a continuum of non-equivalent growth functions and all featuring isospectral Laplacians. 
The group $G$ was originally defined as an interval exchange transformation on $[0,1]$ 
which has the property that the action of the stabilizer of a dyadic subinterval $[k-1,k]/2^n$ 
is after a natural identification with $[0,1]$ the same group action (see
\cite{BartholdiGrigorchukNekrashevych}). It is an example of a finite automatic group and 
a source for fractals. The pro-finite completion of $G$ has fractal Hausdorff dimension $5/8$.
Grigorchuk's original definition on $\{0,1\}^{\mathbb{N}}$ was $(0x)^a=1x,(1x)^a=0x$,
$(0x)^b=0x^a, (1x)^b=1x^c$, $(1x)^c=0x^a, (1x)^c=1x^d$ and $(1x)^d = 0x, (1x)^d = 1x^b$
is an interval exchange transformation exchanging dyadic intervals similar as
the {\bf von Neumann-Kakuktani system} (=adding machine) generated by the single 
transformation $(0x)^a=1x, (1x)^a =0x^a$. 
The first Grigorchuk group not finitely presentable. 

{\bf The Grigorchuk group $G$ and the Gupta-Sidki group are natural.}

Proof: A graphical regular representation argument of \cite{LeemannDeLaSalle2022} tells that
every finitely generated group that is not virtually abelian admits 
a finite-degree Cayley graph whose full automorphism group is just the regular 
translation action. For infinite, finitely generated groups, the only groups without a 
graphical regular representation are the abelian groups and the generalized dicyclic groups. 
Both graups are finitely generated, nonabelian torsion branch groups 
and so neither is virtually abelian. Therefore each admits a finite-valency 
graphical regular representation, a Cayley graph $\Gamma$ with graph metric $d$. 
Every right translation of the original group is an isometry. 
If another group structure $H$ on the same vertex set also 
had all right translations isometric, its right-regular representation would 
be a regular subgroup of ${\rm Isom}(\Gamma)={\rm Aut}(\Gamma)=G$.
But G itself already acts regularly, so a regular subgroup of ${\rm Aut}(\Gamma)$ 
must have the same underlying permutation group: the unique element sending 
the identity vertex to any chosen vertex is already determined. Hence $H\cong G$.

\begin{figure}[!htpb]
\scalebox{0.6}{\includegraphics{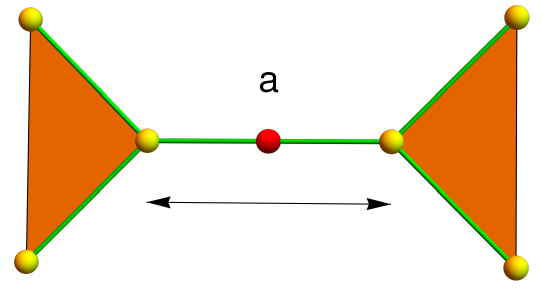}}
\scalebox{0.6}{\includegraphics{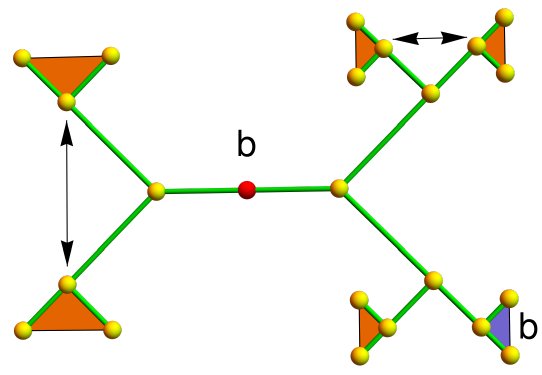}}
\scalebox{0.6}{\includegraphics{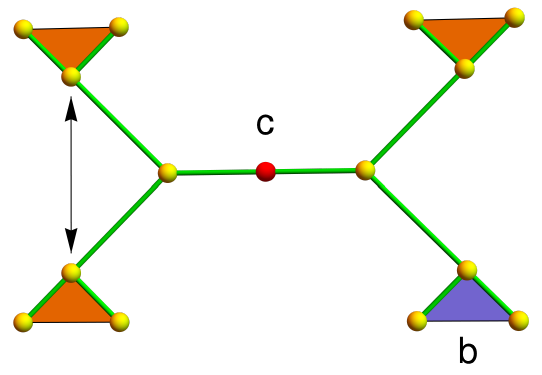}}
\scalebox{0.6}{\includegraphics{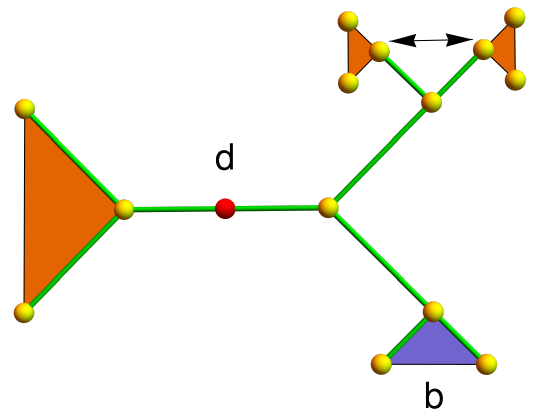}}
\label{Figure 2}
\caption{
The four generators of the Grigorchuk group $G$ acting on a
rooted Bethe lattice. The Grigorchuk group is natural. 
}
\end{figure}

\paragraph{}
Motivated by Grigorchuk is the {\bf Gupta-Sidki group}. It is an infinite $3$-group
which produced a simple example for the {\bf weak Burnside problem} and the first known
$2$-generator infinite $3$-group. The construction obviously generalizes from $p=3$ to any
odd prime $p$. The case $p=3$ is featured prominently in \cite{Baumslag}, 
where Baumslag mentions that he was told that the group is not finitely presentable
but that he would not know of a proof. 
\footnote{Baumslag mentioned this as an open problem in his course from 1987 
and as a student of that course I have at that time made attempts to prove this.
Sidki \cite{Sidki1987} proves in 1984 that Grigorchuk 
gave in 1980 already an argument. \cite{Baumslag} was published
in the year 1993, where the presentation problem is still stated an open problem but this
must have been the state in 1987. It appears as if Sidki's proof needed some time to be
accepted.}
The Gupta-Sidki group acts on a rooted tree of base vertex degree $3$ and otherwise having
constant degree $4$. Pick
an arbitrary root. The first generator $a$ rotates the three main branches $T_1,T_2,T_3$
of the tree. The second generator $b$ rotates the first branch $T_1$ in one direction, 
the second branch $T_2$ in the other direction and plants its action recursively on the
third branch $T_3$. One has $a^3=b^3=1$ and all elements of $G$ have order which is a
power of $3$. To see that the group is natural, pick $G$ itself as the set of generators.
Pick the from the $\{ a,b,a^{-1},b^{-1} \}$ generated geodesic metric 
$d(x,y)=d(0,s)$ where $s$ is written as the shortest word in the letters $a,b,a^{-1},b^{-1}$. 
Now perturb the metric so that all generators have rationally independent lengths. The
geodesic distance in this weighted Cayley graph $X$ produces now a metric space $(X,d)$
which only features the group $G$. The reason is because for odd $p$, the groups $Z_p^n$ are 
all natural. 

\begin{figure}[!htpb]
\scalebox{0.4}{\includegraphics{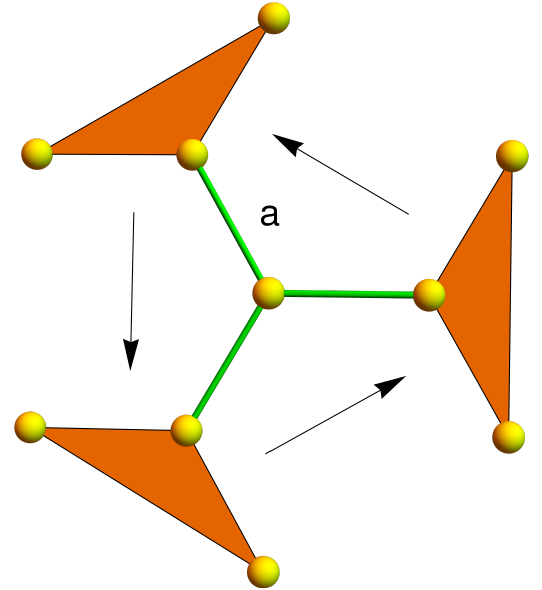}}
\scalebox{0.4}{\includegraphics{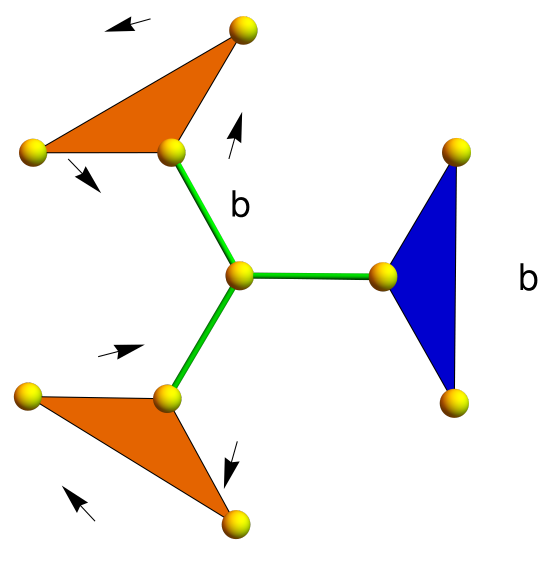}}
\label{Figure 3}
\caption{
The two generators of the Gupta-Sidki group acting on the 
rooted tree of degree $4$ and base degree $3$. 
}
\end{figure}
% MAKE THIS THREE DIMENSIONAL. THE ACTION IS ON A ROOTED TREE AGAIN

\paragraph{}
The group $\mathbb{Z}_2$ of {\bf dyadic integers} is not 
natural. One can generate it as a subgroup of the automorphism group of the tree $X$ as
before. Define $U(x,y) = (y,U(x)), V(x,y)=(V(x),b)$, then define $A(x,y) = (y,x)$
and $B(x,y) = (V(y),U(x))$. 
Then $A^2=1$ and because $VU=UV=Id$ also $B^2=1$. 
(Just check $B^2(x,y) = B( V(y),U(x) ) = (VU(x),UV(y)) = (x,y)$ )
We also have $AB(x,y)=(U(x),U(y))$ and $BA(x,y) = (V(x),V(y))$. 
So, $AB=T$ adds $1$ on both branches while 
$BA$ subtracts $1$ on both branches. 
We have written the von Neumann Kakutani system as a product of 
two involutions. 

\lstset{language=Mathematica} \lstset{frameround=fttt}
\begin{lstlisting}[frame=single]
(* nonabelian group on dihedral group of integers    *)
U[x_]      :=x;         V[x_] :=x;                   n=8; 
U[{X_, Y_}]:={Y ,U[X]}; V[{X_, Y_}]:={V[Y],   X};
A[{X_, Y_}]:={Y , X};   B[{X_, Y_}]:={V[Y],U[X]}; 
P[x_]:=Partition[x, 2]; F[x_]:=ListPlot[Flatten[x]]; 
u=Last[NestList[P,Range[2^(n+1)],n]];
GraphicsGrid[{{F[A[u]],F[B[u]]},{F[A[B[u]]],F[B[A[u]]]}}]
\end{lstlisting}

\begin{figure}[!htpb]
\scalebox{0.5}{\includegraphics{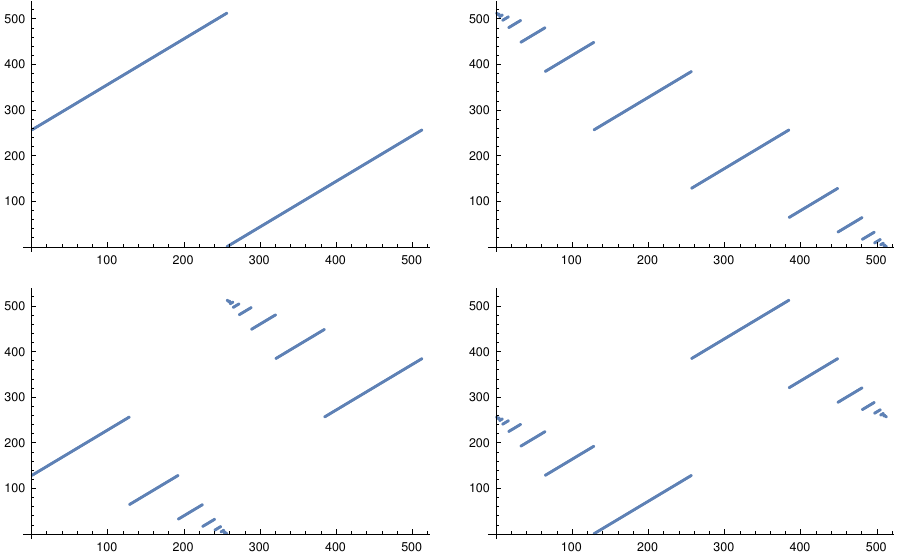}}
\label{Figure 5}
\caption{
The permutations $A,B$ and $T=AB, T^{-1} = BA$ leave the 
spheres $S_k(0)$ in the tree $X$ in variant. The transformations $A,B$
induce on each sphere $S_k(0)$ a dihedral group action. 
The pro-finite limit of dihedral groups $D_{2^n}$ produces the
dihedral dyadic group. }
\end{figure}

\paragraph{}
When trying to find a group which corresponds to Barycentric refinement in two dimensions, 
one can look at the following recursively defined group on a rooted tree $X$ with constant 
degree $4$ except at the origin where the degree is $3$. We have $T(a,b,c)=(T(b),c,a),
S(a,b,c)=(b,S(c),a)$ and $U(a,b,c)= (b,c,U(c))$. We can check that $T(S(U))) = Id$. 
But the group is not commutative, $ST \neq TS$. With $J(a,b,c)=(J(b),J(c),J(a))$ and
its inverse $J^{-1}=K(a,b,c) = (K(c),K(a),K(b))$, we have $J^{-1} S J = U$ and 
$J^{-1} T J = S$. The group is generated by any two members of $\{T,S,J\}$ 
like $\langle S,J \rangle$. If we take the Cayley graph with these generators and 
take different distances $d(0,S(0))$ and $d(0,J(0))$, we have a metric space
which only admits this group action. 

\begin{figure}[!htpb]
\scalebox{0.3}{\includegraphics{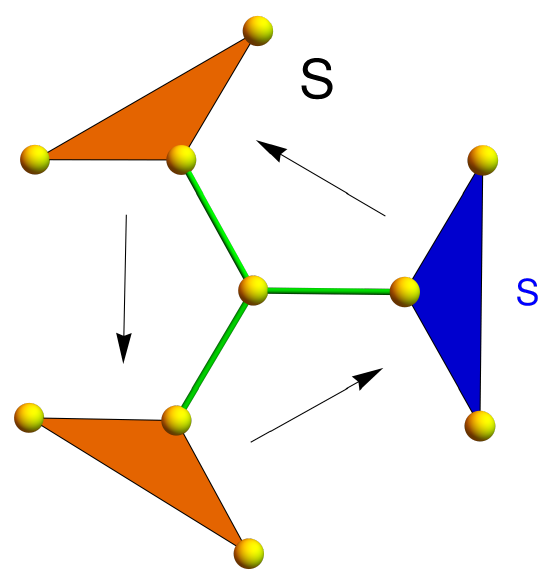}}
\scalebox{0.3}{\includegraphics{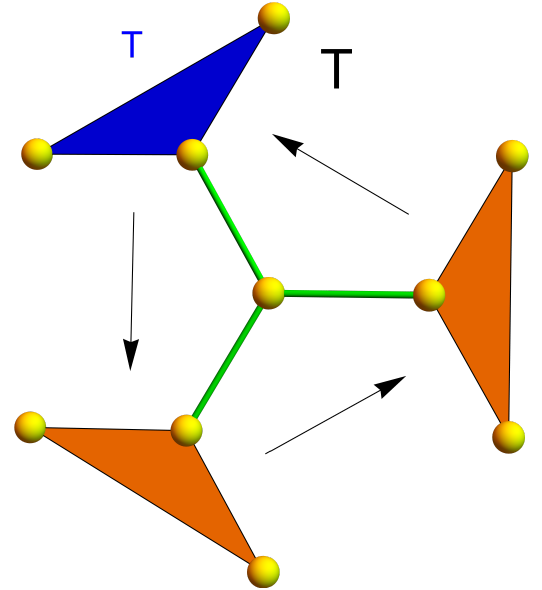}}
\scalebox{0.3}{\includegraphics{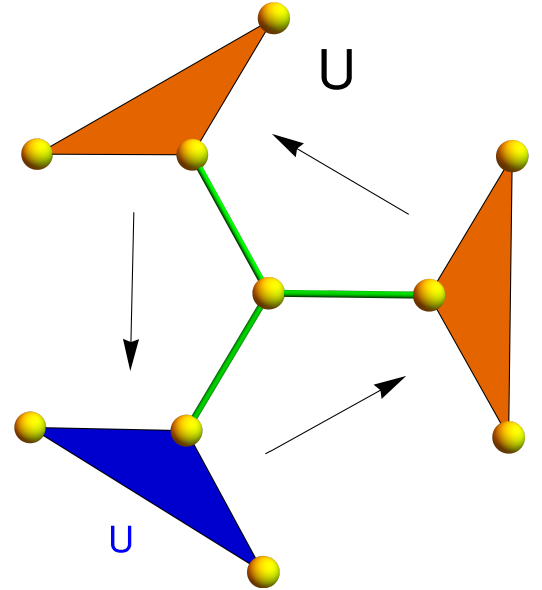}}
\scalebox{0.3}{\includegraphics{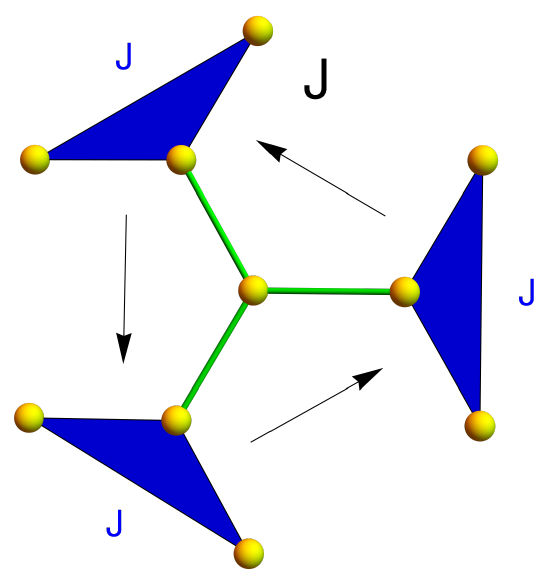}}
\label{Figure 5}
\caption{
$S,T,U$ generate the dihedral 3-adic integer group.
Each rotates the three main branches and induces itself on one of the branches. 
$G$ is also generated by $S,J$ as $J,J^2=J^{-1}$ 
conjugate the three generators. From $STU=Id$ one can see that any two of the 
four actions generate the full group.
}
\end{figure}

\begin{figure}[!htpb]
\scalebox{0.9}{\includegraphics{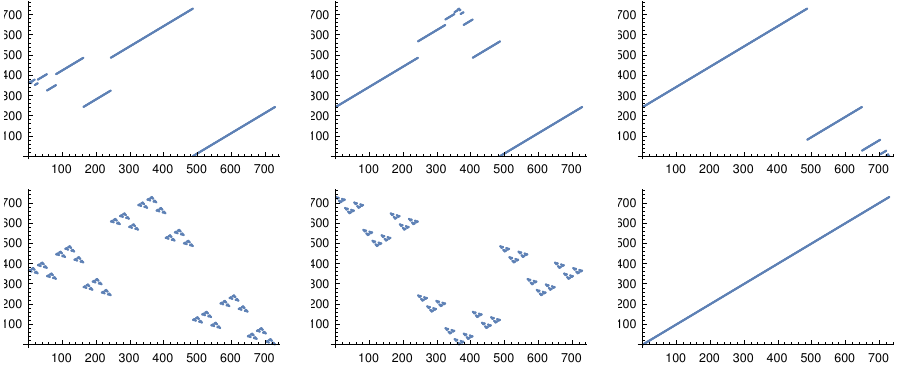}}
\label{Figure 6}
\caption{
The graphs of S,T,U as well as $J$, $J^2$ and $S^3=1$. 
We have $\langle S^3 \rangle \sim \langle T^3 \rangle \sim \langle U^3 \rangle
= \mathbb{Z}_3$ so that $\langle S,T,U \rangle$ is the non-abelian 
$\mathbb{Z}_3 \rtimes Z_3$. It shows that $(\mathbb{Z}_3,d)$ with dihedral metric
is not natural.
}
\end{figure}

\paragraph{}
Here is the Mathematica code which produces the graphs of 
$T,S,U$ as well as of $J,J^2,J^3$. 

\lstset{language=Mathematica} \lstset{frameround=fttt}
\begin{lstlisting}[frame=single]
n = 5; T[x_] := x; S[x_] := x; U[x_] := x; J[x_]:=x; K[x_]:=x;
T[{a_, b_, c_}] := {T[b], c, a};
S[{a_, b_, c_}] := {b, S[c], a};
U[{a_, b_, c_}] := {b, c, U[a]};
J[{a_, b_, c_}] := {J[b],J[c],J[a]};
K[{a_, b_, c_}] := {K[c],K[a],K[b]};
P[x_] := Partition[x, 3];  F[x_]:=ListPlot[Flatten[x]];
R[X_] := Last[NestList[P, X, n]]; u = R[Range[3^(n + 1)]]; 
GraphicsGrid[{{F[T[u]],F[S[u]],F[U[u]]},
             {F[J[u]], F[J[J[u]]], F[J[J[J[u]]]]}}];
\end{lstlisting}

\begin{comment}
GraphicsGrid[{{F[K[S[J[u]]]], F[U[u]]}}]
GraphicsGrid[{{F[K[T[J[u]]]], F[S[u]]}}]
n = 2; T[x_] := x; S[x_] := x; U[x_] := x;
T[{a_, b_, c_}] := {T[b], c, a};
S[{a_, b_, c_}] := {b, S[c], a};
U[{a_, b_, c_}] := {b, c, U[a]};
P[x_] := Partition[x, 3];  F[x_]:=ListPlot[Flatten[x]];
R[X_] := Last[NestList[P, X, n]]; 
u = R[Range[3^(n + 1)]]; uu = Flatten[u];
TT[X_] := Flatten[T[R[X]]]; SS[X_] := Flatten[S[R[X]]]; UU[X_] := Flatten[U[R[X]]]; 
G1 = FindPermutation[uu, TT[uu]];
G2 = FindPermutation[uu, SS[uu]];
G3 = FindPermutation[uu, UU[uu]];
G = PermutationGroup[{G1, G2, G3}]; 
s = UndirectedGraph[CayleyGraph[G]];
GraphPlot3D[UndirectedGraph[Graph[EdgeRules[NeighborhoodGraph[s, 1, 2]]]]]
\end{comment}

\begin{comment}
n = 3; T[x_] := x; S[x_] := x; U[x_] := x;
T[{a_, b_, c_}] := {T[b], c, a};
S[{a_, b_, c_}] := {b, S[c], a};
U[{a_, b_, c_}] := {b, c, U[a]};
P[x_] := Partition[x, 3]; F[x_] := ListPlot[Flatten[x]];
R[X_] := Last[NestList[P, X, n]];
u = R[Range[3^(n + 1)]]; uu = Flatten[u];
TT[X_] := Flatten[T[R[X]]]; SS[X_] := Flatten[S[R[X]]]; 
UU[X_] := Flatten[U[R[X]]];

VV = {uu}; EE = {};
Do[
 XX = {}; Do[XX = Append[XX, TT[VV[[k]]]], {k, Length[VV]}];
 YY = {}; Do[YY = Append[YY, SS[VV[[k]]]], {k, Length[VV]}];
 ZZ = {}; Do[ZZ = Append[ZZ, UU[VV[[k]]]], {k, Length[VV]}];
 WW = Union[XX, YY, ZZ];
 Do[If[VV[[k]] != WW[[l]], EE = Append[EE, VV[[k]] -> WW[[l]]]], {k, 
   Length[VV]}, {l, Length[WW]}]; VV = Union[VV, WW]; 
 EE = Union[EE], {m, 1, 3}];

s = UndirectedGraph[Graph[VV, EE]];
GraphPlot3D[s]
\end{comment}

\begin{figure}[!htpb]
\scalebox{0.8}{\includegraphics{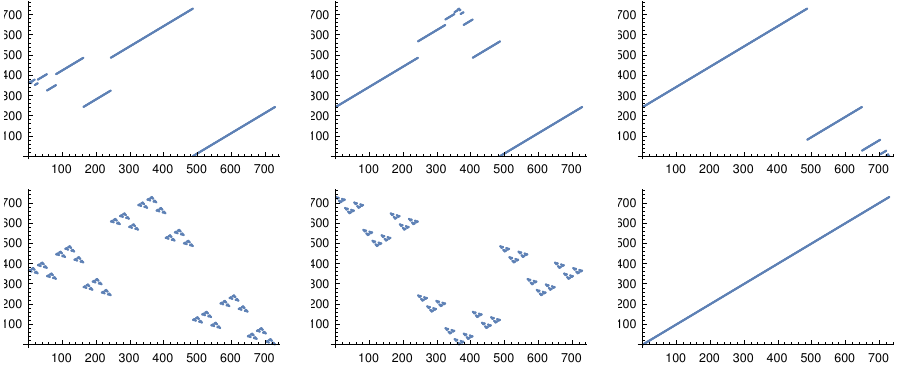}}
\label{Figure 5}
\caption{
The graphs of the permutations $T,S,U$ which generate a non-abelian subgroup 
of ${\rm Aut}(X)$, where $X$ is the rooted Bethe lattice. It is also generated by $T,J$. 
The maps $T,S,U$ induce translations on the triadic group $\mathbb{Z}_3$.
The three fold symmetry $K=J^{-1}$ rotates the alphabet. 
The graphs of $J,K=J^2,J^3=Id$ are shown below. 
}
\end{figure}

\paragraph{}
While the just constructed group did not help us to answer the question what the 
nature of the Barycentric limit in two dimensions is, we can say: 

\begin{propo}
The group generated by $S,T,U$ is a semi-direct product of 
$\mathbb{Z}_3 \rtimes \mathbb{Z}/(3 \mathbb{Z})$.
\end{propo}

\paragraph{}
A result of Bannai \cite{Bannai1977} tells that 
a non-solvable group generated by two elements $x,y$ with $x^2=1$
has a graphical regular representation. It implies a theorem of 
Watkins showing that the simple groups $A_n$ for $n \geq 5$ have
graphical regular representations. Since we can show by hand that $A_2,A_3,A_4$
are natural, we have:

\begin{propo}
All finite alternating groups $A_n$ and symmetric groups $S_n$ are natural.
\end{propo}

\begin{proof}
This has been covered now more elegantly by using that $A_n, n \geq 5$ and 
$S_n$ have presentations in which all generators are involutions.
The old argument here was:
The small cases $A_2=Z_2,A_3=Z_3$ are clear. The group $A_3$ has no graphical
regular representation because the automorphism group of $Z_3$ is $D_3$ which has
$6$ elements. The tetrahedral group 
$A_4$ does not have a graphical regular representation but we can change the 
metric on the truncated tetrahedron graph such that it becomes the automorphism
group. For $n \geq 5$ there is a graphical regular representation so that
we have already naturalness without the need to change the metric. 
\end{proof} 

\paragraph{}
In the context of {\bf graph representation theory} of finite groups $G$, the number
$\mu(G)$ given by the minimal $|V|$ for which $(V,E)$ is a graphical representation 
for $G$ is interesting. The obvious choices can often be improved. There is a 60 point
Cayley graph, a truncated dodecahedron (which 20 vertices truncated to triangles) 
which hosts $A_5$. But there is already $16$-vertex graph for 
example which has the automorphism group $A_5$. 
It is the vertex-face graph of the hemi-icosahedron. For us when we look 
at a metric space to upgrade to a 
group $(A_5,d)$. We need a metric space with $|A_5|=60$ points. The 
truncated dodecahedron does the job. 

\paragraph{}
Related to the notion of naturalness is the {\bf regular representation} of a
group $G$ on the vector space $V$ which is freely generated by basis elements
$b_g$ for elements $g$. Then the right translations are given by 
{\bf permutation matrices}. This is also called the {\bf permutation representation}. 
For the cyclic group $G=C_n$ for example, the matrices are known as 
{\bf circulant matrices}. 

\paragraph{}
The {\bf connection graph} $(V',E')$ of $(V,E)$ is the graph in which the complete
subgraphs $V'$ are the vertices and where two such vertices $x,y$ if their intersection
is non-empty. We have seen (see i.e. \cite{KnillEnergy2020}) that the connection
Laplacian $L$ of a connection graph is always unimodular and has many other 
interesting properties. The inverse of $L$, the Green function, is explicit given.
The entries $g(x,y)$ can be interpreted as potential energy between two cliques $x,y$
in the graph. The sum over all these potential energies is the Euler characteristic 
$\chi(G)$ of the graph $G$ which is $\chi(G) = \sum_x (-1)^{{\rm dim}(x)}$, where the
sum is over all complete sub-graphs. 

\paragraph{}
We claimed in 2022 that the {\bf connection graph}
$C_n'$ of the cyclic graph $C_n$ is always natural reasoning that it has the automorphism group 
$D_n$. But $D_n$ does not act regularly on its vertices because the vertex degrees
are not constant. The $D_4$ symmetry preserves the vertex degrees.

{\bf The connection graph $C_n'$ with geodesic metric is not natural.}

\paragraph{}
What about connection graphs of higher dimensional graphs (meaning graphs having
complete subgraphs $K_3$ already. The simplest two dimensional example is the triangle 
$K_3$. Its connection graph $(V,E)$ is also known as the {\bf Fano plane}, the smallest
finite projective plane. This graph $(V,E)$ has $7$ vertices, the
$7$ simplices of $K_3$. The Barycentric refinement of $K_3$ is the graph in which two 
simplices are connected if one is contained in the other. The connection graph is 
a subgraph of this. It happens to coincide with the {\bf Fano plane}.

\paragraph{}
Note that as a graph, the {\bf Fano plane} (equipped with the geodesic metric on the graph)
can not be a natural metric space because the graph has a prime
number $|V|=6$ of vertices so that it could harbor only the cyclic group $Z_7$.
This does not work since we do not have a cyclic symmetry. The Automorphism group 
of $(V,E)$ is $D_3=S_3$ which has only $6$ elements. 
More generally, the automorphism group of $K_n$ or its connection graph $K_n'$ 
is the symmetric group $S_n$ which has $n!$ elements. As the connection graph 
has $2^n-1$ vertices, it can not be natural. In general, connection graphs or 
Barycentric refinements are examples, where the symmetry group can be much larger
than the number of elements, even so the graphs are not just complete graphs. 

\paragraph{}
A cyclic group $Z_{p^n}$ for an odd prime $p$ 
is natural. So, the group $G=Z_9$ is natural. But $Z_9$ is also example of a 
non-simple but non-split group. We see in examples like $Z_4,\mathbb{Z}, Q_8$ that
non-simple, non-split groups are non-natural.
There are also non-natural groups like $Z_4 \times Z_3$ which are split.
The class of {\bf non-natural} and the class of {\bf non-split, non-simple} groups
has no general relation with each other. But there appears to be a trend. Among non-natural
groups, we see so far many non-split, non-simple groups. 

\paragraph{}
There are exactly two groups of order $6$. The group
$Z_3 \times Z_2 = \langle a,b | a^2=b^3=1, ab=ba \langle$ is not natural
but $D_3 = \langle a,b | a^2=b^2=(ab)^3 \rangle$ is natural.
There are two basic parameters we can choose differently: $d(0,a)$ and 
$d(a,b)$.

\paragraph{}
The metric space $(\mathbb{R},d(x,y) = |x-y|)$ is natural. The metric forces the additive
group structure. The reason is that we know the isometry group of this metric space
to be the non-abelian $\mathbb{R} \rtimes \mathbb{Z}_2$. It is the group generated by point 
reflections. This larger group of isometries is a topological group. But it is not connected
because ${\rm det}(dT) \in \{-1,1\}$ for any isometry $T$ and therefore can not be planted on $\mathbb{R}$. 
It can however be planted on a double cover of $\mathbb{R}$: define
the involution $A(x)=x'$ and $T_y(x)=x+y$ on the first copy and $T_y(x')=x'-y$. 
Then define $B_y(x)=A T_y(x)$. Now $A^2=B_y^2=Id$ and $AB_y(x)=T_y(x)$. This double cover
of the group carries again a group structure. Define $B_x=A T_x$ and $C_x=A T_x^{-1}$. 
If $d(0,B_x(0)) = d(0,C_x(0))$, there are two group structures now. The second one is 
the abelian $\mathbb{R} \times \mathbb{Z}_2$. If 
$d(0,B_x(0)) \neq d(0,C_x(0))$, there is only one group structure. We can force the 
double cover to be natural. Note that also the multiplicative group $\mathbb{R}^+,*)$ 
is natural as it is isomorphic via the exponential function to $(\mathbb{R},+)$. 

\paragraph{}
The non-abelian {\bf dihedral group} $D_n=D/(n \mathbb{Z})$ is natural. 
In order to force the group structure, we need a metric space with $(2n)$ elements
which has a $n$-fold symmetry and a $2$-fold symmetry but we have to avoid a $2n$-symmetry.
In the later case, with a $2n$-fold
symmetry, the product of the groups $\mathbb{Z}/(2n \mathbb{Z})$ 
and $D_n=D/(n \mathbb{Z})$ can be
admitted as a group structure. 
A concrete metric space $(X,d)$ which forces the dihedral groups is a prism graph.
With the metric coming from the embedding of the graph, the
symmetry $Z_n \times Z_2 = \mathbb{Z}/(n \mathbb{Z}) \times \mathbb{Z}/(2 \mathbb{Z})$ is not
allowed as a group structure. The group inversion is 
an isometry but the two neighbors $a,b$ of $0$ are not inverses of each other. There
is a ``twist" or asymmetry in the graph. 

\paragraph{}
We can also take the cyclic graph $X=C_{2n}$, where is the distance $d$ between 
the vertices $2k$ and $2k+1$ is $1$ and the distance between $2k+1$ and $2k+2$ is $2$.
We can also take the {\bf connection graph} of $C_n$. It is natural and
admits only the $D_n$ group structure. The group operations on the edges are
reflections at the edges, the group operations on the vertices are translations. 

\begin{figure}[!htpb]
\scalebox{0.4}{\includegraphics{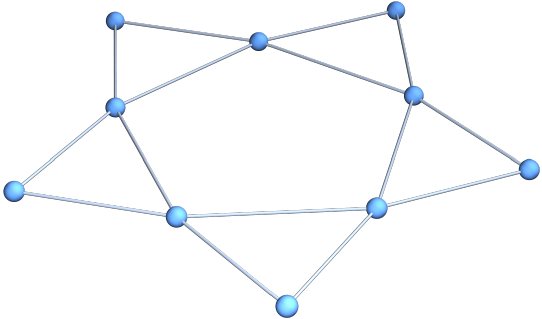}}
\label{Figure 1}
\caption{
The connection graph of a cyclic graph, (here $C_5'$) is not 
natural as there is no transitive group action. 
}
\end{figure}

\paragraph{}
Unlike $(\mathbb{R},|x-y|)$, the
metric space $(\mathbb{Z},|x-y|)$ is not natural. It admits both $(\mathbb{Z},+)$
as well as the infinite dihedral group $(D,+)$ structure and these two groups are 
not isomorphic. The group $(D,+)$ is natural because it comes from a metric on the half
integers which admits only one group structure. Again, as in the finite cyclic case, 
a concrete metric space $(X,d)$ which forces the
group structure $D$ are the half integers $X=\mathbb{Z} + \mathbb{Z}/2$ with the geodesic
distance generated by an asymmetric distance $a$ between $0$ and $1/2$ and $b$ between 
$0$ and $-1/2$. 

\paragraph{}
The {\bf infinite dihedral group} $D_{\infty}$ is the finitely 
presented group $\langle a,b | a^2=b^2=0 \rangle$. Addition in the group is given by 
concatenating  words. The neutral element $0$ is the empty word. 
We can write $1=a,2=ab,3=aba,4=abab,...$ and $-1=b,-2=ba,-3=bab,-4=baba$
showing that no artificial negative elements have to be introduced here. 
One could label elements also as translation and reflection 
symmetries $R_n,S_n$ with $n \in \mathbb{Z}$.
The translations are one one branch, the reflections on the other.

\paragraph{}
In the context of calculus we like to think of $R_n$ as the integers and of $S_n$ as the
half integers. The $R_n$ are translations and the $S_n$ are reflections at the edges of $\mathbb{Z}$. 
If $D_{\infty}$ is seen as acting on itself as a set represented by the half integers $\mathbb{Z}/2$, 
then we can take the reflections $a(x)=-x+1/2$ and $b(x)=-x-1/2$ as generators.
We have $ab(x)=a(-x-1/2)= x+1/2+1/2=x+1$ and $ba(x) = b(-x+1/2) = x-1/2-1/2=x-1$. 
Both $a$ and $b$ switching the pair of integer lines but $ab$ translates us forward. 
Now, if we take a metric on the Cayley graph of 
$D_{\infty} = F_2/\langle a^2,b^2 \rangle = \langle a,b | a^2=b^2= \rangle$
such that $d(0,1/2) \neq d(-1/2,0)$ then this is a natural metric space. 

\paragraph{}
Similarly as the {\bf finite dihedral groups} $D_n = Z_n \rtimes Z_2$ are semi-direct products
of cyclic groups $Z_n$ with $Z_2$, the 
infinite dihedral group is the semi-direct product 
$\mathbb{Z} \rtimes Z_2$ of the integers with the $2$-point group $Z_2$. 
The integers $\mathbb{Z}$ are then a normal subgroup and $Z_2$ is $D_{\infty}/\mathbb{Z}$. 
We like to think about a semi-direct product as a fiber bundle over the base group $Z_2$. The normal
fibers are the integers $\mathbb{Z}$. The group operation is $(s,n) + (t,m) = (s+t,n + (-1)^s m)$
which means that on the first branch (the integers), we have the usual addition, while on the
second branch (the half integers), we have a reverse group operation. 

\begin{propo}
The infinite dihedral group $D_{\infty}$ is natural.
\end{propo}

\begin{proof}
Take the weighted Cayley graph $(X,d)$ of the finitely presented group $\langle a,b | a^2=b^2=1 \rangle$
such that the edges corresponding to $a$ have different lengths than the edges corresponding to $b$. 
Declare $d(0,a) = \alpha >0$ and $d(0,b) = \beta>0$ with $\alpha \neq \beta$. There are no other 
neighbors. This defines a metric space which as a set is just the integers $\mathbb{Z}$ (or half-integers 
$\mathbb{Z}/2$ after a Hilbert hotel identification). We claim that this metric space admits only 
one group structure such that right translations are isometries. 
Its full isometry group is precisely the regular $D_{\infty}$ generated by the two half-integer 
reflections. Since the full isometry group already acts regularly, any subgroup acting transitively 
must be the entire group. Therefore any compatible right-translation group is $D_{\infty}$.
\end{proof} 

\paragraph{}
We can not conclude that the discrete complex plane $D_{\infty} \times D_{\infty}$ 
is natural. The metric space comes from the Shannon product of weighted
graphs. If we take a different scale in the two components, we have no additional
symmetry. The symmetry group is then $D \times D$ and every strict subgroup
is not transitive.

\paragraph{}
The {\bf lamplighter group} is the finitely generated group 
$G=\langle l,r | l^2=1, l r^n l r^{-n} l r^n l r^{-n}=1 \rangle$,
where $l$ stands for ``light change" and $r$ stands for ``walking right". 
We will see below that it is not natural. 
A group element can be represented as a pair $(\omega,n$, where $\omega$ is a
configuration in $\oplus_{k \in \mathbb{Z}} \{0,1\}$ of
configurations with finitely many non-zero elements (lights are on) 
as well as a position position $n \in \mathbb{Z}$ (the lamplighter). 
The addition is $(n,\omega) + (m,\eta) = (n+m,\omega + \eta + e_n$ where $e_n$ is the basis
vector with support at $n$. 
The set of finite configurations can also be written as the countable group 
$\Omega = \oplus_{k \in \mathbb{Z}} Z_2$ or the group of polynomials $\Omega = Z_2[t,t^{-1}]$. 
The generator $l$ lights the lamp at $0$. The generator $r$ walks
the lamplighter one step to the right. The move $r^n l r^{-n}$ lights 
or closes the light at position $n \in \mathbb{Z}$. 
The move $l r^n l r^{-n} l r^n l r^{-n}$ changes a light at $0$, then at $n$, then at $0$ again and
then at $n$ again. After such a double round tour, the lamplighter is at the same 
spot and all the lights in the same state than before the tour. 

\paragraph{}
The lamplighter group can also be considered in finite set-ups, where $\mathbb{Z}$ is 
replaced by the cyclic group $Z_n$. The world on which the lamplighter lives, is now
cyclic. The group is natural. 
Side remark: There is a HNN extensions of the Lamplighter group which appears as the fundamental group 
of a 7-manifold and which shows that bounded torsion is a necessary condition for the 
strong Atiyah conjecture on $L^2$ cohomology of Riemannian manifolds. The conjecture 
of Atiyah that closed manifolds of bounded torsion have rational $L^2$ Betti numbers is still open. 

\begin{figure}[!htpb]
\scalebox{0.4}{\includegraphics{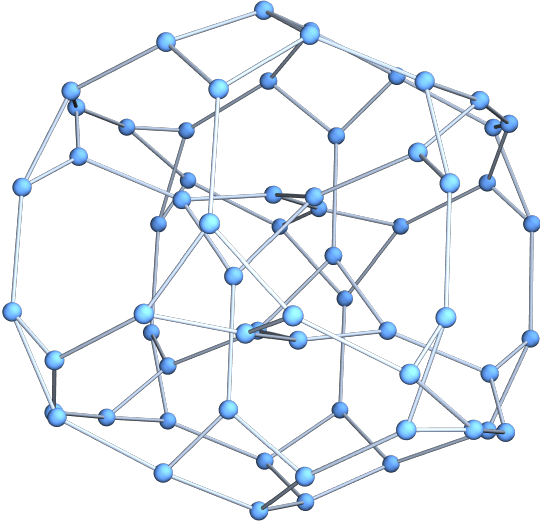}}
\scalebox{0.4}{\includegraphics{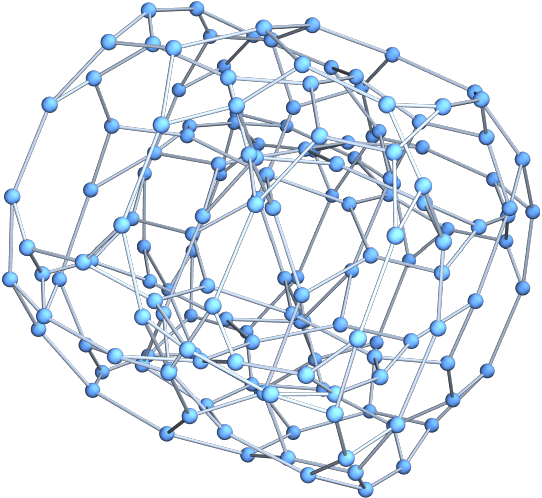}}
\label{Figure 1}
\caption{
Both the lamplighter group over $Z_4$ as well as the lamplighter group
over $Z_5$ are natural. 
}
\end{figure}

\paragraph{}
The lamplighter group is a non-abelian group with countably many elements.
How large is the ball $B_k(0)$ of radius $k$? After $k$ steps, only $k$
lamps can have been lighted and the lighter can only be in distance $k$ away. 
This means that there are maximally $2^k$ group elements in $B_k(0)$. The 
group $G$ has {\bf exponential growth}. The group $G$ is finitely generated 
but not finitely presentable. 

\paragraph{}
The lamplighter group can be written as the semi-direct product 
$\oplus_{k} Z_2 \rtimes \mathbb{Z}$ or the wreath product 
$Z_2 \wr_{\mathbb{Z}} \mathbb{Z}$.
Define the {\bf dihedral lamplighter group} as $Z_2 \wr_{\mathbb{Z}} D_{\infty}$. 
It is obtained by replacing the walk $\mathbb{Z}$ of the lamplighter by the
non-abelian group $D_{\infty} = \langle a,b | a^2=b^2=1 \rangle$. 
Figuratively speaking, the lamplighter's left and right step can now have 
different length. We could call it the {\bf limping lamplighter group}. 
It is presented as
$$ \langle l,a,b | l^2=a^2=b^=1, l w l w' l w l w'=1 \rangle  \; , $$
where $w$ runs over all words in $a,b$ where two $a$ or $b$'s are not neighbors
and where $w'$ is the word read backwards. 

\begin{figure}[!htpb]
\scalebox{0.4}{\includegraphics{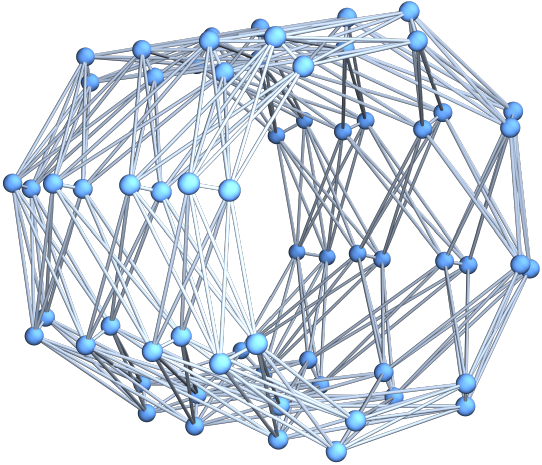}}
\scalebox{0.4}{\includegraphics{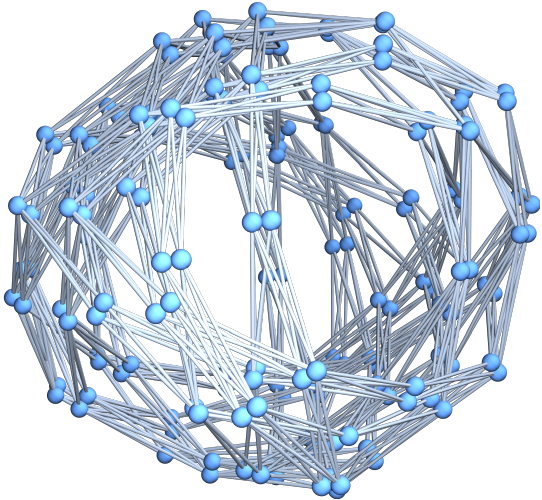}}
\label{Figure 1}
\caption{
The Zig-Zag product of the Cayley graphs in the case of the 
lamplighter group over $Z_4$ and $Z_5$. We use now a different argument 
to show that it is natural.
}
\end{figure}

\paragraph{}
Similarly as we had to go from the abelian group $\mathbb{Z}$ to the non-abelian
group $\mathbb{D}$ to get a natural group, we have to go from the 
lamplighter group to the limping lamplighter group. 
In 2022, we claimed that the lamplighter group is not natural. All lamplighter
groups are actually natural. The mistake was based on our original confusion about
the definition of "natural": right invariance plus 
commutativity forces inversion to be an isometry. For the non-abelian lamplighter group,
a merely right-invariant metric does not force inversion or reflection of 
the underlying line to be an isometry.

{\bf All lamplighter groups are natural.}

\begin{propo}
The lamplighter group $L=Z_2 \wr \mathbb{Z}$,
the dihedral lamplighter group $L_D=Z_2 \wr \mathbb{D}_\infty$ and all 
cyclic lamplighter groups $L_n = Z_2 \wr Z_n$ are natural.
\end{propo}

\begin{proof} 
a) \cite{LeemannDeLaSalle2022} characterizes finitely generated groups admitting a regular graphical
representation. The Leemann-delaSalle theorem assures that every finitely generated group that 
is not virtually abelian admits a finite-degree Cayley graph whose full automorphism group 
consists exactly of the translations. The lamplighter group is not virtually abelian because of 
its growth length. After $m$ steps, the lamplighter can visit m consecutive positions and 
independently choose the states of those m lamps, producing $2^m$ different elements. 
The lamplighter group has exponential growth, whereas finitely generated virtually 
abelian groups have polynomial growth. \\
b) While $L=\oplus_{n \in \mathbb{Z}} Z_2 \rtimes \mathbb{Z}$, the dihedral lamplighter
group is $L_D = \oplus_{n \in \mathbb{Z}} Z_2 \rtimes D_{\infty}$. 
It is generated by $l,a,b$, where $l$ toggles the lamp, and $a,b$ are the generators of 
the infinite dihedral group. We therefore have the group generated by 3 involutions $l^2=a^2=b^2$. 
As seen before this allows the involution-generator lemma to apply. \\
c) The finite cyclic lamplighter groups $L_n = Z_2 \wr Z_n = Z_2^n \rangle Z_n$
has $2^n n$ elements. Parity does not matter. For $n \geq 4$, the group $L_n$ is 
nonabelian, has order at leas 63 and is outside the finite graphical regular representation
exceptions. It therefore is natural. \cite{LeemannDeLaSalle2022} give the GRR classifications.
As for the small $n$ cases: $L_1=Z_2$ is natural, $L_2=Z_2 \wr Z_2 = D_4$ is natural.
$L_3=Z_2 \times A_4$ has a graphical regular representation.
\end{proof} 

\paragraph{}
More general lamplighter groups have been considered \cite{Amchislavska}.
For any commutative ring $R$, one can look at $R[t,t^{-1}] \rtimes \mathbb{Z}$
for example. The action of the base $\mathbb{Z}$ onto the fibers has to be 
specified. One can also make the base higher dimensional, like the
Baumslag-Remeslennikov example $\mathbb{Z}[x,x^{-1}(1+x)^{-1}] \rtimes \mathbb{Z}^2$,
where $(t,s) \in \mathbb{Z}$ acts on the fibers by multiplying the Laurent polynomial
with $x$ or $(1+x)$, the first example of a finitely presented group with an abelian 
normal subgroup of infinite rank. 

\paragraph{}
A {\bf metabelian group} which by definition is a group
for which there is an abelian normal subgroup $H$ such that $G/H$ is abelian. 
The lamplighter group, all dihedral groups, the infinite dihedral group,
all finite groups of order smaller than $24$
are metabelian. Also the limping lamplighter group is metabelian. 
The symmetric group $S_4$ of order $24$ is the smallest non-metabelian finite group. 
A metabelian group $G$ has $[G,G]$ as an abelian normal subgroup and since abelianization 
$G/[G,G]$ is naturally abelian, the metabelian property is equivalent to having the 
commutator subgroup $[G,G]$ as a normal subgroup.

\paragraph{}
The {\bf dyadic group of integers} $\mathbb{Z}_2$ is the dual group
of the {\bf Pruefer group} which is the group of {\bf rational dyadic numbers} $e^{2\pi i k/2^n}$.
The group $\mathbb{Z}_2$ is the unit ball in the field $\mathbb{Q}_2$ of 2-adic numbers 
with $2$-adic norm. A dyadic integer is given by a sequence $x=a_1 a_2 a_3 \dots$ with $a_k \in \{0,1\}$. 
The addition is point-wise with carry-over. If only finitely many 
$a_k$ are non-zero, we can associate with $x$ an integer 
$\sum_k a_k 2^k$. The group $\mathbb{Z}_2$ is not natural unlike $\mathbb{Z}_p$ for odd $p$
which are all natural. 

\paragraph{}
We claimed in 2022 that $\mathbb{Z}_2$ is not natural as it is 
the {\bf pro-finite limit} of non-natural groups $Z_{2^k}$.  The statement that
$\mathbb{Z}_2$ is natural is true but the
argument needs to be modified. We also can not use a density argument. 
Indeed $\mathbb{Q},+$ is natural while $\mathbb{Z}_2$ is not natural.
The statement still remains true: 

\paragraph{}
In ergodic theory, one can represent a dyadic integer also as a real number $x=\sum_k a_k 2^{-k} \in [0,1]$ 
and the presentation is unique for all but a countable set. The addition 
$x \to x+1$ is the {\bf von Neumann-Kakutani system}, a measure preserving transformation
which has the Pruefer group as the spectrum. 
While the {\bf $p$-adic group} of numbers is natural for 
odd primes $p$, this does not apply for the subgroup of integers. This mirrors that
$\mathbb{Z}$ is not natural, while $\mathbb{R}$ is natural. We have seen this already
from the fact that $\mathbb{Z}$ is dense in $\mathbb{Z}_p$ and as $\mathbb{Z}$ is not 
natural, also $\mathbb{Z}_p$ is not natural. 

\paragraph{}
Let us look at the compact {\bf Lie group} $G$, where we have a left and right 
invariant metric $g$ coming from the Killing form. But this does not prove naturality. 
Even if $G$ is abelian, there could be a non-abelian competitive group action. 
We comment on this more in the appendix.

{\bf Every connected Lie group is natural. }

\paragraph{}
It follows that the $0$-sphere $Z_2$, the $1$-sphere $T^1$ and the $3$-sphere $SU(2)$
Let us look at the example of the compact Lie group $G=SU(2)$. 
It is natural coming from the natural metric space
given by the unit sphere $S^3 = \{ x \in \mathbb{R}^4, |x|=1 \}$ with the induced rotational
invariant metric. The group of isometries is the Lie group $H=SU(2) \times SU(2)/Z_2$ 
as the center of the group is $Z_2 = \{ 1, -1 \}$. As the only connected subgroup of $H$ of 
dimension $3$ is $SU(2)$, we see that $SU(2)$ is natural. 

This is related to the classification of list of normed real division algebras 
$\mathbb{R}, \mathbb{C},\mathbb{H}$ classified by Frobenius. 
As they are compact Lie groups the multiplication is unique.
Any other Euclidean sphere does not admit a topological group structure
\cite{Samelson1940}. But this does not settle the naturality question for 
higher dimensional spheres. 

\paragraph{}
One can ask which compact Riemannian manifolds are natural. We have seen that
the zero, one and three dimensional spheres $Z_2,T^1,SU(2)$ are. 

It follows for example that the $2$-sphere $S^2$ is not natural. 
A natural Riemannian manifold needs a large group of isometries. 
The group of isometries of an even-dimensional Riemannian manifold $M$ of 
{\bf positive curvature} is a finite dimensional Lie group and the fixed point set $N$ is again a positive 
curvature manifold with each component having even codimension \cite{Conner1957,Kobayashi1958}.
Interesting is already the case when the group of isometries contains a circle. This restricts the 
structure of even dimensional positive curvature manifolds considerably. For example,
if the fixed point set has co-dimension $2$, then a theorem of Grove and Searle applies
\cite{GroveSearle,GroveSearle2020}. 

\paragraph{}
We mentioned in 2022 already that it looks promising to decide about {\bf dicyclic groups} $Dic_n$.
as this was for $n=2$ the quaternion group $Dic_2=Q_8$ and for $n=1$ is the 
cyclic group ${\rm Dic}_1=C_4$ (which usually is not considered dicyclic as it is already cyclic). 
As both ${\rm Dic}_1$ and ${\rm Dic}_2$ are non-natural, it was only obvious to
suspect that ${\rm Dic}_n$ are non-natural also for $n>2$. As pointed out in the first section, 
this is now all resolved, thanks to \cite{LeemannDeLaSalle2022}.
The {\bf generalized dicyclic groups} ${\rm Dic}(A,y=x^2)$ are
defined by an abelian group $A$ and an element added element $x$. The additional ``imaginary"
element $x$ satisfies $x^2=y \in A$. The group ${\rm Dic(A,y)}$ is then the group generated by $A$ and $x$. 

\paragraph{}
The alternating group $A_4$ is the only 
non-simple alternating group. It is also known as the {\bf tetrahedral group}.
It is a semi-direct product of the {\bf Klein group} $Z_2 \times Z_2$ and 
the cyclic group $Z_3$. A natural metric space can be given by the
Shannon product of Cayley graphs. One can also see directly that it is the 
only group structure which can be planted on the truncated tetrahedron. 
Since the truncated tetrahedron graph $G=(V,E)$ 
has $|V|=12$ vertices and an automorphism group ${\rm Aut}(V,E)$ of order $12$,
the natural nature is established. 
It follows that the graph complement $G^c$ of the 
truncated tetrahedron is natural. (By the way, $G$ has Euler characteristic
$-2$ and constant curvature $-1/6$ everywhere, 
$G^c$ has Euler characteristic $6$ and constant curvature
$1/2$ everywhere. )

\begin{figure}[!htpb]
\scalebox{0.4}{\includegraphics{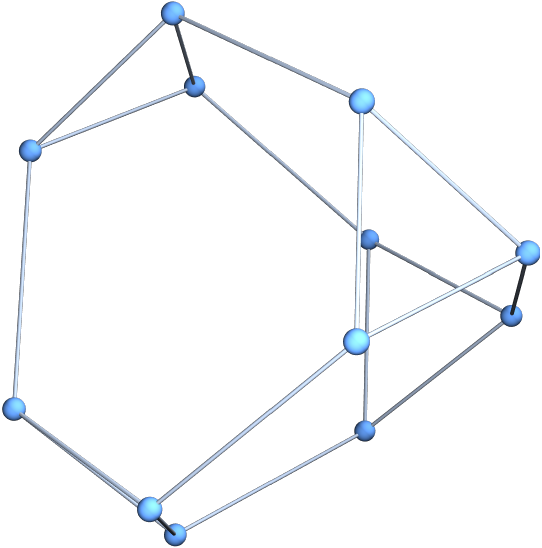}}
\scalebox{0.4}{\includegraphics{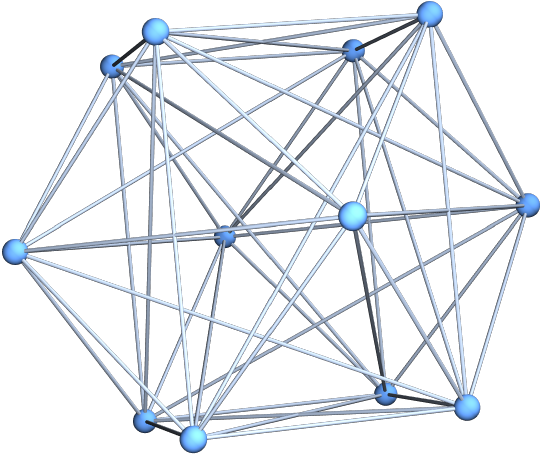}}
\label{Figure 1}
\caption{
The truncated tetrahedron graph is a natural graph. It admits
only one group structure: the alternating group $A_4$. It therefore
is natural. To the right we see the graph complement. 
}
\end{figure}

\paragraph{}
The {\bf alternating group} $A_5$ is the smallest non-solvable group
and the smallest simple non-abelian simple group. As one could
suspect that we need a non-trivial normal subgroup in order that
the group is natural, it is an interesting example. 
Still, the group is natural. It comes from the truncated 
icosahedron graph, a graph with 60 vertices which is an Archimedean solid. 

\begin{figure}[!htpb]
\scalebox{0.4}{\includegraphics{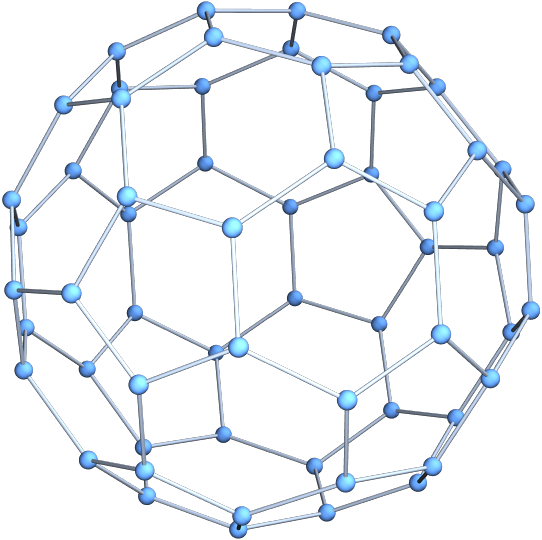}}
\scalebox{0.4}{\includegraphics{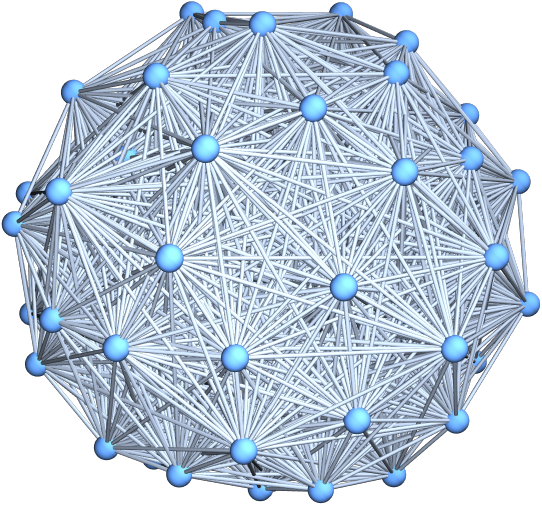}}
\label{Figure 1}
\caption{
The truncated icosahedron graph is natural. It produces the natural 
group $A_5$, the smallest non-abelian simple group. To the right, we 
see the graph complement of the truncated icosahedron. 
}
\end{figure}

\section{Extensions}

\paragraph{}
Given a group $G$ and a normal subgroup $N$, one of the questions in the 
{\bf Jordan-H\"older program} aiming to understand the structure of groups, is to 
determine whether $G$ {\bf splits}, that is whether one can write $G$ as a semi-direct
product $G = N \rtimes H$, where $H=G/N$ is the factor group. The examples 
$Z_4,\mathbb{Z},Q_8$ show that this is not possible in general, even so there 
are normal subgroups in each case. Remarkably, these groups are also non-natural. 

\paragraph{}
Given a group $(G,+,0)$ on can look at the new set $G \times Z_2$
on which there is a {\bf dihedral doubling}: let $a: x \to x'$ denote the identification from
one branch to the other satisfying $a^2=1$. Also assume $a g a = g^{-1}$. 
For every group element $g$ define the new group element $b_g= a g$. It is an involution.
The new dihedral upgrade $\phi(G,d,+,0)$ of the group is the group generated by these 
involutions. It is an example of a {\bf semi-direct product}.
Now, if $d(0,a) \neq d(0,b_x)$ for all $x$, then the metric alone determines the 
action of $a$. For now, we need the technical assumption that there is a distance in $G$
not assumed: $\bigcup_{x \in G} d(0,x) \neq [0,\infty)$. 
This is the case if $G$ is discrete or compact. 

\begin{propo}
If $(G,d,+,0)$ is a natural discrete or compact group,
then $\phi(G,d,+,0)$ is a natural group. 
\end{propo}

\begin{proof}
We have to show that with a suitable metric on $(G,d)$,
there is no other group structure possible on $\phi(G,d)$
than the dihedral doubling construction. By assumption we can 
find a distance $k$ which does not appear in $G$. Define $d(0,0')=k$. 
By group invariance, we have $d(x,x')=k$. Given an element $x \in G$ we
can find its dihedral double $x'$ because the sphere of radius $k$ consists
of only one element. For the subgroup $(G,d)$, we have
no choice. On the conjugate part $(G',d')$, we have an 
isomorphic structure $G'$. There is an image $0'$ of $0$
and the only possible identifications are $x \to x'$ or $x \to -x'$. 
With the former identification, there are two group structures
which work, with the later, we have to define a metric on the double 
cover which selects the twisted one. If we choose a metric satisfying 
$d(0,ag) \neq d(0,ga)$, this excludes the product group structure 
$G \times Z_2$. 
\end{proof}

\paragraph{}
This is related to the {\bf generalized dicyclic group construction}.
If $y$ is an element of order $2$ in an abelian group, define an
other element with $x^2=y$ and postulate $x^{-1} a x = a^{-1}$. 

\paragraph{}
There are many questions:
what are the fixed points of the doubling operation $\phi$ on the 
category of locally compact metric groups? In the compact case, we also 
have a renormalization of dynamical systems:
given a natural metric space $(X,d)$, there is a unique probability measure
$m$, the Haar measure on the compact topological group $(X,d,+,0)$. 
The corresponding $\sigma$-algebra $\mathcal{A}$ produces now a 
probability space $(X,\mathcal{A},m)$ which is a Lebesgue space. We
can therefore arrange that $X=[0,1]$ and that $G$ is a group of 
measure preserving transformations on $X$. Under the doubling map, 
the operations are rewritten on $X_1=[0,1/2], X_2=[1/2,1]$ and a symmetry 
$A: X_1 \to X_2, A: X_2 \to X_1$ added. As the set of measure preserving 
transformations of $[0,1]$ has a complete metric 
$d(T,S) = m( \{ x \in [0,1], T(x) \neq S(x) \}$ we can look at limit points or
fixed points. In particular, we can find, whether any fixed point of 
$\phi$ is natural. 

\paragraph{}
As a comparison to the dicyclic extension, the 
2:1 integral extension for measure theoretical dynamical system has
a unique fixed point, the {\bf von Neumann-Kakutani system} which is conjugated 
to a group translation on the group $\mathbb{Z}_2$ of dyadic integers. This operation
also works on the class of compact metric spaces: define the distance between
two points $x,y'$ as $1+d(x',y')/2$ and normalize so that the maximal
distance is $1$. 

\paragraph{}
Let us write $Z_n=\mathbb{Z}/(n \mathbb{Z})$ for the {\bf cyclic group} with $n$ elements
and $\mathbb{Z}_p$ for the {\bf group of p-adic integers}.
The p-adic integers are a {\bf pro-finite group} obtained as an
{\bf inverse limit} $\varprojlim Z_{p^n}$ of finite cyclic groups. 
We tried to prove in general that if we take the pro-finite limit $G$ of
non-natural abelian groups $G_n$ is non-natural. The question of completion 
is tougher in general as pointed out before. 

\paragraph{}
The $2$-adic group of integers is non-natural despite that the pro-finite completion 
argument did not show. 
The {\bf dihedral p-adic group} of integers
$$  D_p = \varprojlim D_{p^n}   \; . $$
however were natural. The self-similar nature does not allow natural groups to 
survive the pro-finite limit. 
It was tempting to use a pro-finite argument as 
cylic groups $Z_{p^n}$ are natural for $n \geq 1$ unlike $Z_{2^n}$ 
are not natural for $n>1$.

%AAA
Is the pro-finite limit of natural groups is not necessarily natural?

\paragraph{}
The group $\mathbb{Z}$ is the {\bf free group} with one generator $F_1$. We can also see
it generated by $\{ a,a^{-1} \}$. With both $a$ and its inverse as generators, the 
Cayley graph s the undirected graph $(V=\mathbb{Z}, E=\{ (n,n+1), n \in \mathbb{Z} \})$. 
We can also see it as the {\bf Bethe lattice} $B_2={\rm Bethe}_2$ for which the vertex
degree is constant $2$. 

\section{Barycentric renormalization}

\paragraph{}
Given a finite graph $G_0$ we can look at the sequence of {\bf Barycentric refinements}
$G_n$. The Graph $G_n$ is obtained from $G_{n-1}$ by taking the complete subgraphs as vertices
and connecting two if one is contained in the other. 
The eigenvalues $L(G_n)$ of the Kirchhoff Laplacian scaled to the interval $[0,1]$ converge
to a universal limit function $F$ which only depends on the maximal dimension $d$ of $G_0$ 
\cite{KnillBarycentric,KnillBarycentric2}. In the case $d=1$, where $F(x) = 4 \sin(\pi x/2)^2$
satisfies $T(F(x))= F(2x)$ with the quadratic map $T(z)=z(4-z)$, the {\bf integrated
density of states} $F^{-1}(x) = (2/\pi) \arcsin(\sqrt{x}/2)$ and density of states
$f(x) = \pi^{-1}/\sqrt{x(4-x)}$ on $[0,4]$ which is the {\bf arcsin-distribution}. 
In the case $d=2$, the function $F(x)$ appears to be non-smooth, even have a fractal nature. 
We expect that there is an almost periodic operator on a compact topological group
which produces this spectrum. We still have not yet identified this group. 
% T[x_]:=4x-x^2; F[x_]:=4 Sin[Pi x/2]^2; Simplify[T[F[x]]-F[2x]]
% G[y_]:=(2/Pi) ArcSin[Sqrt[y]/2]; f[x_]:=Simplify[G'[x]]; 

\paragraph{}
Let us revisit the Barycentric renormalization story in the one-dimensional case
but focus on its relation to {\bf group theory}. Let $X$ be the 
{\bf rooted Bethe lattice}. It is an {\bf infinite tree} with constant vertex degree $3$ except at the 
root $0 \in V(X)$, where the vertex degree is $2$. Let $T$ be the recursively defined 
{\bf graph automorphism} of the tree $X$ which flips the two main 
branches and on one of the two branches induces the same transformation $T$. The transformation
preserves the {\bf spheres} $S_k(0)$ of radius $k$ centered at the origin $0$ of $X$. It induces there a 
cyclic permutation of the $2^k$ vertices in $S_k(0)$. The boundary of the tree $X$ is the 
{\bf dyadic group of integers} and $T$ extends to this {\bf compactification} where it induces
the translation $x \to x+1$, the {\bf adding machine}. This dynamical system is measure
theoretically conjugated to an interval map with a countable set of intervals of $[0,1]$ 
which is also called the {\bf von Neumann-Kakutani system}. One can abstractly get this from
the ergodic system by computing the group of eigenvalues $\{ e^{2\pi k/2^n} \}$ which is called
the {\bf Pr\"ufer group} $\hat{\mathbb{Z}_2}$, which is the group of dyadic rationals modulo $1$
and the {\bf dual group} of the compact topological group $\mathbb{Z}_2$ of dyadic integers. 
In general, any ergodic automorphism of a probability space is conjugated to a group translation on 
a compact topological group (e.g. \cite{CFS}). 

\paragraph{}
From a group theoretical point of view, the permutation induced on the spheres $S_k(0)$ are
cyclic permutations. The corresponding {\bf Cayley graphs} are also known as {\bf Schreier coset graphs} 
(Nebengruppenbilder) of the stabilizer groups $G_k$ which fix the trees on level $k$. They are
cyclic graphs of length $2^k$. One can visualize the permutation $T_k$ on these $2^k$ points by 
plotting the graph of $p$ mapping $S_k(0)=\{ 1,2,3, \cdots, 2^k\}$ to $T_k(p) \in S_k(0)$. 
For every $k$, we have a discretization of the {\bf interval map} of von-Neumann-Kakutani. In the 
limit we get the interval map on $[0,1]$. 
The group generated by $T$ restricted to the boundary $\mathbb{Z}_2$ of $X$ is a 
{\bf pro-finite limit} of finite cyclic groups. 
In the case $d=2$, we expect a similar thing to happen. The problem is however
that the graphs belonging to Barycentric refined graphs are in dimension larger than $1$ never
Cayley graphs, nor Schreier graphs because the both Cayley and Schreier graphs are vertex transitive
and so have constant vertex degree. 

\lstset{language=Mathematica} \lstset{frameround=fttt}
\begin{lstlisting}[frame=single]
(* The adding machine on the dihedral group of integers  *)
k=8; T[x_]:=x;T[{X_,Y_}]:={Y, T[X]}; P[x_]:= Partition[x,2];
p=Last[NestList[P,Range[2^(k+1)],k]];ListPlot[Flatten[T[p]]]
\end{lstlisting}

\begin{figure}[!htpb]
\scalebox{0.5}{\includegraphics{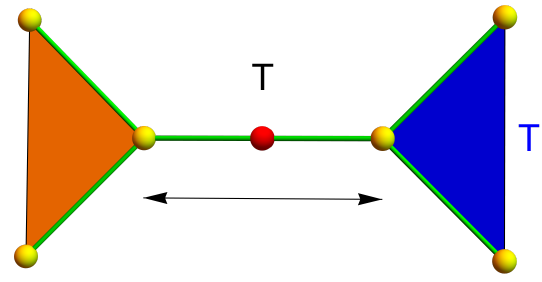}}
\scalebox{0.5}{\includegraphics{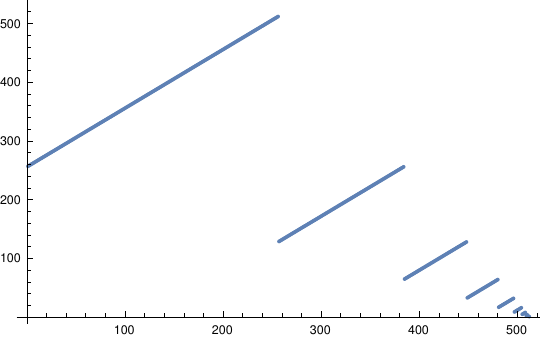}}
\label{Figure 4}
\caption{
The inductively defined transformation $T$ on the rooted binary tree $X$
induces cyclic permutations $T_k$ on the spheres $S_k(x)$. In the limit $k \to \infty$
it produces the von Neumann-Kakutani system $T$ which is conjugated to a group
translation on the dyadic group of integers $\mathbb{Z}_2$. 
}
\end{figure}

\begin{figure}[!htpb]
\scalebox{0.4}{\includegraphics{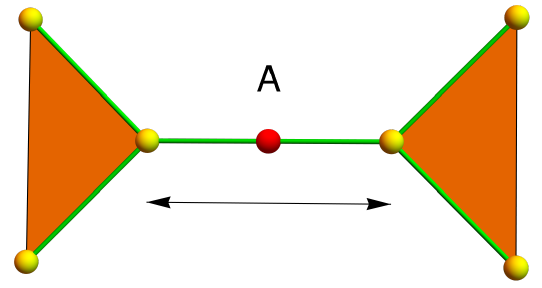}}
\scalebox{0.4}{\includegraphics{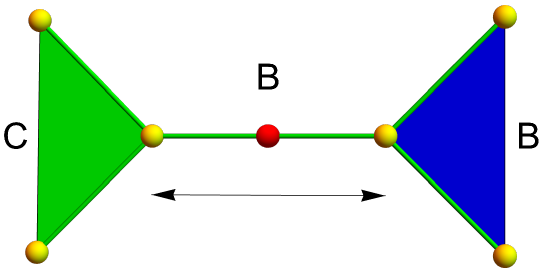}}
\scalebox{0.4}{\includegraphics{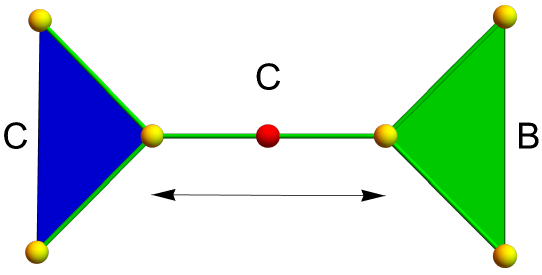}}
\label{Figure 5}
\caption{
The transformations $a,b,c$ which generate a non-abelian group 
structure on the metric space $\mathbb{Z}_2$. 
We have $a^2=1$ and $c=b^{-1}$. Now, $b^2$ generates the von-Neumann
Kakutani adding machine $x \to x+1$ on the right branch
and the inverse on the left branch. 
}
\end{figure}

\pagebreak

\section{Remarks}

\paragraph{}
The first investigation on natural spaces was triggered by an observation we made in 2008 that
the group structure on the integers $(\mathbb{Z},+)$ is not determined 
by the metric as it also allows for the infinite dihedral group. 
This clashed with the case $(\mathbb{R},+)$, where the metric determines the group. 
A second group structure on the integers was then seen as the dihedral infinite group $(\mathbb{D},+)$. 
This group then turned out to be natural and fixed by a metric space. The question of naturalness 
appeared to have affinities with the question which non-simple groups do 
not split. The topic was for us also an opportunity to learn
more about groups and visit old friends like the Rubik puzzle groups or Gupta-Sidki
groups (a group, I learned about from a course of Baumslag \cite{Baumslag}).
To see that the Rubik group is natural, just note that its commutator subgroup is non-abelian so 
that it is not generalized dicyclic. 
We also hoped (later in 2022) that the group theory of finitely generated but not
necessarily finitely presented groups could also help us in the future to answer a
question about the Barycentric limit of simplicial complexes \cite{KnillBarycentric,KnillBarycentric}
which in the one-dimensional case is an {\bf abelian story}, with an underlying 
group given by the dyadic group of integers. In the general case, the limiting
dynamical system is almost certainly a compact non-commutative self-similar group.

\paragraph{}
Having seen that all symmetric groups $S_n$ or alternating groups are natural
and that the semi-direct group $N \rtimes H$ of two natural groups
$N,K$ natural establishes that the {\bf Rubik cubes} is natural.
The most standard $3 \times 3 \times 3$ Rubik cube $G$ \cite{Singmaster1981,CubedRubik,Joyner} is
$$ G = N \rtimes H $$, where $N=Z_3^7 \times Z_2^{11}$ is the normal subgroup   
(the fiber) and where $H=S_8$ is the base. The $2 \times 2 \times 2$ Rubik
cube is $N \rtimes H$, where $N=\mathbb{Z}_3^7 \times \mathbb{Z}_2^{11}$
and $H = (A_8 \times A_{12}) \rtimes \mathbb{Z}_2$.  Games are
a fascinating source for graphs and groups \cite{Games2026}. Almost exclusively, all these
puzzles coming from groups are natural groups. 

\paragraph{}
By the way, a $2 \times 2 \times 2$ version of the Pocket cube was suggested by Larry Nichols in 1970 while
he was a graduate student at Harvard. It was based on magnets and did not feature colors.
Nichols then founded the Moleculon Research Corporation.
He lived from 1959 to 2022 in Arlington, MA and passed away in 2022, 3 months after the
first version of this paper was posted. 

\paragraph{}
We got interested in the subject of group actions also in the context of 
{\bf dynamical systems theory}, the area of mathematics we have worked in 
initially. Dynamical systems theory is naturally 
a {\bf theory of time} as a group $G$ acting on a space $X$ defines a {\bf time evolution} on $X$. 
In the abelian case $G=\mathbb{R}$, one deals with evolution equations like 
ordinary or partial differential equations.
In the case $\mathbb{Z}$, one deals with invertible discrete dynamical systems. 
One can now ask what happens if $G$ acts on itself as automorphism (by right translations!). 
Are there more cases where only one group structure appears? It came to a surprise to us that the 
integers $\mathbb{Z}$ are not natural, as any metric admitting one group structure must admit 
two non-isomorphic group structures. 
But then we realized that the {\bf dihedral integers} $D=\langle a,b, a^2=b^2=1 \rangle$ are natural. 
This group is very close to the integers and one can see it as a pair of two integer lines
or {\bf half integers}. The group $D_{infty}$ has a translation subgroup generated by $ab$, 
a product of two non-commuting involutions. As an index $2$ subgroup, 
$\mathbb{Z}$ is automatically normal in $D_{\infty}$. 
We can see $D$ as a vector bundle over $\mathbb{Z}_2$. This generalizes to higher dimensions
like $DF_2= \langle a,b,c,d | a^2=b^2=c^2=d^2=1 \rangle \cong C_2 \times C_2 \times C_2 \times C_2$ 
which is a dihedral version of the free group $F_2$ with two generators. 
The Grigorchuk group is a $2$-group which is a factor of this group $DF_2$. 

\paragraph{}
Changing from $\mathbb{Z}$ to $D_{\infty}$ 
originates from reflection geometry. One can see {\bf reflections} as {\bf generalized points}, 
similarly as {\bf ideals} can be viewed as {\bf generalized numbers}. 
This also works in Euclidean spaces. Rather than talking
about points $x$ in $\mathbb{R}^n$, we can talk about the reflections $a_x$ at these points. This
defines a non-abelian group satisfying $a_x a_y = - a_y a_x$. Similarly as $D_{\infty} = \mathbb{Z} \times C_2$
is a {\bf double cover} of $\mathbb{Z}$, the {\bf Euclidean reflection group} $\mathbb{R}^n \rtimes O(n)$ 
covers the translation group.
Also similarly as the {\bf Hilbert hotel map} defines a bijective correspondence between 
$\mathbb{Z}$ and $D_{\infty}$, the reflections $a_x$ correspond in a 2-1 way to points $x$ but the 
group structure is different. The subgroup generated by all pairs $a_x a_y$ is the usual translation group,
an abelian group. Remarkably, the reflection group is a {\bf non-abelian additive group} structure. 
Reflection groups are important in mathematics in general. Coxeter groups and in particular Weyl groups
are generated by reflections. The theory of {\bf Coxeter groups} started with Coxeter in 1934. See
\cite{Goodman2004}. 

\paragraph{}
One can see the space $D_{\infty}$ 
also as the set of {\bf half-integers} $-1=ba,-1/2=b,0,1/2=a,1=ab,3/2=aba...$ ,
a group where {\bf word concatenation} is the addition. This non-abelian group is more natural 
written as a multiplication group. If $d(0,a)\neq d(0,b)$ the weighted graph is natural 
allowing the dihedral structure. 
Because $x \to x^{-1}$ is not an isometry in the group, there is an {\bf arrow of time}.

\paragraph{}
The {\bf arrow of time} is a one-dimensional feature.
It disappears in higher dimensions as we can now have different reflections $a_i$. A combination of
three reflections produces in a flat space the zero translation. Doing this is now possible without
violating a Pauli principle. This generalizes on the dihedral free group $DF_n$ where we
can not return to the origin. On the factor group $DZ^n = DF_n/\{ a_i a_j =a_j a_i\}$ however for $n>1$ 
we can turn around without violating that Pauli principle: $abababcbababac$ is the $0$ element. 
The phenomenon of an "arrow of time" is much more clear now as group inversion is 
not an isometry on $D_{\infty}$ if we chose a naturalizing metric on $D_{\infty}$.

\paragraph{}
Writing motion as a product of reflections can generalized to 
Riemannian manifold with non-positive curvature. As then, the reflection $T$ at a point $a$
is well defined thanks to the non-existence of caustics (Hadamard). To define $T(x)$
build the geodesic to $a$ then continue that geodesic through $a$ until the same distance is
done again. The combination of two such reflections is now a translation on the manifold. 

\paragraph{}
In the compact case, where the {\bf Haar measure} is a unique
group bi-invariant measure a natural metric space $(X,d)$ alone defines a unique
{\bf measure theoretical dynamical system} $(X,\mathcal{A},G,m)$ in which the 
group $G$ emerges as the ``time evolution" on itself by right translation.
We do not need bi-invariance to talk about a Haar measure. A right invariant transitive
action suffices. A natural group so upgrades (choosing a defining metric) to a {\bf probability space}.
The metric space especially defines a unique representation of the unique $G$ in the unitary 
group of $L^2(X,\mathcal{A},m)$. The interplay between topological dynamics
and measure theoretical dynamics is rich \cite{DGS}. It is especially interesting
if a topological dynamical system has a unique invariant measure, is uniquely ergodic. 
One can then use both the measure theoretical world as well as the topological world. 
In the first case one can be interested in things which are true almost everywhere
(like convergence of Birkhoff sums) or then in things which are true Baire generically. 
Also recurrence questions can be asked both in the topological as well as in the 
ergodic setting \cite{FurstenbergRecurrence}. 

\paragraph{}
One can now attach properties to points of the metric space. For example, define
the {\bf ergodic set} in a natural compact metric space $(X,d)$
as the set of points $x$ for which the corresponding 
transformation $T_x(y) = y*x$ defines an ergodic 
dynamical system $(X,\mathcal{A},T_x,m)$. For the circle $X=\mathbb{T}_1$
for example with that natural metric, the ergodic set is set of points $x$ for which 
$d(x,0)/diam(X)$ is irrational. We can define spectral types
to every point $x$ depending on the spectral type of the unitary $f \to f(T_x)$
on $L^2(X,\mathcal{A},m)$. 

\paragraph{}
The question of whether a given structure determines an other categorical structure
can be asked for any pair of categories. The set-up can not be reversed in our case: 
we can not switch topology and algebra and ask whether the group
determines a metric. The reason is that there
are in general many metrics on a group for which all the group operations 
are isometries; the {\bf trivial metric} $d(x,y) = \delta_{x-y}$ is always an 
example. It is a metric that is invariant on for any group.

\paragraph{}
Finally, we should mention that the non-abelian dihedral structure of the integers is already 
present in the very first mathematical artifacts we know. The marks of a 
{\bf tally stick} can be identified with $a$. To make the marks visible, they have to 
be separated. The gaps can be identified with $b$. Our 10 thousand year old ancestors have 
not seen the need for negative numbers, but a gap is needed if we want to determine the sign.
A number like $1$ given by $ab$ gives the direction. The symbol $baba$ with the gap to the left
would mean that the mark is red backwards and means $-2$.
Reading a number backwards gives negative numbers $bababa$ for example is $-3$.
In the case of the tally sticks, the gaps are of course only implicit, but in some sense,
the quartz stone planted at the end could (in principle at least) determine in 
which direction one has to read the bone.

\paragraph{}
The {\bf infinite dihedral number line} $D_\infty$ is in in many ways much more natural than 
the {\bf integer number line} $\mathbb{Z}$. We have seen a mathematical explanation by 
identifying the former as a natural group and the later as a non-natural group. 
But there are other ways why the structure
is natural: we do not have rather awkwardly introduce a group completion $\mathbb{N} \to \mathbb{Z}$ 
to get the negative numbers, we have them already given as words, without having to 
introduce a negative sign. We have $ab=1$ and $ba=-1$ for example. We could replace 
the letters $a,b$ with other symbols. In computer science, it is the binary notation with $0$
and $1$ which is the basic way to write but there, no Pauli principle disallowing $aa$ or $bb$
is in place so that much more information can be filled. But such a system would have been 
impractical for the cave mathematicians writing on bones as it would be hard to count 
the number of gaps. Our binary encoding is a place value system which 
historically was first introduced by Sumerians who used a hexagesimal system when writing 
Clay tablets. The dihedral number line given as alternating sequences of two symbols
faded and only reappeared when geometers started to look at symmetry groups. 

\begin{figure}[!htpb]
\scalebox{0.1}{\includegraphics{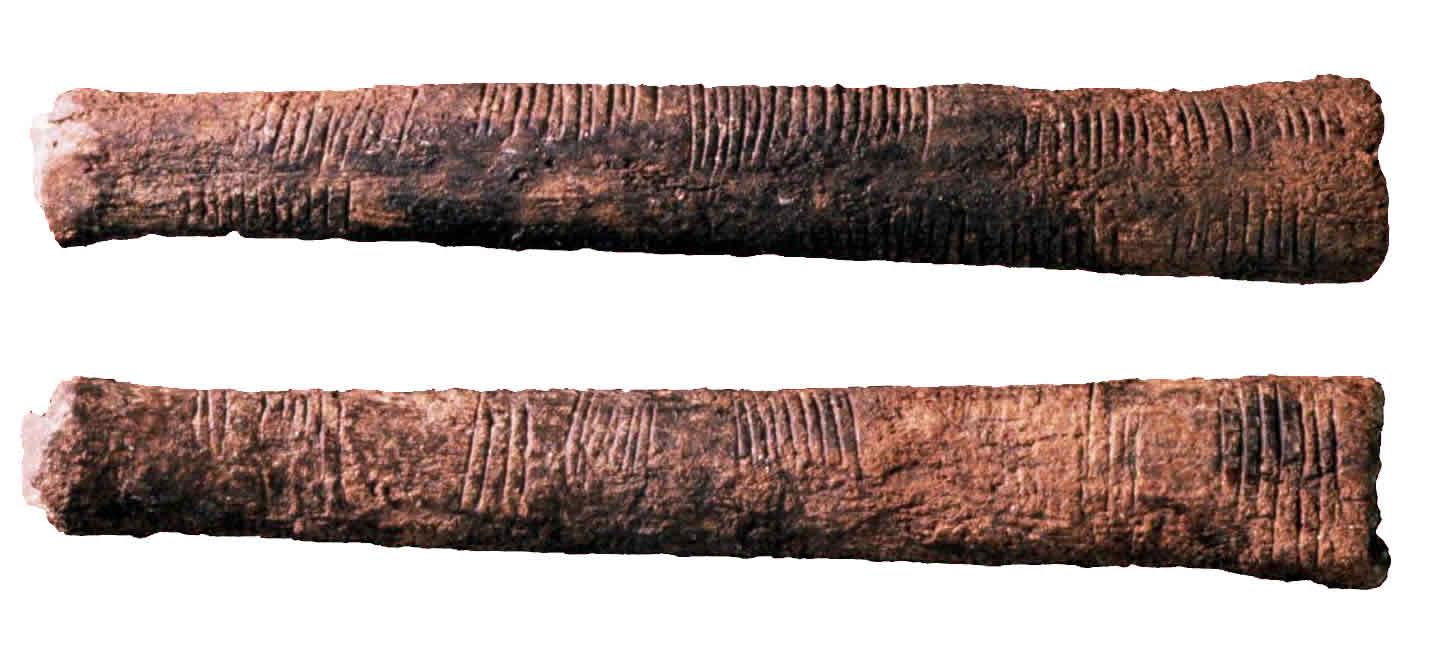}}
\label{Figure 2}
\caption{
The earliest known mathematical documents are tally sticks like the 
Ishango bone. If we think of a scratch as an $a$ and as a gap as a symbol $b$,
then we have a dihedral description. Reading things backwards gives negative
numbers $ababab=3$ and $bababa=-3$. Of course the dihedral non-abelian group 
structure was not on the mind of the earliest mathematicians but marks on 
bones are a way of counting which in some sense already involves two type of
symbols, the mark and the gap between two marks. 
}
\end{figure}

\paragraph{}
It was Felix Klein who shifted in his Erlangen program the interest from the geometric 
metric spaces to its {\bf group of symmetries}. Euclidean symmetries are generated
by reflections leading to ``reflection geometry" \cite{Bachmann1959}.
One can identify for example the set of translations in Euclidean space
with the space itself. This group of translations has a cover in the
form of the group of all point reflections. It is a non-abelian group 
which should be thought of as a dihedral double cover of Euclidean space.

\paragraph{}
We can look at a dihedral versions of the Cayley graph of $\mathbb{Z}^n$ for example,
with generating set $S=S^{-1} = \{ a_1,\dots,a_n, a_1^{-1},\dots,a_n^{-1} \}$ with a
generating set $\{a_1,\dots,a_n, b_1,\dots,b_n \}$ and adding the relations
$a_i^2=b_i^2=1,a_i a_j = a_j a_i, a_i b_j = b_j a_i$ for all $0 \leq i \neq j \leq n$. 
If we postulate that an eligible path in this lattice is one where no Pauli principle
violation occurs (no unnecessary $a_i a_i$ terms for example), then in one dimensions
we are forced to move in one direction, while in for $n \geq 2$ we can pretty much move
around freely in the space. 

\paragraph{}
Our first experimentations in this area relating metric spaces and groups
origins from 2008. Sparks to it was motivated by {\bf inquiry based learning pedagogy}.
Already years before 2008, after seeing experts teaching with
such methods and attending conferences on the subject of IBL, I felt at that time 
that inquiry based learning can only work well if {\bf new and unstudied topics} are
covered. In a well connected world with search engines and social network discussion
groups like stack exchange, maintaining an inquiry based climate has become 
very difficult if not impossible. (Today, in our brave new world of artificial 
intelligence, it is even more challenging.) The pioneers of inquiry based learning worked 
initially with new set theoretical topology, which at that time was not saturated and 
also had a gentle entry point. Today, for many classical topics, one can find answers to 
known questions by just asking the question verbatim on the internet. There are data 
bases of graphs, groups, topological spaces, books on counter examples etc etc.
14 years have passed since 2008 (and 18 since 2026) and the power of search (and AI)
has increased enormously since. We are at a point, where a standard math 
question can be entered as such into the search. (Now, in 2026, we can even pose research 
problems and can expect reasonable answer. Rather than thinking about
a proof for hours, days or weeks, a known proof emerges in a fraction of a second. So, we need
a topic which is completely unexplored and the theme studied here would have been 
a good example. Of course, once this file is on ArXiv, a search for a question might bring
it up. We will need an other new topic. (now in 2026, the web has learned this subject already
well).

\paragraph{}
Arithmetic usually starts with the {\bf integers}, a commutative additive
group. From it, commutative or non-commutative ring or more sophisticated 
structures are built. Could there be something more fundamental than 
the integers? We argue here that the answer could be yes. In order to do so, we have
to specify what ``fundamental" is. We propose to look at arithmetic
structures which emerge uniquely from a metric space, a set equipped with a distance
function. It turns out that the non-commutative dihedral group of 
integers is more fundamental than the integers. Unlike the integers, the
group structure (and so multiplication table) of the dihedral numbers is determined 
from a metric space. The same also holds for the finite dihedral groups. 
While the cyclic groups $C_{2n}$ are not natural in the above sense, the dihedral groups are. 

\paragraph{}
Leopold Kronecker was probably the first {\bf ultra-finitist}. His famous statement that
the {\it "dear lord has created the integers while every thing else is made by humans"},
indicates this. Challenging this statement, one can ask now whether 
there is a more fundamental additive arithmetic than the usual
"cave" or "pebble" arithmetic $(\mathbb{Z},+)$ we know from elementary school.
Of course, this question depends on what one considers "natural" or "fundamental".
But the dihedral structure (as we tried to explain above) is actually very natural. 
It introduces counting in a minimalistic way without number system, without place-value
just using two symbols and defining $1=ab,2=abab,3=ababab,\dots$ having negative
numbers without group completion built in too as $-1=ba,-2=baba,-3=bababa,\dots$. The
numbers are doubled with half integers $1/2=a,3/2=aba, etc$. One of the
surprises was that the group structure is forced from the metric, if one takes a metric
space in the form of a weighted graph such that $d(0,a) \neq d(0,b)$. 

\paragraph{}
An other motivation comes if one looks at the overall structure of mathematics.
One can muse about {\bf which categories are the champions in a major field}.
If one looks at algebra, the concept of a {\bf ``group"} appears to be the most intriguing.
In a topological setting, the concept of 
{\bf ``metric space"} is most fundamental. When looking at order structures
the {\bf ``posets"} are among the top. Partially ordered sets immediately define 
``graphs", where the nodes are the sets and where two are connected if one 
is contained in the other. Alternatively, one can look at the graph in which 
the connections are taken if two sets intersect. This is interesting from the spectral 
point of view for example \cite{KnillEnergy2020}.
What happens if one asks whether one champion
defines an other? For Lie groups for example, 
one can look at the problem of finding a metric structure in the form of a
bi-invariant metric. Now, in 2026, it appears to me that the one sided case is actually 
much easier. 

\paragraph{}
{\bf Double cover structures} are often observed in mathematics. An example is given by 
{\bf projective spaces} which are double covered by {\bf spheres}, {\bf spin groups} 
are double covers of {\bf special orthogonal groups}. 
{\bf Pin groups} are double covers of {\bf orthogonal groups}.
It can make sense therefore to generalize ``points" like the integers $\dots, -1,-1,0,1,2, \dots $ 
to objects generated by ``involutions" $a,b$ leading to $\dots baba,ba,0,ab,abab \dots$.
There is an analogy in number theory, where one transitioned from 
``points" to ``prime ideals". Involutions are very natural objects.
When going from $\mathbb{Z}$ to $\mathbb{D}$
one moves from points (numbers) to elements of the symmetry group (reflections). 
Natural means that the elements of the symmetry group of a space can be 
identified with the space itself. 

\paragraph{}
An early draft of this document from 2008 gained my attention more recently 
when looking at complex analysis in a {\bf discrete setting}.
Natural arithmetic comes up when looking at {\bf ultra-finite single variable 
calculus models} in a situation, where {\bf time} is one-dimensional. 
For ultra-finite calculus on $\mathbb{Z}$ with derivative $df(x)/dx=f(x+1)-f(x)$
and polynomials $[x]^n=x(x-1) \cdots (x-n+1)$, one has $d[x]^n/dx = n [x]^{(n-1)}$. 
If we want to make this algebra work well also on the negative axes,
we need $df(x)/dx = f(x+1/2)-f(x-1/2)$ and use 
$[x]^n = (x-n/2)(x-n/2+1) \dots (x+n/2-1) (x+n/2)$. 
Again, one has $d[x]^n/dx = n [x]^{(n-1)}$. 
The use of half integers reminds of dihedral structures. A doubled lattice appears naturally
when constructing the dihedral integers $D$. It appears also as a Barycentric refinement. 

\paragraph{}
There is an other link. In geometry, other dihedral type structure appear:
if we have a geometry and an exterior derivative $a=d^*$ for example we can also look at the 
dual $b=d$, then $a b = d^* d$
and $ba= d^* d$. While in geometry these are Laplacians, in the case of dihedral numbers, both are 
considered a translation. In the algebra of {\bf creation and annihilation operators} in quantum field theories
we have $a a^* = -a^*a$ and $a a^* - a^* a = 2$. The effect of the
combination $a a^*$ is to increase the number of particles by $1$ (mathematically
the eigenvalue increases). This is exactly what the word $a b$ does in the 
dihedral group. $x \to x+ab$ means adding $1$. The words of even length form a {\bf normal subgroup}
of the dihedral group of integers isomorphic to $\mathbb{Z}$. 
Indeed, $\mathbb{D} = \mathbb{Z} \rangle \mathbb{Z}_2$ is a bundle over the two point space, where
each fiber is the integer line. 

\section{Dihedral time} 

\paragraph{}
This section is a bit different. It touches on a topic which is very popular in the physics literature
because it has also philosophical connections. It is the question {\bf ``What is time?"}.
For a mathematician, the answer is simple: time is a group or semi-group acting on a space. 
The puzzling riddle is when we look at the role of time how it is observed in nature: 
while we can move freely in space, we can not move freely in time. There is something fundamentally different
about time which geometric notions of {\bf space-time manifolds} do not address. Postulating a 
{\bf semi-group} like $\mathbb{N}$ rather than $\mathbb{Z}$ 
does not solve this problem as it just puts the non-invertibility as an axiom. This is Prigogine's 
approach but it feels like a cheap stunt similarly as 
``deus ex machina", ``antropic principle" or ``multiverse landscape" ideas \cite{LandscapeAndMultiverse}. 

\paragraph{}
A major enigma in physics is the {\bf arrow of time}. We observe the inability to influence the
development of time, despite the fact that all fundamental physical laws are invariant 
under group inversion possibly with modification like in the CPT symmetry case, where a time change
comes with charge and parity change. The phenomenon is often formulated within a thermodynamic or cosmological
frame works. But it is something, we experience in small laboratories and in daily life. Already on a
small scale experiments, we have no ``force field" which allows us to change the direction in time.
There are external force fields which allow us to change the direction of a moving particle, we can use physical
obstacles to prevent going an object to a specific place but we can not build obstacles which slow down, stop or
reverse time. We can slow down time by moving around or to placing the experiment in a gravitational field but
that is just a consequence of relativity. 

\paragraph{}
Explanations which have been proposed are entropy \cite{Eddington1928} or 
postulating irreversibility \cite{Prigogine}.
Both are problematic if we accept that fundamental interactions can be
reversed. There is no evidence that any of the known fundamental forces lead to dynamics that is impossible to 
reverse. There is a CPT symmetry observed for example. 
{\bf Time machines} are possible in principle in the frame work of 
general relativity \cite{Thorne} if one postulates the existence of exotic matter producing a mass tensor for which the 
metric of the Einstein equations allows travel backwards in time. Similarly, exotic small matter (like in an ant-men
(Marvel comics) universe) would allow to build a {\bf Maxwell demon}, violating the second law of thermodynamics. A 
more practical implementation of the Maxwell Demon idea is the Feynman-Smuluchowski ratchet, first proposed by 
Marian Smoluchowski. As Feynman illustrated in his lectures, the idea does not work because the ratchet itself
is part of the system. As a mathematician however, we can think about a model in which shrinking of objects like
in the Marvel universe is possible. While ludicrous from a physics standpoint, it is completely consistent with
most of the physics we know. An other argument against entropy is Poincar\'e recurrence. We just have to 
wait long enough to have a system come back close to the situation it has started with. 
Already in 1896, Ernst Zermelo argued against such a mechanism but the argument was just that the time
span is unreasonably long. This is not a mathematical argument, it is an engineering argument. If we have to 
wait trillions of year until something comes back, it is an argument which fails for practical reasons and not
fundamentally. Non-recurrence can occur in infinite dimensional Hamiltonian systems, where no invariant
probability measure might exist. An example is Vlasov dynamics as formulated in a Hamiltonian frame work as
a dynamical system on diffeomorphisms \cite{knillprobability}.

\paragraph{}
While postulates like the irreversibility postulate of Prigogine are bold, they resemble the Greek atomic 
hypothesis. Also there, the suggestion had been based on experiments and experience
(we can not crush matter arbitrarily small) with limited effort)
\footnote{A gypsum particle has the size of maybe 25 -35 micro meters which is 0.025 mm. This is far from 
the size of an atom which is of average size 0.1 nano meters which is about $10^{-7}$ mm.} 
but these limitations were removed more and more until reaching the length scale of the constituents of the
nucleus. We have perfect reversibility for the fundamental equations which describe the motion of 
particles but have to deal with friction phenomena breaking this. Democritus was not able
to crush very small quartz sand particles and postulated that there are smallest units. 
While we can not read the mind of Democritus any more, it is likely that his postulate was
motivated by the empirical fact like that sufficiently refined powder became impossible to be refined any 
further. There was no evidence at the time of Democritus which 
supported an atomic hypothesis. Experimental evidence came only much, 
later with experiments like Millikan's establishing that some particles like the 
electron have quantized fundamental charges.
%  Update 2025: should probably mention Brownian motion, Einstein's work

\paragraph{}
Could a dihedral time-line explain the special nature of time? Once
we have a direction of time and the time evolution is obtained by stacking
reflective operations on top of each other, then we can only move in the direction 
in which we are already going. In some sense, the sequence of letters $a,b$ 
is like the evolution of a wave, the $a$ playing the role of the position and 
$b$ of momentum. In Hamiltonian systems one has
both position and momentum. A wave does not change direction because it has also 
momentum. If we accept that fundamentally (similarly as the Pauli principle), that no pairs 
$aa$ or $bb$ are allowed to occur, then we can only move in one direction. The finitely 
presented Coxeter group $D_{\infty} = \langle a,b | a^2=b^2=1 \rangle$ 
is in some sense more fundamental in comparison to the free group $\mathbb{Z} = \langle a | \rangle$.
Maybe this is too philosophical but in some sense the group $\mathbb{D}$ is more natural than 
$\mathbb{Z}$ as the generator structure is simpler, consisting of involutions. 
The fact that a natural metric on $D_{\infty}$ actually only allows right translation as
isometry makes the idea even more appealing.

\paragraph{}
Much as been written about the nature of space of time. Already when looking at the book literature
we see an obsession with the topic: 
\cite{Eddington1928,Hawking,Rovelli,Mackey92,Bardon,Moxey,Prigogine,Oriti,Reichenbach,Albert2000, 
CarrollTime,HandbookTime, Lineweaver,Davies,ScientificAmericanTime,SmolinTime,GleickTime,Damour,RovelliTime,JanusPoint}.
Much of these writings are also philosophical, cosmological or interpretative. 
But they show what a vast pool of ideas have been proposed. Some examples:

\begin{center} \begin{tabular}{ll}
Zeno        &  Motion needs to be explained properly \\
Laplace     &  The future is determined from the present \\
Newton      &  Time is an absolute and given by a real time axes  \\
Eddington   &  Entropy drives time, makes processes irreversible \\
Einstein    &  Time is part of a space-time geometry  \\
Prigogine   &  Determinism is flawed on a fundamental level \\
Mueller     &  Expansion of the universe creates more time  \\
Smolin      &  Time is an illusion. Physics evolves with time   \\
Rovelli     &  Time emerges in a statistical context \\
Barbour     &  Time evolves towards rich complexity  
\end{tabular} \end{center}

\paragraph{}
As pointed out in \cite{JanusPoint}, the problem why we can not reverse time is not 
compatible with the fundamental laws of physics for the simple reason that the later laws
are all reversible. Laws in statistical mechanics are not really fundamental 
as they emerge as an idealized limit of a {\bf microcanonical description} that is a reversible
n-body problem. All fundamental n-body problems are deterministic. Indeterministic
phenomena like collisions in the Newtonian n-body problem or simultaneous multi-particle
collision singularities for billiards are just limitations of the model. Gas particles are
not billiard models and celestial bodies have a radius. As for the evolution of space-time
singularities we have even no model which tells for example if a black hole evaporates to 
the point that it loses the property of being a singularity. This is a question of t'Hooft. 
It illustrates very well, what difficulties we have when combining the quantum and the 
gravitational world. 

\paragraph{}
The problem of the arrow of time is neither a cosmological nor a statistical 
mechanics problem because we can look at the problem also in a very small
laboratory and study it involving only very few particles. There is no fundamental force which allows
us to change the direction of time. While we can slow down time by placing an apparatus
into a gravitational field or observe the experiment in motion, 
we can not reverse it without some exotic matter. Not that imagination lacks: there are many 
stories and movies dealing with time machines and time travel. 

\paragraph{}
The puzzle is why we have no ability to speed up or slow down time willfully for
a part of physical space, while keeping time the same in the rest.
This is especially enigmatic if one looks at mathematical
models of space-time which treat space and time on the same footing as geometry. 
We can for example use electric or magnetic field to reverse the velocity of a moving
particle but we can not use force fields to reverse the velocity of the time evolution
for the particle. The inability to reverse time is something we can experience 
in small scale laboratories already. 
The loss of information due to motion happens in any direction, whether we move 
forward or backwards. 

\paragraph{}
Eddington explained the arrow of time with ``entropy". One must note
however that a mathematical definition of entropy of the ``universe"
is beyond mathematical rigor. In a cosmological setting it is questionable 
as the universe is not a system in equilibrium. 
Different parts of space which have no causal connection. Mathematically, entropy 
has been defined by Bolzmann or Shannon for finite probability spaces like
finite $\sigma$-algebras of a probability space. The notion 
can also make sense for smooth enough probability distributions as 
{\bf differential entropy} $-\int f(x) \log(f(x)) \; dx$
in probability theory, but entropy is already not defined if $f$ is not smooth enough.
Entropy can become negative for continuous distributions.
Most notions of entropy assume the entropy of some sort of finite $\sigma$-algebra 
$\mathcal{A} = \{ A_1, \dots, A_n\}$ with $P[A_i]=p_i$ or quantum mechanically as von Neumann entropy
$-{\rm tr}(\rho \log(\rho))$, where $\rho$ is a finite dimensional {\bf density matrix}. 
The mathematics has been pushed to dynamical systems \cite{RuelleThermo} using
dynamically defined entropy in an ergodic setting (metric entropy of Sinai)
or topological setting (topological entropy). 

\paragraph{}
It is an empirical fact that any fundamental dynamical system we know is reversible.
Any dynamical system defined by fundamental process therefore {\bf  preserves entropy}.
All fundamental particle or wave motions are {\bf Hamiltonian 
systems} and so reversible, whether they are classical, relativistic or quantum. If one has an invariant
probability measure, one has Poincar\'e-recurrence. A good model problem is the {\bf Vlasov
gas} in a finite container. This is an infinite dimensional Hamiltonian system
of the form $f''(x) = \int \nabla V (f(x)-f(y)) \; dy$. As a Hamiltonian system, it 
is obviously a reversible system and generalizes the {\bf $n$-body problem} with potential $V$. 
The $n$-body problem is the case when the measure
$\mu=dy$ is replaced with a measure supported by finitely many points. What happens is that if we put such
a gas into a container like a box, then initially this macroscopic system will move 
(the gas shakes the container leading to a non-trivial coupled system of a infinite and finite dimensional
Hamiltonian system). The system however will settle down asymptotically (at least this is what we expect. 
No mathematical proof has been done even so there is no doubt about this, especially if the billiard dynamics
in the container is ergodic). The reason is that time evolution $f(t)$ in a
space of smooth functions moves towards parts of space where $f$ becomes more and more
complicated. This happens even in integrable situations like a piston in a
cylindrical container on which we have a Vlasov gas on both sides. We observe an arrow of
time, even-so the system is Hamiltonian. For more details see \cite{Vlasovgasboundary}
or the Chapter on Vlasov in \cite{knillprobability}. For a Vlasov gas, we can look at the
particle density which is $\rho(x) = \int_Y f(x,y) \; dy$. This is a conditional expectation
and of course, the entropy of $\rho_t(x)$ in time will increase to the maximum where the 
density is constant. The reason for the increase of entropy is that we take a conditional
expectation. The particle {\bf position density} $\rho(x)$ contains much less information than the 
particle {\bf phase space density}. Looking at $\rho$  is looking at the infamous shadows of
Plato's cave but the allegory is different as we {\bf know} what the phase space system is. It is 
perfectly nice infinite dimensional Hamiltonian system. 

\paragraph{}
Despite all the statistical mechanical or cosmological diversions,
the conundrum remains why we we can freely move around in 
space, but not in time. The statistical mechanical point of view is obviously 
based on the fact that smooth quantities evolving under simple evolutions become
more and more complicated. Mathematically this can be expressed that their derivatives grow over 
time. This leads to {\bf loss of information} as we have to go to finer and finer
$\sigma$-algebras to keep the initially known information. This mechanism is often dubbed ``entropy increases".
But fundamentally, it is just our inability to keep track. This losing track happens
very fast if the system shows {\bf sensitive dependence on initial conditions}. Iterating
a map $T(x,y)=(2x-y + 2 \sin(x),x) \; {\rm mod} \; 1$ a 100
times forward and then 100 times backwards does not bring
us to the same point if the calculation is done in real arithmetic.
Time is obviously built in a completely different way than space 
even so general relativity treats space time as a {\bf space-time manifold}.
Let us mention, besides {\bf Vlasov dynamics}, an other Hamiltonian system which features
some sort of arrow in time. It deals with a Hamiltonian system which features
an expansion of space.

\begin{comment}
T[{x_, y_}] := Mod[{2 x - y + 2 Sin[x], x}, 1];
S[{x_, y_}] := Mod[{y, -x + 2 y + 2 Sin[y]}, 1];
X = Last[NestList[T, {0.3, 0.4}, 10]]
Y = Last[NestList[S, X, 10]]
\end{comment}

\paragraph{}
A spontaneous expansion of a geometry like a Riemannian manifold or a graph can be explained by 
moving in the symmetry group of the geometry: if we allow the Dirac operator $D=d+d^*$ 
of a Riemannian manifold or simplicial complex deform isospectrally, this leads to an 
expansion if one postulates distances is coming from the Dirac operator $d+d^*$.
The deformed operator $D(t)= d + d^* + b$ will develop
a diagonal part but still feature exterior derivatives $d(t)$ which keeps the
same cohomology. As $d$ defines distance, the deformation affects distance. 
The origin of the inflation (the Janus point in the terminology of 
Barbour) is derived mathematically as the time in which the Dirac operator $D$
has no diagonal part. It does not matter in what direction the system moves. 
It always expands when moving away from the Janus point. If we pick the time $0$
when the Dirac operator has no diagonal parts, we can evolve in any direction
and experience expansions in any direction \cite{IsospectralDirac2}. Why does the
Dirac operator move at all? Because not moving in the isospectral set has probability
zero (similarly that there is zero probability that a stone moving in space does not rotate
or that any system with continuum symmetry is in a certain particular state.
Still, also this Dirac expansion model 
postulates that time is the real axes. Time could be multi-dimensional
like in a Kaluza-Klein situations. Also, having time as a discrete non-commutative
structure similar as in {\bf Connes-Lott model} of the standard model is not
excluded. Going to a dihedral time is a non-commutative notion of time which actually 
is very close to a discrete time $\mathbb{Z}$. 

\paragraph{}
In physics, one has seen again and again that the role of {\bf symmetries} is 
extremely important. By {\bf Noether's theorem}, symmetries are related to {\bf invariants}.
Translation symmetry is related to momentum conservation, time symmetry related to 
energy conservation, rotational symmetry to angular momentum conservation etc. 
Also, symmetries seem to dictate what processes are possible and what processes
are not. We also have to accept the empirical fact
that all {\bf fundamental processes} are reversible. Postulates like Prigogine's that
fundamentally things are not reversible is not supported yet by any fundamental model: 
as pointed out by Barbour for example, 
all fundamental models in classical or relativistic mechanics, quantum mechanics, 
general relativity, electrodynamics, quantum field theories and especially 
the standard model allow for a time reversal.
In fundamental interaction, one knows also the {\bf TCP symmetry} but no fundamental 
process is known which does not work backwards in time (possibly changing charge 
and parity). What happens with the {\bf dihedral time hypothesis} is that we have
a time structure which is postulated to be the infinite dihedral group. The
justification so far to look at this postulate is purely mathematical and not based
on physics. The dihedral group is a natural metric space in which right group translations
is isometries, but where time reversal, nor left translation nor the inverse operation
are isometries. There is an arrow of time. 

\paragraph{}
The {\bf dihedral time hypothesis} is the suggestion that time reversal is not possible,
because doing so would {\bf break of symmetry}. 
The use of a {\bf non-abelian time in the discrete},
changes from an unnatural space $\mathbb{Z}$ to a natural space $D_{\infty}$.
There are two steps which are related: The postulate of a {\bf quantized time} 
replacing the Lie group $\mathbb{R}$ with the discrete group $\mathbb{Z}$. 
But this loses naturalness. In order for it to become natural, it requires to use 
$\mathbb{D}$ and having {\bf time to become non-commutative}. 
This breaks the time reversal symmetry leading to an {\bf arrow of time}. In the 
group $\mathbb{D}=\{ a^2=b^2=1 \}$, going from $x$ to $x^{-1}$ is not an isometry,
if the group comes from a natural metric (forcing the algebraic structure uniquely from
the metric).  This is something I learned in 2026 from Jakub: 
Time reversal symmetry $x \to x^{-1}$ is not an isometry with 
respect to a natural metric. 

\section{About this document}

\paragraph{}
In June of 2008 we wondered for which metric groups, the topology determines the
arithmetic. It was a bit of a surprise to see hat the group $\mathbb{Z}$ of rational
integers does not qualify in this respect as it also carries also a non-abelian
$\mathbb{D}$ group structure.
However, the non-commutative dihedral group $\mathbb{D}$ is natural in this respect.
There is a metric on $\mathbb{D}$ as a set such that the group structure $\mathbb{D}$
is determined. The general question then appeared which groups are natural in this
respect. Finite groups of prime order are natural simply because there is only one
group structure. Already the cyclic groups $\mathbb{Z}_n$ are not natural in general
for $n \geq 4$ while all non-commutative dihedral groups $\mathbb{D}_n$ are natural.
The real line $\mathbb{R}$ is natural. There is only one group operation on the usual
Euclidean metric $d(x,y)=|x-y|$ such that all right translations are isometries
isometries because the group of isometries of the line is $\mathbb{R} \rtimes Z_2$
generated by translations and reflections.

\paragraph{}
The additive group $(\mathbb{Z},+)$ is the foundation 
for virtually all the arithmetic we know. Since the beginning of mathematics  
in the Paleolithic era, illustrated by artefacts like the Ishango bone found near the Semliki River,
the addition of numbers and especially the commutativity of the basic additive arithmetic has been
developed all over the world. The pair of arithmetic operations, addition and multiplication,
defines a ring, bound together by distributivity. The additive operation on integers as the most fundamental
operation has never been questioned and Kronecker was the first who formulated ultra finite ideas. 
The quotation of Kronecker has been popularized in \cite{Bell} and the book title of the 
collection \cite{HawkingGodCreatedIntegers}. 

\paragraph{}
Could it be that some non-commutative additive structure is more fundamental than the 
known additive structure on integers? I could argue here that the answer is ``yes", 
despite the fact that historically, serious steps into non-commutative mathematics \cite{Connes} 
have emerged much, much later
than the counting on integers. Non-commutative multiplicative structures most likely first appeared as
{\bf symmetry groups of geometries} and especially in permutation groups. The symmetry group of a 
polygon is already non-commutative. 

\paragraph{}
But much earlier, already Euclid must have known (at least intuitively as Euclid had no group
theory at hand) that the affine group of the plane which contains rotations and translations, 
is not commutative. It is generated by all reflections on lines as the composition of two 
reflections at parallel lines produce a translation perpendicular to the lines and the 
composition of two reflections at intersecting lines produces a rotation around the 
intersection point. In both cases, the order in that Coxeter 
group picture matters. There have been attempts to use reflection geometry early in geometry classes
\cite{Bachmann1959,JegerTransformationGeometry} 

\paragraph{}
When {\bf matrix algebra} and later also {\bf representation theory} emerged, one has linked groups 
with operator algebras  Non-commutativity in the multiplication turned out to be fundamental.
In the very small, the multiplication of observables in quantum theories is no more commutative. 
This is illustrated by the commutation relations in quantum mechanics which is now part of 
pop culture and motivated fundamental developments like non-commutative differential 
geometry \cite{Connes} by pushing ideas from non-commutative measure theory and commutative
topology to differential geometry. Non-commutative measure theory is an extended theory of von-Neumann
algebras and non-commutative topology is an extended theory of $C^*$ algebras. Non-commutative
differential geometry uses {\bf spectral triples} $(D,A,H)$, in the form of a Hilbert space $H$, 
an algebra of operators and a self-adjoint operator $D$ satisfying for which $[a,D]$ has finite
norm for all $a \in A$. The operator $D$ allows a polar decomposition $D=|D| U$ which then defines a 
metric space $(X,d)$ where $X$ is the set of pure states on the norm closure of $A$.

\paragraph{}
The question of a more fundamental arithmetic than the``cave" arithmetic $(\mathbb{Z},+)$
of course depends on what is considered ``natural" or ``fundamental". 
We do not mean to weaken or change the axioms but insist on keeping groups and 
metric spaces. Using only a {\bf monoid structures} for example
and to look for the natural numbers $(\mathbb{N},+)$, would not be an option for us.% be an option
We want to have a criterion which makes a {\bf group} natural and try to achieve
this by linking the group structure with a {\bf different category of objects} 
outside of arithmetic. This can be an {\bf order structure}, a {\bf topological structure}, a 
{\bf measure theoretical structure}, a {\bf simplicial complex structure},
a {\bf differential structure} of a {\bf symplectic structure}. 
In {\bf Lie group theory} for example, one asks the group structure to be 
compatible with a given {\bf differential structure}
on the manifold: the addition and inversions need to be a smooth operation. 

\paragraph{}
What we start to explore here could be done between any two algebraic and topological
categories. It turns out however that things often do not lock together nicely. 
Given a too large algebraic structure like the category of 
{\bf groupoids}, it does not force a link in the sense that the algebraic structure determine the
topological structure or that the topological structure decides the algebraic structure. 
Also taking something more general than a metric space like {\bf topological spaces} is harder:
the reason is that the {\bf symmetry group of a topological space}, the 
{\bf group of homeomorphisms} of a topological space is in general 
much larger than the space itself. The group of homeomorphisms on a circle $\mathbb{T}^1$ 
for example is infinite-dimensional, while the circle itself is finite-dimensional. 
But the group of transformations which are compatible with the metric, the {\bf isometries} of
the metric space have a chance to be a finite-dimensional space. The group of isometries 
on any {\bf compact Riemannian manifolds} for example is known to be a Lie group 
(possibly discrete or even trivial) and so finite-dimensional. 

\paragraph{}
This file, like many other projects has not seen the light for more than a decade. 
In the summer of 2021, after a minor health scare,
I became aware that I have reached an age, where for pure statistical reasons 
the chances of termination increase 
and I started to work already then on this file after realizing the odds. 
\footnote{In the US, the death rate of men per year in my age group is about 
$p=1/100$ which is 10 times higher than for a 20 year old. (https://www.ssa.gov)} 
It would be a nightmare to observe from a possible after-life 
(probably working as a junior data research assistant to the Almighty who tries 
desperately to decrease the odds that the ``humans on earth experiment" 
would end up in failure), that rather interesting files like this one would have 
been overwritten by a teenager playing GTA 5 on a recycled SSD drive of my 
workstation, that has ended up in a garage sale. I post this on the ArXiv. 
There are no plans to send this document for publication. 
It definitely would need a complete rewrite, starting entirely from scratch. 
Any attempts in the last months to shorten the paper made things only longer.
Any comments or feedback are welcome, especially since I have not seen the 
theme appear anywhere else. 

\pagebreak 

\section*{Appendix}

\paragraph{}
In the 2022 version of this paper, we assumed also that left translation and the inverse
operation are isometries. {\bf Jakub Gismatullin} alerted me in July 2026 that this original definition 
is not consistent with some of the results (and so provided an other example illustrating the
"proofs and refutation" picture of Lakatos \cite{Lakatos}).
In \cite{Gismatullin}, having both left and right translation invariance was referred as {\bf strongly 
natural}. Gismatullin shows for example that the direct product $G_1 \times G_2$ of two strongly 
natural groups is strongly natural and that a finite abelian group is strongly natural if and only 
if it is not of the form $Z_4 \times Z_2^n$ for $n \geq 0$. Naturality under the revised
right-translation definition is more common. It turned many results around. 

\paragraph{}
The paper \cite{LeemannDeLaSalle2022} illustrates that 
one can check many groups to be natural by studying their graphical regular respresentation. It covers
$F_n$ for $n \geq 2$, lamplighter groups, $S_n$ for $n \geq 4$ and $A_n$ for $n \geq 5$
and $D_n$ for $n \geq 5$, it even covers every finitely generated non-abelian simple group. 
The GRR exceptions in the finite non-abelian groups are $A_4,Q_8 \times C_3, Q_8 \times C_4$ and $B(2,3)$.
But we have seen that failure to have a graphical regular representation does not not imply 
non-naturality. The cases $Q_8,A_4$ or $B(2,3)$ illustrate this. But $Q_8$ is non-natural. 

\paragraph{}
The results for natural and strongly natural are sometimes complementary. We wanted to have that a metric
space defines a unique group structure up to isomorphism. Already for the 
prototype situation $\mathbb{Z}$ or $\mathbb{D}_{\infty}$ the results are different: 
the infinite dihedral group is natural but not strongly natural, the integers are not 
natural but strongly natural. 

\begin{center} 
\begin{tabular}{|c|c|c|} \hline
   Group                &   Natural?   & Strongly natural?  \\  \hline
  $\mathbb{Z}$          &    NO        &   YES              \\ 
  $\mathbb{D_{\infty}}$ &    YES       &   NO               \\ \hline
\end{tabular} \end{center} 

\paragraph{}
The weaker condition also does not allow the 
semi-direct product and so the direct product of finite natural groups to remain natural.
We had claimed that Zig-Zag products of suitable Cayley graphs would do the job.
So, also the wreath products do not preserve natural groups in general. Here is an example:

\begin{center}
\begin{tabular}{|c|c|c|} \hline
   Group  (m is odd)   &   Natural?   & Strongly natural?  \\  \hline
  $Z_{2m}=Z_m \times Z_2$    &    NO for $m>2$   &   YES for $m > 2$    \\
  $D_m = Z_m \rtimes Z_2$    &    YES            &   NO               \\ \hline
\end{tabular} \end{center} 

\paragraph{}
Lets look at the example of $Z_6 = Z_2 \times Z_3$. It is not natural. 
Take any right-invariant metric $d$ on $Z_6$.
Because $Z_6$ is abelian, the metric is also left-invariant so that 
inversion $x \to -x$ is an isometry (as seen in the proof below). 
Therefore, also $a(x)=1-x$ is an isometry. Together with translation $t(x)=x+2$
These generate $\langle t,a \rangle$, which is $D_3=S_3$. Every
right-invariant metric for $Z_6$ also supports $D_3$.
A similar argument shows that $Z_{2m}$ for $m \geq 2$ is non-natural as it also 
admits the dihedral group $D_m$. For $m>2$ however the group $D_m = Z_m \rtimes Z_2$ 
is natural as when we break the symmetry in geodesic metric of the
zig-zag product of the Cayley graphs, we break the symmetry and $Z_{2m}$ does no more
work. 

\paragraph{}
Including left translation invariance or including the inverse operation to be
isometries are equivalent. Here is the statement: 

\begin{lemma}
If a metric group $(G,d)$ satisfies  
$d(y,z) = d(yx,zx)$ for all $x,y,z$, then the following two statements are equivalent: \\
a)  $d(y,z) = d(x y,x z)$ for all $x,y,z$. \\
b)  $d(y,z) = d(y^{-1},z^{-1})$ for all $y,z$. 
\end{lemma}

The rather pedantic nature of the proof was an initial reactions of the 
shock of having learned that my original paper of 2022 had used definitions 
that were over-stated.  Lakatos \cite{Lakatos} has considered this a 
genuine part of the epistemology of science. It should not happen however. 

\begin{proof}
{\bf $a \Rightarrow b$:}
\begin{eqnarray*}
                  d(x,y) &=& d(1,x^{-1}  y)      \; {\rm  (left \; invariance})  \\
                         &=& d(y^{-1},x^{-1})    \; {\rm  (right \; invariance}) \\
                         &=& d(x^{-1},y^{-1})    \; {\rm symmetry \; of \;  metric}) \\
\end{eqnarray*}

{\bf $b \Rightarrow a)$:} \\

Part (i)  $d(x,1) = d(x^{-1},1)$.

\begin{eqnarray*}
       d(x,1) &=& d(1,x^{-1})  \;   {\rm right \;  invariance} \\
              &=& d(x^{-1},1)  \;   {\rm symmetry \; of \;  metric}  \\
\end{eqnarray*}

(ii) $d(x y,1) = d(y x,1)$

\begin{eqnarray*}
  d(x y,1) &=& d(y^{-1} x^{-1},1)  \;  {\rm by \; (i)} \\
           &=& d(y^{-1},x)         \;  {\rm by \; right \; invariance} \\
           &=& d(y,x^{-1})         \;  {\rm by \; inverse \; invariance \;  assumption} \\
           &=& d(y x,1 )           \;  {\rm by \; right \; invariance} \\
\end{eqnarray*}

(iii) Claim:

Proof:
\begin{eqnarray*}
   d(x y,x z) &=& d(x y (x z)^{-1},1)     \; {\rm right \; invariance} \\
              &=& d((y (x z)^{-1}) x, 1)  \; {\rm part (ii)}  \\
              &=& d(y z^{-1}, 1)  \; {\rm associativity} \\
              &=& d(y, z))         \; {\rm right invariance} \\
\end{eqnarray*}
\end{proof}

\paragraph{}
Lets remark about the difference between naturality and iso-group (\cite{Niemiec2014}). 
Niemiec's iso-group notion is a notion for topological transformation groups.
On the other hand, naturality is an abstract-group/metric notion.
Niemiec calls a group $G$ of homeomorphisms of a locally compact
Polish space $X$ an {\bf iso-group} of transformations if 
there is a metric $d$ on $X$ such that $G = {\rm Iso}(X,d)$. This is different
from naturality as we ask ask for a metric $(G,d)$ {\bf on the space $G$ itself!}
that possesses an up to isomorphism unique group structure $(G,+)$
for which all right group translations are isometries.

\paragraph{}
Are iso-groups automatically natural? The answer is no:
Niemiec shows that essentially every locally compact
Polish group can be realized as a full isometry group of 
some suitable locally compact Polish space.
We know however that there are locally compact Polish groups such like
$\mathbb{Z}$ and $\mathbb{R} \times C_2$ that are not natural.

\paragraph{} 
There are no example among locally compact Polish  natural groups
that are not an iso-group. It is natural to ask therefore: 
Are natural groups automatically iso-groups? Also here the answer is no:
the p-adic integers $\mathbb{Z}_p$ with $p>2$ are natural group for which the
naturalizing metric cannot make the usual metric is not the natural metric.

\paragraph{}
Are there groups that both natural and iso-groups? Yes.
There are many examples, like $\mathbb{R}$.
The infinite dihedral group $D_{\infty}$ is a locally compact 
Polish group and therefore is an iso-group. 

\paragraph{}
Is there are group that is neither natural nor an iso-group: Yes. 
$(\mathbb{Z},+)$ equipped with the p-adic topology is not natural.
Equipped with the $p$-adic metric, it is not complete and its completion
is $\mathbb{Z}_p$. Niemiec's general characterization says that a topological group $G$ is
topologically isomorphic to the full isometry group of a complete metric
space precisely when it coincides with its $G_{\delta}$-closure in its
Raikov completion. [A topological group G has two canonical uniformities: 
the left uniformity, determined by $x^{-1} y$ being close to $1$
and the right uniformity, determined by $x y^{-1}$ being close to $1$.
The Raikov or two-sided uniformity is the supremum of the left and right uniformities.
$\overline{G}$ is then defined as the completion of G with respect to this two-sided uniformity.
Now multiplication and inversion extend to the completion, so that $G$ is again a topological group.
If $G$ is abelian, left and right uniformities coincide and two-sided completion 
is the ordinary uniform completion.] 

\paragraph{}
In the context of {\bf graphical regular representations}, one is interested in
the following question: 
How much of the algebraic structure an unlabelled Cayley graph can remember?
The {\bf naturality question} is different, as it asks 
How much of the group law can a metric space remember?
This question is often related to the previous one but it is different and also is interesting
in the context of Lie groups, where we have seen that every connected Lie group
is natural and geometry determines algebra. 

\paragraph{}
And finally an update rather philosophical remarks about the role of time.
The fact that time $D_{\infty}$ comes from a natural metric space and that we
have an {\bf arrow of time} fits now even better after correcting the definition.
We can now take the picture of Emmy Noether more seriously: 
As the left group translation is
not a symmetry on the natural group $D_{\infty}$, this group is a time structure with an
{\bf arrow of time}. 

\paragraph{}
Only going forward in time produces unitary 
operations on observables. Since 2022, we have worked on discrete time wave
dynamics (see for example \cite{ConnectionDirac} and discrete time geodesic flows
\cite{geodesics1,geodesics2}. Discrete time wave dynamics on a discrete space was
forced by causality. If we look at a differential equation $u_{tt} = -Lu$ for some
discrete Laplacian $L$ then a change of the geometry at the end of the universe immediately
influences the dynamics $u(0)$ even for small times because eigenfunctions are non-local.
If we take locality serious, we are forced to look at discrete time, if we have
discrete space. It is no problem to make discrete time dihedral. 
Also the discrete geodesic flow was by nature dihedral. We wrote the dynamics as $T=AB$ with 
two involutions $A,B$  because it is then easier to write then the inverse as $T^{-1} = BA$.

\bibliographystyle{plain}

\end{document}